\documentclass{amsart}

\usepackage{amssymb}
\usepackage{graphicx}

\newtheorem{theorem}{Theorem}[section]
\newtheorem{lemma}[theorem]{Lemma}
\newtheorem{proposition}[theorem]{Proposition}
\newtheorem{corollary}[theorem]{Corollary}

\newtheorem{_conjecture}[theorem]{Conjecture}

\newtheorem{_problem}[theorem]{Problem}

\newtheorem{_algorithm}[theorem]{Algorithm}

\newtheorem{_subroutime}[theorem]{Subroutine}

\newtheorem{_claim}[theorem]{Claim}
\newenvironment{claim}{\begin{_claim}\rm}{\end{_claim}}

\newtheorem{_subclaim}[theorem]{Sub-claim}

\newtheorem{_definition}[theorem]{Definition}
\newenvironment{definition}{\begin{_definition}\rm}{\end{_definition}}

\newtheorem{_question}[theorem]{Question}
\newenvironment{question}{\begin{_question}\rm}{\end{_question}}

\newtheorem{_observation}[theorem]{Observation}
\newenvironment{observation}{\begin{_observation}\rm}{\end{_observation}}

\newtheorem{_remark}[theorem]{\it Remark}
\newenvironment{remark}{\begin{_remark}\rm}{\end{_remark}}

\newtheorem{_example}[theorem]{Example}

\numberwithin{equation}{section}
\numberwithin{table}{section}
\numberwithin{figure}{section}

\newcommand{\A}{\mathord{\mathbb A}}

\newcommand{\F}{\mathord{\mathbb F}}
\renewcommand{\P}{\mathord{\mathbb  P}}
\newcommand{\Q}{\mathord{\mathbb  Q}}
\newcommand{\R}{\mathord{\mathbb R}}
\newcommand{\Z}{\mathord{\mathbb Z}}

\newcommand{\CCC}{\mathord{\mathcal C}}

\newcommand{\EEE}{\mathord{\mathcal E}}

\newcommand{\GGG}{\mathord{\mathcal G}}

\newcommand{\III}{\mathord{\mathcal I}}

\newcommand{\LLL}{\mathord{\mathcal L}}

\newcommand{\OOO}{\mathord{\mathcal O}}
\newcommand{\PPP}{\mathord{\mathcal P}}

\newcommand{\UUU}{\mathord{\mathcal U}}
\newcommand{\VVV}{\mathord{\mathcal V}}

\newcommand{\ZZZ}{\mathord{\mathcal Z}}

\font\mathgot=eufm10

\newcommand{\AAAA}{\mathord{\hbox{\mathgot A}}}
\newcommand{\SSSS}{\mathord{\hbox{\mathgot S}}}
\newcommand{\moduli}{\mathord{\hbox{\mathgot M}}}

\newcommand{\maprightsp}[1]{\; \smash{\mathop{\; \longrightarrow \; }\limits\sp{#1}}\; }

\newcommand{\mapdownsurj}{
\hbox{$\bigm\downarrow$}
\llap{\hbox{\raise 2pt\hbox{$\bigm\downarrow$}}}%
\vstrechmapdown
}

\newcommand{\inj}{\hookrightarrow}

\newcommand{\isom}{\smash{\mathop{\;\to\;}\limits\sp{\sim\,}}}
\newcommand{\isomto}{\smash{\mathop{\;\to\;}\limits\sp{\sim\,}}}

\newcommand{\set}[2]{\{\; {#1} \; \mid \; {#2} \;  \}}

\newcommand{\shortset}[2]{\{ {#1}\mid{#2}\}}

\newcommand{\map}[3]{ #1 \, : \, #2 \, \to \, #3}
\newcommand{\mapisom}[3]{ #1 \, : \, #2 \; \isom \; #3}

\newcommand{\sm}{\setminus}
\newcommand{\st}{\subset}

\newcommand{\sprime}{\sp\prime}

\newcommand{\spar}[1]{\sp{(#1)}}
\newcommand{\spprime}{\sp{\prime\prime}}

\newcommand{\sptimes}{\sp{\times}}

\newcommand{\dual}{\sp{\vee}}

\newcommand{\inv}{\sp{-1}}

\renewcommand{\qed}{\hfill {$\Box$}}

\renewcommand{\Im}{\operatorname{\rm Im}\nolimits}
\newcommand{\Ker}{\operatorname{\rm Ker}\nolimits}

\newcommand{\pr}{\operatorname{\rm pr}\nolimits}
\newcommand{\id}{\operatorname{\rm id}\nolimits}

\newcommand{\Pic}{\operatorname{\rm Pic}\nolimits}
\newcommand{\GL}{\operatorname{\it GL}\nolimits}
\newcommand{\PGL}{\operatorname{\it PGL}\nolimits}

\newcommand{\Aut}{\operatorname{\rm Aut}\nolimits}
\newcommand{\Hom}{\operatorname{\rm Hom}\nolimits}

\newcommand{\Sing}{\operatorname{\rm Sing}\nolimits}
\newcommand{\Spec}{\operatorname{\rm Spec}\nolimits}
\newcommand{\Supp}{\operatorname{\rm Supp}\nolimits}

\newcommand{\Stab}{\operatorname{\rm Stab}\nolimits}

\newcommand{\closure}[1]{\overline{#1}}

\newcommand{\rmand}{\textrm{and}}

\newcommand{\quand}{\quad\rmand\quad}

\newcommand{\DK}{\mathord{\rm DK}}

\newcommand{\Pt}{\P\sp 2}
\newcommand{\card}[1]{|#1|}

\newcommand{\XG}{X\sb G}
\newcommand{\YG}{Y\sb G}

\newcommand{\Pow}{\operatorname{\rm Pow}\nolimits}

\newcommand{\NS}{\mathord{{\it NS}}}

\newcommand{\Bs}{\mathord{\rm Bs}}

\newcommand{\lb}{[}
\newcommand{\rb}{]}

\newcommand{\oo}{\omega}
\newcommand{\bo}{\bar\omega}

\newcommand{\va}{\mathord{\bf a}}
\newcommand{\vx}{\mathord{\bf x}}

\newcommand{\Hesse}{\mathord{\rm Hesse}}
\newcommand{\GH}{G\sb{\Hesse}}

\newcommand{\pol}{\LLL}
\newcommand{\Xpol}{(X, \LLL)}
\newcommand{\ZXpol}{Z\sb{\Xpol}}

\newcommand{\typesharp}{(\sharp)}

\newcommand{\Phipol}{\Phi\sb{|\pol|}}

\newcommand{\gcode}{\CCC}
\newcommand{\abscode}{\mathord{\bf C}}
\newcommand{\CA}{\abscode\sb A}
\newcommand{\CB}{\abscode \sb B}
\newcommand{\CC}{\abscode\sb C}
\newcommand{\CT}{\abscode\sb T}

\newcommand{\Ps}{\PPP}
\newcommand{\gs}{\GGG}

\newcommand{\gcodeXpol}{\gcode(X,\pol,\gamma)}

\newcommand{\Ff}{\F\sb 4}

\newcommand{\CDK}{\abscode\sb{\DK}}

\newcommand{\Psb}[1]{P\sb{#1}}

\newcommand{\PG}{\mathord{PG}}
\newcommand{\AF}{\mathord{AF}}

\newcommand{\cremona}{\mathord{\hbox{\rm{CT}}}}
\newcommand{\Cremona}{\widetilde{\cremona}}
\newcommand{\mcremona}{\mathord{\hbox{\rm{ct}}}}
\newcommand{\ZZ}[1]{Z(#1)}

\newcommand{\spT}{{}\sp t\hskip -1pt}

\newcommand{\GA}{\mathord{G\hskip -.2pt A}}
\newcommand{\GB}{\mathord{G\hskip -.2pt B}}
\newcommand{\GC}{\mathord{G\hskip -.1pt C}}
\newcommand{\GT}{\mathord{G\hskip -.1pt T}}

\newcommand{\wtil}[1]{\widetilde{#1}}

\newcommand{\blu}{\mathord{\beta}}
\newcommand{\tildeblu}{\mathord{\tilde\beta}}
\newcommand{\bld}{\mathord{\beta\sprime}}
\newcommand{\ED}{M}
\newcommand{\EED}{N}
\newcommand{\Line}{H}
\newcommand{\HL}{\sqrt{L(H)}}

\newcommand{\JT}{J\sb{T\sprime}}
\newcommand{\JS}{J\sb{T}}
\newcommand{\trivial}{\Delta=0}

\newcommand{\ideal}{\mathord{\hskip.5pt \III\hskip.5pt}}
\newcommand{\intsct}{}

\newcommand{\JA}{J\sb{A}}
\newcommand{\JB}{J\sb{B}}
\newcommand{\JC}{J\sb{C}}
\newcommand{\JAsprime}{J\sb{A}\sprime}
\newcommand{\JBsprime}{J\sb{B}\sprime}
\newcommand{\JCsprime}{J\sb{C}\sprime}

\newcommand{\coreq}[4]{ #1 & #2 & #3 &\parbox[t]{10cm}{$#4$\phantom{\vrule depth 5pt }}  \\ \hline }

\begin{document}

\title[Moduli curves of supersingular $K3$ surfaces]{%
Moduli curves of supersingular $K3$ surfaces\\ 
in  characteristic $2$ 
with Artin invariant $2$}

\author{Ichiro Shimada}
\address{
Department of Mathematics,
Faculty of Science,
Hokkaido University,
Sapporo 060-0810,
JAPAN
}
\email{shimada@math.sci.hokudai.ac.jp
}


\dedicatory{Dedicated to Professor Igor~V.~Dolgachev
 for his 60th birthday}

\begin{abstract}
We construct explicitly  moduli curves of polarized supersingular $K3$ surfaces 
in characteristic $2$ with Artin invariant $2$.
As an application,
we detect a ``jump" phenomenon in a family of automorphism groups of supersingular $K3$ surfaces
with a constant N\'eron-Severi lattice.
\end{abstract}

\maketitle
\section{Introduction}\label{sec:intro}
A $K3$ surface is called \emph{supersingular}
if  its  numerical N\'eron-Severi lattice   is of   rank $22$.
Supersingular $K3$ surfaces exist only in positive characteristics.
Artin showed in~\cite{Artin74} that,
in characteristic $p>0$,
the discriminant of 
the numerical N\'eron-Severi lattice of a supersingular $K3$ surface $X$  is of the form $-p^{2\sigma(X)}$,
where $\sigma (X)$ is a positive integer $\le 10$.
This integer $\sigma (X)$ is called the \emph{Artin invariant} of $X$.
\par
\medskip
We work over an algebraically closed field $k$
of characteristic $2$.
\begin{definition}
Let $X$ be a supersingular $K3$ surface,
and let $\pol$ be a line bundle on $X$ with $\pol\sp 2=2$.
We say that $\pol$ is a \emph{polarization of type~$\typesharp$}
if the following conditions are satisfied:
\begin{itemize}
\item  the complete linear system $|\pol|$ has no  fixed components, and
\item the set of curves contracted by  the morphism $\Phipol : X\to\Pt$ 
defined by $|\pol|$  consists of  $21$ disjoint $(-2)$-curves.
\end{itemize}
\end{definition}
In~\cite{Shimada2003},
we have shown that every supersingular $K3$ surface $X$ in characteristic $2$ has a polarization of type~$\typesharp$,
and that, if $\pol$ is a  polarization of type~$\typesharp$ on $X$,
then the morphism $\Phipol$ is purely inseparable.
In~\cite{Shimada2004},
we have constructed a $9$-dimensional moduli space $\moduli$ of polarized supersingular $K3$ surfaces of type~$\typesharp$.
In this paper,
we investigate the locus $\moduli\sb{2}$ of $\moduli$ corresponding to
supersingular $K3$ surfaces with Artin invariant $2$.
As Artin~\cite{Artin74} showed,
this locus is of dimension $1$.
We will show that the curve  $\moduli\sb{2}$ is a disjoint union of three affine lines
punctured at the origin.
We will also
construct explicitly the universal family of
polarized supersingular $K3$ surfaces over certain finite covers 
of these punctured affine lines.
The construction involves investigations
of   configurations of lines and   conics 
on the  projective plane in characteristic $2$.
These configurations are encoded by certain binary codes.
In order to construct the moduli curve, 
we have to determine the automorphism groups 
of these codes.
The automorphism group of the polarized $K3$ surface is also obtained from 
the  automorphism group of  the corresponding code.
\par
\medskip
Let us briefly review the construction of the moduli space $\moduli$ in~\cite{Shimada2004}.
For a non-zero homogeneous polynomial $G \in H\sp 0 (\Pt, \OOO\sb{\Pt} (6))$
of degree $6$,
we denote by
$$
\map{\pi\sb G}{Y\sb G}{\Pt}
$$
the purely inseparable double cover of $\Pt$ defined by
$W^2=G(X,Y,Z)$.
\begin{definition}
Let $\UUU$ denote the locus of all non-zero homogeneous polynomials $G \in H\sp 0 (\Pt, \OOO\sb{\Pt} (6))$ 
such that   the surface $Y\sb{G}$ 
has  
$21$ ordinary nodes as its only singularities.
\end{definition}
The locus $\UUU$
is Zariski open dense in $H\sp 0 (\Pt, \OOO\sb{\Pt} (6))$.
Indeed, in characteristic $2$,
the differential  $dG$ of $G$ can be defined  as a global section of $\Omega\sp 1\sb{\Pt} (6)$
for any homogeneous polynomial $G\in H\sp 0 (\Pt, \OOO\sb{\Pt} (6))$,
because,
by the isomorphism $\OOO\sb{\Pt}(6)\cong \OOO\sb{\Pt} (3) \sp{\otimes 2}$,
we can assume that the transition functions of the line bundle corresponding
to $\OOO\sb{\Pt} (6)$ are all squares.
Since $c\sb 2 (\Omega\sp 1\sb{\Pt} (6))=21$,
the subscheme $\ZZ{dG}$ defined by $dG=0$
 is reduced of dimension $0$ if and only if  it consists of $21$ points.
 The singular locus  $\Sing (Y\sb{G})$ of $Y\sb G$ is equal to $\pi\sb {G}\inv (\ZZ{dG})$,
 and 
the singular point of $Y\sb G$ lying over a reduced point of $\ZZ{dG}$ is an ordinary node.
Hence 
the condition that $G$ be  a point of  $ \UUU$ is equivalent to the open condition 
that $\ZZ{dG}$ be  reduced of dimension $0$.

Let $\Xpol$ be a polarized supersingular $K3$ surface of type~$\typesharp$.
Then there exists a homogeneous polynomial
$G\in \UUU$
 such that the Stein factorization of $\Phipol$ is written as
$$
X\;\maprightsp{\rho \sb G}\; \YG\;\maprightsp{\pi\sb G}\;\Pt.
$$
Conversely, suppose that we are given  $G\in \UUU$.
Let $\rho \sb G :\XG\to \YG$
be the minimal resolution of the surface $\YG$.
Then $\XG$ is a supersingular $K3$ surface,
and the invertible sheaf 
$$
\pol\sb G:=(\pi\sb G\circ \rho\sb G)\sp *\OOO\sb{\Pt}(1)
$$
on $X\sb G$ is a polarization of type~$\typesharp$.

We put 
$$
\VVV:=H\sp 0 (\Pt, \OOO\sb{\Pt}(3)).
$$
Because we have $d(G+H^2)=dG$ for any $H\in \VVV$,
the additive  group  $\VVV$ acts on the space $\UUU$ by
$$
(G, H)\in \UUU\times\VVV\;\;\mapsto\;\; G+H^2 \in \UUU.
$$
\begin{proposition}\label{prop:equivG}
Let $G$ and $G\sprime$ be homogeneous  polynomials in $\UUU$.
Then the following  conditions are equivalent:
\begin{itemize}
\item[{\rm (i)}]
$Y\sb G$ and $Y\sb{G\sprime}$ are isomorphic over $\Pt$,
\item[{\rm (ii)}]
$\ZZ{dG}=\ZZ{dG\sprime}$, and 
\item[{\rm (iii)}]
there exist $c\in k\sptimes$ and $H\in \VVV$ such that $G\sprime =c\,G +H^2$.
\end{itemize}
\end{proposition}
See \S\ref{sec:review} for the proof.
\par
\medskip
Therefore the moduli space $\moduli$ of polarized supersingular $K3$ surfaces of type $\typesharp$ is constructed by
$$
\moduli =\PGL(3, k)\backslash  \P\sb * (\UUU/\VVV).
$$
For $G\in \UUU$,
let $[G]$ denote the point of $\moduli$ corresponding to $G$,
which corresponds to the isomorphism class of 
the polarized supersingular $K3$ surface $(X\sb G, \pol\sb G)$
of type $(\sharp)$.
By Proposition~\ref{prop:equivG},
the automorphism group $\Aut (X\sb G, \pol\sb G)$ of the polarized supersingular $K3$ surface is canonically 
identified with 
$$
\set{g\in \PGL (3, k)}{ g (\ZZ {dG})=\ZZ{ dG}}.
$$
The moduli space $\moduli$ is stratified by the Artin invariant $\sigma (\XG)$
of  $\XG$.
We put
$$
\moduli\sb{\sigma}:=\set{[G]\in \moduli}{\sigma (X\sb G ) = \sigma}
\quand
\moduli\sb{\le\sigma}:=\set{[G]\in \moduli}{\sigma (X\sb G ) \le \sigma}.
$$
\par
\medskip
As was shown in~\cite{Shimada2004},
the locus $\moduli\sb{\le 1}=\moduli\sb{1}$
consists of a single point $[G\sb{\DK}]$,
where
$$
G\sb{\DK}:= XYZ(X^3+Y^3+Z^3)
$$
is the homogeneous polynomial discovered by Dolgachev and Kondo in~\cite{DK}.
The points $\ZZ{dG\sb{\DK}}$ coincide with the $\F\sb 4$-rational points of $\Pt$,
and hence the group $\Aut (X\sb{G\sb{\DK}}, \pol\sb{G\sb{\DK}})$ is equal to $\PGL (3, \F\sb 4)$.
We call $[G\sb{\DK}]$
the \emph{Dolgachev-Kondo point}.
\par
\medskip
Now we can state our main results.
\begin{theorem}\label{thm:main}
The locus $\moduli\sb{\le 2}$ is a union of three irreducible curves 
$\closure{\moduli}\sb {A}$, $\closure{\moduli}\sb {B}$ and $\closure{\moduli}\sb{C}$.
In $\moduli$,
they are situated in such a way that, set-theoretically, 
$$
\closure{\moduli}\sb{A} \cap \closure{\moduli}\sb{B} =
\closure{\moduli}\sb{B} \cap \closure{\moduli}\sb{C} =
\closure{\moduli}\sb{C} \cap \closure{\moduli}\sb{A}=
\{ [G\sb{\DK}] \}.
$$
\end{theorem}
For $T=A, B$ and $C$, we put
$$
\moduli\sb{T}:=\closure{\moduli}\sb T \sm \{[G\sb{\DK}]\}.
$$
Hence $\moduli\sb 2$ is the  disjoint union of $\moduli\sb{A}$, $\moduli\sb B$ and $\moduli\sb C$.
\begin{theorem}
For $T=A, B$ and $C$,
the curve  $\moduli\sb{T}$
is isomorphic to an affine line punctured at the origin.
\end{theorem}
%
%
%
%
%
We will describe the curves $\moduli\sb T$ more explicitly.
Let $\omega\in \F\sb 4$ be a primitive third  root of unity,
and let ${\bar\omega}$ be $\omega+1=\omega^2$.
\begin{theorem}\label{thm:A}
Let $\Gamma\sb A$ be the group
$$
\Bigl\{ \lambda, \lambda+1, \frac{1}{\lambda}, \frac{1}{\lambda+1}, \frac{\lambda}{\lambda+1}, \frac{\lambda+1}{\lambda}\Bigr\}
$$
acting on the punctured $\lambda$-line $\A\sp 1\sm\{0, 1\}=\Spec k [\lambda, 1/\lambda(\lambda+1)]$.
We put
$$
J\sb A:=\frac{(\lambda^2+\lambda+1)^3}{\lambda^2\, (\lambda+1)^2 }
$$
so that 
$ k \left[\lambda, {1}/{\lambda(\lambda+1)} \right] \sp{\Gamma\sb A} = k [J\sb A]$ holds.
We also put
\begin{equation*}\label{eq:eqA}
\GA [\lambda]:=XYZ\left (X+Y+Z\right )\left ({X}^{2}+{Y}^{2}+\left ({\lambda}^{2}+\lambda\right
){Z}^{2}+XY+YZ+ZX\right ).
\end{equation*}
Then there exists an isomorphism 
$$
\moduli\sb{A}\cong \Spec k [J\sb A, 1/J\sb A]
$$
such that the family $W^2=\GA [\lambda]$ of sextic double planes
over the finite Galois cover
$\A\sp 1\sm\{0,1,\omega,\bar\omega\} =\Spec k [\lambda, 1/(\lambda^4+\lambda)]$
of the moduli curve $\moduli \sb{A}$
yields the universal family of polarized supersingular $K3$ surfaces.
The points $\ZZ{d\GA [\lambda]}$ are given in Table~\ref{table:gammalambdaA}.
  The  origin $J\sb A=0$ corresponds to the Dolgachev-Kondo point.
\par
For $\alpha\in k\sm\{0, 1, \omega, \bar\omega\}$,
$\Aut (X\sb{\GA [\alpha]}, \pol\sb{\GA [\alpha]})$ is equal to
the group
\begin{equation}\label{eq:groupA}
\left\{
\;
\left.
\left[
\begin{tabular}{c|c}
$A$ & $\vcenter{\vbox {\hbox{$a$}\vskip 4pt\hbox{$b$} \vskip 4pt}}$ \\
\hline
$0$ $0$ & $1$
\end{tabular}
\right]
\in \PGL(3, k)
\;
\right|
\;
\hbox{$\vcenter{\vbox{\hbox{$A\in \GL(2, \F\sb 2)$,}\hbox{$a, b\in \{0, 1, \alpha, \alpha+1\}$}}
}$}
\;
\right\}
\end{equation}
of order $96$.
\end{theorem}
\begin{theorem}\label{thm:B}
We put
$$
Q\sb\lambda:= ({\bar\omega}\lambda+\omega )\,{X}^{2}+{\bar\omega}\,{Y}^{2}+\omega\lambda\,{Z}^{2}+
(\lambda+1 )\,XY+ ({\bar\omega}\lambda+\omega )YZ+ (\lambda+1)\,ZX ,
$$
and
$$
\GB [\lambda]:=XYZ\left (X+Y+Z\right ) Q\sb{\lambda}.
$$
Let $\Gamma\sb B$ be the group
\begin{multline*}
\left\{\,\,
\lambda,\,\,
\omega\,\lambda+1,\,\,
\frac{1}{\lambda+1},\,\,
{\frac {\lambda+{\bar\omega}}{\lambda+1}},\,\,
{\frac {{\bar\omega}\lambda+\omega}{\lambda}},\,\,
{\frac {{\bar\omega}}{\lambda}},\,\,
{\frac {\omega}{\lambda+{\bar\omega}}},\,\,
\right.\\ \left.
{\frac {{{\bar\omega}}\left (\lambda+1\right )}{\lambda+{\bar\omega}}},\,\,
{\frac {{\bar\omega}\lambda}{\lambda+1}},\,\,
{\frac {\lambda}{\lambda+{\bar\omega}}},\,\,
{\frac {\lambda+1}{\lambda}},\,\,
{\bar\omega}(\lambda+1)\,\,
\right\}
\end{multline*}
acting on the punctured $\lambda$-line
$\A\sp 1\sm\{0,1,\bar\omega\}=\Spec k [\lambda,1/\lambda(\lambda+1)(\lambda+{\bar\omega})]$. We put
$$
J\sb B:=\frac{(\lambda+\omega)^{12}}{\lambda^3(\lambda+1)^3(\lambda+{\bar\omega})^3}
$$
so that 
$k \left[\lambda, {1}/{\lambda(\lambda+1)(\lambda+{\bar\omega})}\right] \sp{\Gamma\sb B} = k [J\sb B]$ holds.
Then there exists an isomorphism 
$$
\moduli\sb {B}\cong \Spec k [J\sb B, 1/J\sb B]
$$
such that the family $W^2=\GB [\lambda]$ of sextic double planes over the finite Galois cover
$\A\sp 1\sm\{0,1,\omega,\bar\omega\} =\Spec k [\lambda, 1/(\lambda^4+\lambda)]$ of the moduli curve $\moduli\sb {B}$
yields the universal family of polarized supersingular $K3$ surfaces.
The points $\ZZ{d\GB [\lambda]}$ are given in Table~\ref{table:gammalambdaB}.
The origin  $J\sb B=0$ corresponds to the Dolgachev-Kondo point.
 \par
For any $\alpha\in k\sm\{0, 1, \omega, \bar\omega\}$,
$\Aut (X\sb{\GB [\alpha]}, \pol\sb{\GB [\alpha]})$ is equal to
the subgroup of $\PGL(3, k)$  generated by
\begin{equation}\label{eq:groupB}
		\left[
				\begin{matrix}
						0 & \omega & 0 \\
						\bar\omega & 1 & 0\\
						1 & 1 & 1
				\end{matrix}
		\right],
\quad
		\left[
				\begin{matrix}
						0 & 0 & \bar\omega \\
						1 & 1 & 1\\
						\omega & 0 & 1
				\end{matrix}
		\right]
\quand
		\left[
				\begin{matrix}
						1 & 0 & 0 \\
						\bar\omega & 1 & 0 \\
						\omega & 0 & 1
				\end{matrix}
		\right].
\end{equation}
In particular, $\Aut (X\sb{\GB [\alpha]}, \pol\sb{\GB [\alpha]})$ is isomorphic to the extended Heisenberg group 
of order $18$.
\end{theorem}
\begin{theorem}\label{thm:C}
Let $\Gamma\sb C$ be the group
$$
\set{\alpha\lambda+\beta}{\alpha\in \Ff\sptimes, \beta\in \Ff}
$$
of order $12$
acting on the  $\lambda$-line
$\A\sp 1= k[\lambda]$.
We put
$$
J\sb C:=(\lambda^4+\lambda)^3
$$
so that $k[\lambda]\sp{\Gamma\sb C}=k [J\sb C]$ holds.
We also put
\begin{equation*}\label{eq:eqC}
\GC  [\lambda]:=XYZ\left (X^3+Y^3+Z^3\right )+(\lambda^4+\lambda) X^3 Y^3.
\end{equation*}
Then there exists an isomorphism 
$$
\moduli\sb{C}\cong \Spec k [J\sb C,1/J\sb C]
$$ 
such that the family $W^2=\GC  [\lambda]$ of sextic double planes
over the finite Galois cover
$\A\sp 1\sm\{0,1,\omega,\bar\omega\} =\Spec k [\lambda, 1/(\lambda^4+\lambda)]$
of the moduli curve $\moduli\sb{C}$ yields the universal family of polarized supersingular $K3$ surfaces.
The points $\ZZ{d\GC  [\lambda]}$ are given in Table~\ref{table:gammalambdaC}.
  The origin $J\sb C=0$ corresponds to the Dolgachev-Kondo point.
\par
For $\alpha\in k\sm\{0, 1, \omega, \bar\omega\}$,
$\Aut (X\sb{\GC  [\alpha]}, \pol\sb{\GC [\alpha]})$ is equal to
\begin{equation}\label{eq:groupC}
\left\{\;\;
		\left.
				\left[
						\begin{matrix}
									a& b & 0 \\
									c& d& 0\\
									a^2 c^2 \alpha +e & b^2 d^2 \alpha + f & 1 
						\end{matrix}
				\right]\;\in\;\PGL(3, k)\;\;
		\right|\;
				\hbox{
						$\vcenter{
									\hbox{ $a,b,c,d,e,f \in \F\sb 4$,
									}
         \hbox{ $ad+bc=1$
									}
						}$
				}
\right\}
\end{equation}
of order $960$.
\end{theorem}
Next we consider the isomorphism classes of \emph{non-polarized} supersingular $K3$ surfaces
with Artin invariant $2$.
\begin{definition}
A reduced (possibly reducible)  curve $D$ in $\moduli\sb{T}\times \moduli\sb{T\sprime}$ is called 
a \emph{correspondence between $\moduli\sb{T}$ and $\moduli\sb{T\sprime}$}.
For a correspondence $D\subset \moduli\sb{T}\times \moduli\sb{T\sprime}$,
let $\spT D$ denote the correspondence in $\moduli\sb{T\sprime}\times \moduli\sb{T}$
obtained from $D$ by interchanging the first and the second factors.
When $D$ is a union of two curves $D\sb 1$ and $D\sb 2$ without common irreducible components,
we write $D=D\sb 1+D\sb 2$ and $D\sb 2=D-D\sb 1$.
Let 
$D\sb 1 \st \moduli\sb{T}\times \moduli\sb{T\sprime}$ and 
$D\sb 2 \st \moduli\sb{T\sprime}\times \moduli\sb{T\spprime}$
be  correspondences.
The \emph{composite} $D\sb 1 * D\sb 2 \st \moduli\sb{T}\times \moduli\sb{T\spprime}$
of $D\sb1$ and $D\sb2$ is defined to be the image of 
$$
(D\sb 1 \times \moduli\sb{T\spprime})\cap (\moduli\sb T \times D\sb 2) 
\;\;\st\;\; 
\moduli\sb T \times \moduli\sb{T\sprime}\times \moduli\sb{T\spprime}
$$
by the natural projection to $\moduli\sb T \times \moduli\sb{T\spprime}$.
\end{definition}
\begin{definition}
A correspondence  $D$ in $\moduli\sb{T}\times \moduli\sb{T\sprime}$
is called an \emph{isomorphism correspondence}
if, for every point $([G], [G\sprime])$ of $D$,
the supersingular $K3$ surfaces $X\sb G$ and $X\sb{G\sprime}$
(without polarization)
are isomorphic.
An isomorphism  correspondence 
$D\st \moduli\sb{T}\times \moduli\sb{T\sprime}$
is said to be \emph{trivial}
if $T$ is equal to $T\sprime$ and $D$ is the diagonal 
$\Delta\sb T$ of $\moduli\sb{T}\times\moduli\sb{T}$.
\end{definition}
Using  Cremona transformations
by quintic curves,
which played a central role in the study of $\Aut (X\sb{G\sb{\DK}})$  in~\cite{DK},
we have obtained   examples of non-trivial isomorphism correspondences.
\begin{definition}\label{def:center}
Let $G$ be a homogeneous polynomial in $\UUU$.
We say that a subset $\Sigma\subset \ZZ{dG}$ of cardinality $6$
is  a \emph{center of Cremona transformation for $(X\sb G, \pol\sb G)$ {\rm or} for $G$}
if $\Sigma$ satisfies the following  conditions:
\begin{itemize}
\item no three points of $\Sigma$ are collinear, and
\item for each   $p\sb i\in \Sigma$,
there exists a conic curve  $\EED\sprime\sb i\st \Pt$
such that $\EED\sprime\sb i\cap \ZZ{dG} =\Sigma\sm \{p\sb i\}$.
\end{itemize}
Note that the conic curve  $\EED\sprime\sb i$ is necessarily nonsingular.
\end{definition}
Let $\Sigma=\{p\sb 1, \dots, p\sb 6\}$ be a center of Cremona transformation for $(X\sb G, \pol\sb G)$.
Consider the linear system
$|\ideal \sb{\Sigma}^2 (5)| \st |\OOO\sb{\Pt} (5)|$
of quintic curves that pass through all the points of $\Sigma$ and are singular at  each point of $\Sigma$.
Then $|\ideal \sb{\Sigma}^2 (5)|$ is of dimension $2$,
and defines a birational map
$$
\cremona\sb{\Sigma}   \;:\;  \Pt   \;\cdots\to\;  \Pt.
$$
The birational map
$\cremona\sb{\Sigma}$
is the composite of the blowing up $\blu: S\to \Pt$ of the  points of $\Sigma$
and the blowing down $\bld: S\to\Pt$ of the strict transforms $\EED\sb i$ of the  conic curves $\EED\sprime\sb i$.
We denote by $p\sprime \sb i $
the image of $\EED\sb i$ by $\bld$.
Note that, if $p\in \Pt\sm \Sigma$, then the point $\cremona\sb{\Sigma} (p) \in \Pt$ is well-defined.
\begin{proposition}[Dolgachev-Kondo~\cite{DK}]\label{prop:Zsprime}
We put
$$
Z\sprime\;:=\;
\set{\cremona\sb{\Sigma} (p)}{ p\in \ZZ{dG} \sm \Sigma}\;\cup\;
\{p\sprime\sb 1, \dots, p\sprime\sb 6\}.
$$
Then there exists a homogeneous polynomial $G\sprime\in \UUU$
such that $Z\sprime=\ZZ{dG\sprime}$.
The birational map $\cremona\sb{\Sigma}$ of $\Pt$ lifts to an isomorphism
$$
\Cremona\sb{\Sigma} \;:\; X\sb{G}\isom X\sb{G\sprime}
$$
of supersingular $K3$ surfaces.
\end{proposition}
See also \S\ref{sec:isomcor} of this paper for the proof of Proposition~\ref{prop:Zsprime}.
Note that the polynomial $G\sprime$ is not uniquely determined, but the point $[G\sprime]\in\moduli$
is uniquely determined by $G$ and $\Sigma$.
We call $\Cremona\sb{\Sigma}$ the \emph{Cremona transformation of $X\sb G$ with center  $\Sigma$}.
\par
\medskip
Let $T$ be $A$, $B$ or $C$.
As Tables~\ref{table:gammalambdaA},~\ref{table:gammalambdaB} and~\ref{table:gammalambdaC} show,
the family
$$
\set{(p, \lambda)  }{p\in \ZZ{d\GT[\lambda]}}\;\;\st\;\;
\Pt\times (\A\sp 1\sm\{0,1,\omega,\bar\omega\} )
$$
of the points $\ZZ{d\GT[\lambda]}$ 
consists of $21$
connected components, each of which is \'etale of degree $1$ 
over the punctured 
$\lambda$-line $\A\sp 1\sm\{0,1,\omega,\bar\omega\}$. 
Therefore it makes sense to talk
about a family  $\Sigma[\lambda]$ of subsets of $ \ZZ{d\GT[\lambda]}$ that depends on $\lambda$ continuously.
It can be shown  that,
if $\Sigma[\alpha]$ is a center of  Cremona transformation for $\GT[\alpha]$ 
at  one $\alpha\in k\sm \{0,1,\omega,\bar\omega\}$,
then so is $\Sigma[\alpha]$
at every $\alpha\in k\sm \{0,1,\omega,\bar\omega\}$.
In this case,
we say that $\Sigma[\lambda]$
is a center of  Cremona transformation for $\GT[\lambda]$
or for $(X\sb{\GT[\lambda]},\pol\sb{\GT[\lambda]})$.
\par
\medskip
Suppose that $\Sigma[\lambda]$
is a center of  Cremona transformation for $\GT[\lambda]$.
Then there exist a family $G\sprime[\lambda]$ of homogeneous polynomials in $\UUU$
 and a family of  isomorphisms
$$
\mapisom{\Cremona\sb{\Sigma[\lambda]}}{X\sb{\GT[\lambda]}}{X\sb{G\sprime[\lambda]}}
$$
depending on the parameter  $\lambda$.
The points $[G\sprime[\lambda]]$ are of course contained in $\moduli\sb 2 =\moduli\sb A\sqcup \moduli\sb B\sqcup \moduli\sb{C}$.
Suppose that $[G\sprime[\lambda]]\in \moduli\sb{T\sprime}$.
Then the curve
$$
\set{([\GT[\lambda]], [G\sprime[\lambda]])\in \moduli\sb{T}\times\moduli\sb{T\sprime}}{\lambda\in \A\sp 1 \sm \{0,1,\omega,
\bar\omega\}}
$$
is an irreducible 
isomorphism correspondence  between $\moduli\sb{T}$ and $\moduli\sb{T\sprime}$.
\begin{theorem}\label{thm:correspondence}
{\rm (1)}
There exist $1644$ centers of Cremona transformation for the family  
$(X\sb{\GA [\lambda]},\pol\sb{\GA [\lambda]})$. 
They yield the following  isomorphism
correspondences:
\begin{itemize}
\item $156$ of them give the trivial correspondence $\Delta\sb A$, 
\item $144$ of them give the correspondence
$$
D\sb{A,A,1}\;\; : \;\;1+J\sb{A}\,J\sb{A}\sprime+{J\sb{A}}^{2}{J\sb{A}\sprime}^{2}+
{J\sb{A}}^{2}{J\sb{A}\sprime}^{3}+{J\sb{A}}^{3}{J\sb{A}\sprime}^{2}=0
$$
in $\moduli\sb A\times \moduli\sb A$,
\item $720$ of them give the correspondence
$$
D\sb{A,A,2}:= D\sb{A,A,1}*D\sb{A,A,1}-\Delta\sb A\;\;\st\;\; \moduli\sb A\times \moduli\sb A,
$$
\item $576$ of them give the correspondence
$$
D\sb{A, B, 1}\;\; : \;\;J\sb{B}+J\sb{A}\,J\sb{B}+J\sb{A}\,{J\sb{B}}^{2}+{J\sb{A}}^{2}J\sb{B}+{J\sb{A}}^{4}=0
$$
in $\moduli\sb A\times \moduli\sb B$,
\item $48$ of them give the correspondence
$$
D\sb{A, C, 1}\;\; : \;\;J\sb{C}+J\sb{A}+J\sb{A}\,J\sb{C}+{J\sb{A}}^{2}J\sb{C}+{J\sb{A}}^{4}{J\sb{C}}^{2}=0
$$
in $\moduli\sb A\times \moduli\sb C$.
\end{itemize}
{\rm (2)}
There exist $1374$ centers of Cremona transformation for $(X\sb{\GB [\lambda]}, \pol\sb{\GB [\lambda]})$.
They yield the following  isomorphism correspondences:
\begin{itemize}
\item $798$ of them give the trivial correspondence $\Delta\sb B$, 
\item $216$ of them give the  correspondence 
$$
D\sb{B, A, 1}:=\spT D\sb{A, B, 1}\;\;\st\;\;\moduli \sb{B}\times \moduli\sb{A},
$$
\item $360$ of them give the  correspondence
$$
D\sb{B,B,1}:=D\sb{B, A, 1}* D\sb{A, B, 1}-\Delta\sb{B}  \;\;\st\;\; \moduli \sb{B}\times \moduli\sb{B}.
$$
\end{itemize}
{\rm (3)}
There exist $2224$ centers of Cremona transformation for $(X\sb{\GC  [\lambda]}, \pol\sb{\GC  [\lambda]})$.
They yield the following  isomorphism correspondences:
\begin{itemize}
\item $1200$ of them give the trivial correspondence $\Delta\sb C$, 
\item $960$ of them give the  correspondence 
$$
D\sb{C, A, 1}:=\spT D\sb{A, C, 1}  \;\;\st\;\; \moduli \sb{C}\times \moduli\sb{A},
$$
\item $64$ of them give the  correspondence
$$
D\sb{C, C, 1}:=D\sb{C, A, 1}* D\sb{A, C, 1}-\Delta\sb{C}  \;\;\st\;\; \moduli \sb{C}\times \moduli\sb{C}.
$$
\end{itemize}
\end{theorem}
Starting from the isomorphism correspondences by Cremona transformation above,
making transposes and composites,
and taking irreducible components,
we obtain  non-trivial  irreducible isomorphism correspondences given in Table~\ref{table:coreqs}.
\begin{table}
$$
\renewcommand{\arraystretch}{1.4}
\begin{array}{|c|c|c|c|}
\hline
\coreq{T}{T\sprime}{\textrm{name}}{\textrm{the equation}}
\coreq{A}{A}{D_{A,A,1}}{{{\JA}}^{3}{{\JAsprime}}^{2}+{{\JA}}^{2}{{\JAsprime}}^{3}+
{{\JA}}^{2}{{\JAsprime}}^{2}+{\JA}\,{\JAsprime}+1
}
\coreq{A}{A}{D_{A,A,2}}{{{\JA}}^{6}{{\JAsprime}}^{2}+{{\JA}}^{4}{{\JAsprime}}^{4}+
{{\JA}}^{2}{{\JAsprime}}^{6}+{{\JA}}^{4}{{\JAsprime}}^{3}+
{{\JA}}^{3}{{\JAsprime}}^{4}+{{\JA}}^{4}{{\JAsprime}}^{2}+
{{\JA}}^{3}{{\JAsprime}}^{3}+{{\JA}}^{2}{{\JAsprime}}^{4}+
{{\JA}}^{4}{\JAsprime}+{\JA}\,{{\JAsprime}}^{4}+{{\JA}}
^{3}{\JAsprime}+{{\JA}}^{2}{{\JAsprime}}^{2}+{\JA}\,{{\JAsprime}}^{3}+{{\JA}}^{3}+{{\JA}}^{2}{\JAsprime}+{\JA}\,{
{\JAsprime}}^{2}+{{\JAsprime}}^{3}
}
\coreq{B}{B}{D_{B,B,1}}{{{\JB}}^{4}{\JBsprime}+{{\JB}}^{3}{{\JBsprime}}^{2}+{{\JB}}^{2}{{\JBsprime}}^{3}+{\JB}\,{{\JBsprime}}^{4}+{{\JB}}
^{3}{\JBsprime}+{{\JB}}^{2}{{\JBsprime}}^{2}+{\JB}\,{{\JBsprime}}^{3}+1
}
\coreq{C}{C}{D_{C,C,1}}{{{\JC}}^{4}{{\JCsprime}}^{4}+{{\JC}}^{3}{\JCsprime}+{{\JC}}^{2}{{\JCsprime}}^{2}+{\JC}\,{{\JCsprime}}^{3}+{{\JC}}
^{3}+{{\JC}}^{2}{\JCsprime}+{\JC}\,{{\JCsprime}}^{2}+{{
\JCsprime}}^{3}
}
\coreq{A}{B}{D_{A,B,1}}{{{\JA}}^{4}+{{\JA}}^{2}{\JB}+{\JA}\,{{\JB}}^{2}+{\JA
}\,{\JB}+{\JB}
}
\coreq{A}{B}{D_{A,B,2}}{{{\JA}}^{6}{\JB}+{{\JA}}^{5}{\JB}+{{\JA}}^{4}{{\JB}}
^{2}+{{\JA}}^{3}{{\JB}}^{3}+{{\JA}}^{4}{\JB}+{{\JA}}^{2
}{{\JB}}^{2}+{{\JA}}^{2}{\JB}+{\JA}\,{\JB}+1
}
\coreq{B}{C}{D_{B,C,1}}{{\JB}\,{\JC}+1
}
\coreq{B}{C}{D_{B,C,2}}{{{\JB}}^{4}{{\JC}}^{3}+{{\JB}}^{3}{{\JC}}^{3}+{{\JB}}^{
3}{{\JC}}^{2}+{{\JB}}^{2}{{\JC}}^{2}+{{\JC}}^{4}+{{\JB}
}^{2}{\JC}+{\JB}\,{\JC}+{\JB}
}
\coreq{C}{A}{D_{C,A,1}}{{{\JC}}^{2}{{\JA}}^{4}+{\JC}\,{{\JA}}^{2}+{\JC}\,{\JA}+{\JC}+{\JA}
}
\coreq{C}{A}{D_{C,A,2}}{{{\JC}}^{2}{{\JA}}^{6}+{{\JC}}^{2}{{\JA}}^{5}+{{\JC}}^{
2}{{\JA}}^{4}+{\JC}\,{{\JA}}^{4}+{{\JC}}^{2}{{\JA}}^{2}
+{{\JC}}^{3}+{{\JC}}^{2}{\JA}+{\JC}\,{{\JA}}^{2}+{{\JA}}^{3}
}
\end{array}
$$
\caption{Non-trivial  irreducible isomorphism correspondences}\label{table:coreqs}
\end{table}
%
%
When $T\ne T\sprime$,
we denote by $D\sb{T\sprime, T, \nu}$ the correspondence $\spT D\sb{T, T\sprime, \nu}$ for $\nu=1$ and $2$.
They have the relations in Table~\ref{table:correls} at the end of~\S\ref{sec:isomcor}.
\begin{question}
Are there any non-trivial irreducible  isomorphism correspondences 
other than the ones in Table~\ref{table:coreqs} and their transposes?
\end{question}
The Cremona transformations that yield the trivial isomorphism
correspondence are also interesting,
because they give automorphisms of the supersingular $K3$ surface $X$
that may not be contained in $\Aut\Xpol$.
See Remark~\ref{rem:CRasAut}.
\begin{observation}
Consider a Cremona transformation $\Cremona\sb{\Sigma}$ on $(X\sb{\GA [\lambda]}, \pol\sb{\GA [\lambda]})$
that yields the non-trivial isomorphism correspondence $D\sb{A,A, 1}$.
The curve $D\sb{A,A, 1}$ intersects the diagonal $\Delta\sb A$ at two points
$(J\sb A, J\sb A\sprime)=(\omega,\omega)$ and $(\bar\omega,\bar\omega)$.
Let $\eta$ be an element of $k$ such that the $J\sb A$-invariant of
$(X\sb{\GA [\eta]}, \pol\sb{\GA [\eta]})$ is $\omega$ or $\bar\omega$;
that is, $\eta$ is a  root of
$$
(\lambda^4+\lambda^3+1)(\lambda^4+\lambda+1)(\lambda^4+\lambda^3+\lambda^2+\lambda+1)=0.
$$
The Cremona transformation $\Cremona\sb{\Sigma}$ gives rise to an automorphism of $X\sb{\GA [\eta]}$,
which cannot be deformed to any automorphisms of $X\sb{\GA [\lambda]}$ for a generic $\lambda$.
In other words,
the automorphism group $\Aut (X\sb{\GA[\lambda]})$ of the non-polarized supersingular $K3$ surface $X\sb{\GA[\lambda]}$ 
\emph{jumps} at $\lambda=\eta$,
even though the numerical N\'eron-Severi lattice of $X\sb{\GA[\lambda]}$ is constant around $\lambda=\eta$.
Note that the automorphism group of a supersingular $K3$ surface
is always embedded into
the orthogonal group of its numerical N\'eron-Severi lattice (\cite[\S8, Proposition 3]{RS}).
\end{observation}
The plan of this paper is as follows.
In \S\ref{sec:review}, we recall from~\cite{Shimada2004} 
the definition of the binary code associated with a polarized supersingular $K3$ surface of type $(\sharp)$.
We stratify the moduli space $\moduli$ according to the isomorphism classes $[\abscode]$
of the codes,
and give a method to construct the  stratum $\moduli\sb{[\abscode]}$
from the code $\abscode$.
In \S\ref{sec:1and2},
we present three isomorphism classes $[\CA]$, $[\CB]$ and $[\CC]$ of codes 
that are associated with polarized supersingular $K3$ surfaces of type $(\sharp)$
with Artin invariant $2$.
In \S\ref{sec:A}, \S\ref{sec:B} and \S\ref{sec:C},
we carry out the method of the construction of $\moduli\sb{[\abscode]}$ for $\abscode=\CA, \CB$ and $\CC$,
and prove Theorems~\ref{thm:A},~\ref{thm:B} and~\ref{thm:C}, respectively.
In \S\ref{sec:cremona},
we review from~\cite{DK} the theory of Cremona transformations by quintic curves.
In \S\ref{sec:isomcor},
we explain the algorithm to  calculate
the isomorphism correspondences
given by Cremona transformations,
and prove Theorem~\ref{thm:correspondence}.
%
%
%
\par
\medskip
The isomorphism classes of codes associated with polarized supersingular $K3$ surfaces  
of Artin invariant $\sigma\ge 3$ are
also given in~\cite{Shimada2004}. For $\sigma=3$, there are $13$ isomorphism classes, and
for $\sigma=4$, there are $41$ isomorphism classes.
It would be a challenging problem in computational algebraic geometry
to construct explicitly the moduli spaces of dimension $\sigma-1$ corresponding to 
these  isomorphism classes of codes,
and to investigate the relations between them.
\par
\medskip
In~\cite{RS_char2},
Rudakov and Shafarevich gave  
explicitly families of supersingular $K3$ surfaces in characteristic $2$
for  Artin invariants $\sigma=1, \dots, 10$.
The equation of the family for $\sigma=2$ is
$$
y^2=x^3+\mu\, t^6 x + t^5 (t+1)^4,
$$
where $\mu$ is the ``modulus".
We would like to know the relation between $\mu$ and our moduli $\JA$, $\JB$ and $\JC$. 
\par
\medskip
The polarized supersingular $K3$ surface of type $(\sharp)$ is 
an example of \emph{Zariski surfaces}.
A general theory of Zariski surfaces has been developed in~\cite{BL}.
\par
\bigskip
{\bf Notation and terminologies.}
\par
(1)
Let $A$ be a commutative ring, and $S$  a set.
We denote by $A\sp S$ the $A$-module of all maps from $S$ to $A$.
\par
(2)
Let $S$ be a finite set.
The full symmetric group of $S$ is denoted by $\SSSS (S)$,
which acts on $S$ from left.
We denote by $\Pow (S)$ the power set of $S$.
A canonical identification  between $\Pow (S)$ and $\F\sb 2\sp S$ is given by 
$f\in \F\sb 2 \sp S\mapsto f\inv (1) \st S$. 
Hence $\Pow (S)$ has a structure of the $\F\sb 2$-vector space by the symmetric difference
$$
T\sb 1 + T\sb 2 =(T\sb 1 \cup T\sb 2)\sm (T\sb 1 \cap T\sb 2)\qquad (T\sb 1, T\sb 2\st S).
$$
A linear subspace of $\F\sb 2 \sp S=\Pow (S)$ is called a \emph{code},
and an element of a code is called a \emph{word}.
A word is expressed either as  a vector of dimension $\card{S}$ with coefficients in $\F\sb 2$,
or as a subset of $S$.
The cardinality $\card{A}$ of a word $A\st S$ is called the \emph{weight} of $A$.
The automorphism group $\Aut (\abscode)$
of a code $\abscode \st\Pow (S)$ 
is the subgroup  of $\SSSS (S)$
consisting of all permutations
preserving $\abscode$.
Two codes $\abscode$ and $\abscode\sprime$ in $\Pow (S) $ are said to be
\emph{isomorphic} if there exists a permutation $\sigma \in \SSSS (S)$
such that $\sigma (\abscode)=\abscode\sprime$.
The  isomorphism class of codes represented by a code $\abscode$ is denoted by
$[\abscode]$.
\par
(3)
A \emph{lattice} is a free $\Z$-module $\Lambda$ of finite rank
equipped with a non-degenerate symmetric bilinear form
$\Lambda\times \Lambda\to \Z$.
A lattice is called \emph{even} if $v^2\in 2\Z$
holds for every $v\in \Lambda$.
A lattice is called \emph{hyperbolic} if the signature of the symmetric bilinear form 
on $\Lambda\otimes \R$
is $(1, r-1)$, where $r$ is the rank of $\Lambda$.
The \emph{dual lattice} $\Lambda\dual$ of $\Lambda$
is the $\Z$-module $\Hom (\Lambda, \Z)$.
There exists  a canonical embedding
$\Lambda\inj \Lambda\dual$
of finite cokernel.
Hence $\Lambda\dual$ can be regarded as a submodule of $\Lambda\otimes\sb{\Z}\Q$. 
We have a natural $\Q$-valued symmetric bilinear form on $\Lambda\dual$
that extends the $\Z$-valued bilinear form on $\Lambda$.
An \emph{overlattice} of $\Lambda$ is a submodule $\Lambda\sprime $ of $\Lambda\dual$
containing $\Lambda$
such that the canonical $\Q$-valued symmetric bilinear form on $\Lambda\dual$
takes values in $\Z$ on $\Lambda\sprime$.
\section{The codes associated with the supersingular $K3$ surfaces}\label{sec:review}
First we give a proof of Proposition~\ref{prop:equivG}.
\begin{proof}[Proof of Proposition~\ref{prop:equivG}]
The equivalence of (i) and (iii) follows from the structure
of the graded ring $\oplus\sb{m\ge 0} H\sp 0 (X, \pol\sp{\otimes m})$,
where $X$ is a $K3$ surface and $\pol$ is a line bundle of degree $2$.
(See \cite[\S7]{Shimada2004}.)
By~\cite[Theorem 2.1]{Shimada2004},
$\ZZ{dG}=\ZZ{dG\sprime}$ holds if and only if $dG=c\cdot dG\sprime $
for some $c\in k\sptimes$.
Since the kernel of $G\mapsto dG$ is equal to $\{H^2|H\in\VVV\}$,
the equivalence of (ii) and (iii) follows.
\end{proof}
\subsection{Definition of the code $\gcodeXpol$}
Let us fix a finite set 
$$
\Ps:=\{ P\sb 1, \dots, P\sb{21}\}
$$
consisting of $21$ elements,
on which the full symmetric group $\SSSS (\Ps)$ acts from left.
\begin{definition}
We denote by $\gs$ the space  of all injective maps
$\gamma: \Ps \inj\Pt$
such that there exists a homogeneous polynomial $G\in \UUU$
satisfying
$\gamma (\Ps)=\ZZ{dG}$.
\end{definition}
The space $\gs$ are constructed as follows.
For $G\in \UUU$,
let $\langle G\rangle\in \P\sb * (\UUU/\VVV)$ denote the point corresponding to $G$.
We denote by 
$$
\ZZZ:=\set{(p, \langle G\rangle)\in \Pt\times \P\sb * (\UUU/\VVV)}{p\in \ZZ{dG}}\;\;\to\;\; \P\sb * (\UUU/\VVV)
$$
the family of $\ZZ{dG}$,
which is finite and \'etale of degree $21$ over $\P\sb * (\UUU/\VVV)$.
We prepare $21$ copies of $\ZZZ$ and make the fiber-product  $\ZZZ\sp{(21)}$ of them over $\P\sb * (\UUU/\VVV)$.
Then $\gs$ is the union of irreducible components
of $\ZZZ\sp{(21)}$ that do not intersect the big diagonal.
\begin{remark}
We fix a base point $\langle G\sb 0\rangle \in \P\sb * (\UUU/\VVV)$,
and consider the monodromy action 
$$
\map{\mu}{\pi\sb 1(\P\sb * (\UUU/\VVV), \langle G\sb 0\rangle)}{\SSSS (\ZZ{dG\sb 0})}
$$
of the algebraic fundamental group of $\P\sb * (\UUU/\VVV)$ on $\ZZ{dG\sb 0}$.
Then the number of irreducible components of $\gs$ is equal to the index of the image of $\mu$ in $\SSSS (\ZZ{dG\sb 0})$.
It was shown in~\cite[Chapter 4, Appendix 2]{BL} that 
the monodromy group on the singular points  of a generic Zariski surface in characteristic $\ge 5$
is equal to the full-symmetric group.
\end{remark}
The group $\SSSS (\Ps)$ acts on $\gs$ from right,
and $\PGL (3, k)$ acts on $\gs$ from left.
By Proposition~\ref{prop:equivG}, we have
$$
\moduli =\PGL (3, k)\backslash  \gs/\SSSS (\Ps).
$$
\par
\medskip
Let 
$$
N\sb 0:=\Z\sp{\Ps} \oplus \Z h= \bigoplus\sb{i=1}\sp{21} \Z e\sb i \oplus \Z h
$$
be  a free $\Z$-module of rank $22$ generated by  vectors
$e\sb 1, \dots, e\sb {21}$ corresponding to
$P\sb 1, \dots, P\sb{21}\in \Ps$
and a vector $h$.
We equip $N\sb 0$ with a structure of the even hyperbolic lattice by
$$
e\sb i^2=-2, \quad
h^2=2, \quad
e\sb i e\sb j=0 \;\;(i\ne j), \quad
h e\sb i=0.
$$
The dual lattice
$$
N\sb 0 \dual =\Hom (N\sb 0, \Z) \;\;\subset\;\; N\sb 0 \otimes\sb{\Z} \Q
$$
is generated by $e\sb i/2$ $(i=1, \dots, 21)$
and $h/2$.
Thus we have a canonical isomorphism
$$
N\sb 0\dual /N\sb 0 \;\cong\; \F\sb 2 \sp{\Ps} \oplus \F\sb 2 \;=\; \Pow (\Ps) \oplus \F\sb 2.
$$
Hence 
we can write an element of $N\sb 0\dual/N\sb 0$  in the form
$(A, \alpha)$, 
where $A$ is a subset of $\Ps$ and $\alpha\in \F\sb 2$.
We denote by
$$
\map{\pr}{N\sb 0\dual}{N\sb 0\dual/N\sb 0=\Pow (\Ps) \oplus \F\sb 2}
$$
the natural projection. 
We also denote by
$$
\map{\rho}{N\sb 0\dual/N\sb 0=\Pow (\Ps) \oplus \F\sb 2}{\Pow (\Ps)}
$$
the natural projection onto the first factor.
The following is obvious:
\begin{lemma}\label{lem:even}
Let $\widetilde{\gcode}$ be a subspace of the $\F\sb 2$-vector space $\Pow (\Ps)  \oplus \F\sb 2$.
Then the submodule $\pr\inv (\widetilde{\gcode})$ of $N\sb 0\dual$  is an even overlattice of $N\sb 0$
if and only if
$$
\card{A}\equiv\begin{cases}
0\, \bmod 4  & \textrm{if $\alpha=0$}\\
1\,  \bmod 4  & \textrm{if $\alpha=1$}
\end{cases}
$$
holds for every $(A, \alpha)\in \widetilde{\gcode}$.
\qed
\end{lemma}
\par
\medskip
Let $\Xpol$ be a polarized supersingular $K3$ surface of type $\typesharp$,
and let $\NS (X)$ denote the numerical N\'eron-Severi lattice of $X$.
There exists $G\in \UUU$ such that  $\Phipol: X\to\Pt$
factors through $\pi\sb G: Y\sb G\to \Pt$.
We put
$$
\ZXpol:= \ZZ{dG}=\pi\sb{G} (\Sing \YG).
$$
There also exists a point $\gamma:\Ps\inj\Pt$ of $\gs$,
unique up to the action of $\SSSS(\Ps)$, 
that induces a bijection
from  $\Ps$ to $\ZXpol$.
We fix such a point $\gamma\in \gs$.
Let $E\sb i$ be the $(-2)$-curve on $X$ such that
$\Phipol (E\sb i)$ is the point $\gamma (P\sb i)\in \ZXpol$.
Then 
we obtain an embedding
$$
\iota\sb{\gamma}\;\;:\;\; N\sb 0\;\; \inj \;\;  \NS (X)
$$
of the lattice $N\sb 0$ into $\NS(X)$  by $e\sb i \mapsto [E\sb i]$ and $h\mapsto [\pol]$.
By the embedding $\iota\sb{\gamma}$,
we can regard $\NS (X)$ as a submodule of $N\sb 0\dual$.
We put
\begin{eqnarray*}
\widetilde{\gcode} (X, \pol, \gamma)&:=&\NS(X)/N\sb 0 \;\;\st\;\; \Pow (\Ps) \oplus \F\sb 2,
\quand \\
\gcodeXpol &:=&\rho(\widetilde {\gcode} (X, \pol, \gamma))\;\;\st\;\; \Pow (\Ps).
\end{eqnarray*}
Since $\NS (X)$ is an even overlattice of $N\sb 0$,
the code $\widetilde {\gcode} (X, \pol, \gamma)$  is uniquely recovered  from 
$\gcodeXpol$ by Lemma~\ref{lem:even},
and hence the lattice $\NS (X)$ is  also uniquely recovered from the code
$\gcodeXpol$.
In particular, the Artin invariant $\sigma (X)$ of $X$ is given by 
$$
\sigma (X) = 11-\dim \sb{\F\sb 2} \gcodeXpol.
$$
Note that the isomorphism class of the code $\gcodeXpol$
does not depend on the choice of $\gamma$.
The following is one of the main results of~\cite{Shimada2004}:
\begin{theorem}\label{thm:codes}
For an isomorphism class $[\abscode]$ of codes in $\Pow (\Ps)$,
the following two conditions are equivalent.

{\rm (i)} There exists a polarized supersingular $K3$ surface $\Xpol$ of type $\typesharp$ such that,
for a {\rm(}and hence any{\rm)} bijection $\gamma$
from  $\Ps$ to $\ZXpol$,
the code $\gcodeXpol$ is in the isomorphism class $[\abscode]$.

{\rm (ii)} A {\rm (}and hence any{\rm)} code $\abscode \in [\abscode]$ satisfies the following:
\begin{itemize}
\item $\dim \abscode \le 10$,
\item the word $\Ps\in \Pow (\Ps)$ is contained in $\abscode$, and
\item$\card{A}\in \{ 0,5,8,9,12,13,16,21\}$ for every word $A\in \abscode$.
\end{itemize}
\end{theorem}
\subsection{Geometry of $\ZXpol$ and the code $\gcodeXpol$}\label{subsec:geom}
Let $\Xpol$ be a polarized supersingular $K3$ surface of type $\typesharp$.
We fix a bijection $\gamma$ from $\Ps$ to $\ZXpol$.
Let $G\in \UUU$ be a homogeneous polynomial such that
$\Phipol$ factors through $\YG$,
or equivalently,
such that $\ZZ{dG}=\ZXpol$ holds.
For the proofs of  the facts stated in this subsection,
we refer  the reader to \cite[\S6 and \S7]{Shimada2004}.
\begin{definition}
Let $C\subset\Pt$ be a reduced irreducible curve.
We say that $C$ \emph{splits in $\Xpol$}
if the proper transform of $C$ by $\Phipol : X\to \Pt$ is non-reduced.
We say that a reduced (possibly reducible) curve $C\sprime$
\emph{splits in $\Xpol$}
if every irreducible component of $C\sprime$ splits in $\Xpol$.
\end{definition}
Since $\Phipol$ is purely inseparable of degree $2$,
 the proper transform of a splitting curve $C$ by $\Phipol$ is written as 
$2F\sb C$,
where $F\sb C$ is a reduced divisor of $X$.
We denote by $w(C)\in \Pow (\Ps)$ 
the image of the numerical equivalence class $[F\sb C]\in \NS (X)$ by 
$$
\NS (X) \;\; \maprightsp{}\;\; \NS (X)/N\sb 0 \;\;\inj\;\; N\sb 0\dual /N\sb 0 \;\;\maprightsp{\rho}\;\;\Pow (\PPP),
$$
where $N\sb 0\inj \NS (X)$ is obtained from the fixed bijection $\gamma:\Ps \isom \ZXpol$.
By definition, we have 
$$
w(C)\in \gcodeXpol.
$$
It is easy to see that
\begin{equation*}\label{eq:word}
w (C)=\set{P\sb i\in \PPP}{\hbox{the multiplicity of $C$ at $\gamma (P\sb i)$ is odd}}.
\end{equation*}
If $C$ is a nonsingular curve splitting in $\Xpol$,
then 
$$
w(C)=\gamma\inv (C\cap \ZXpol).
$$
If $C\sb 1$ and $C\sb 2$ are two splitting curves without common irreducible components,
then $w (C\sb1 \cup C\sb2)=w(C\sb1) +w(C\sb2)$ holds.
\begin{proposition}\label{prop:quintic}
Let $\ideal \sb{\ZZ{dG}}\st\OOO\sb{\Pt}$ be the ideal sheaf defining the subscheme $Z(dG)$.
The linear system $|\ideal \sb{\ZZ{dG}} (5)|$ of 
quintic curves  passing through all the points of $\ZZ{dG}$
is of dimension $2$, and spanned by the curves defined by 
$$
{\partial G}/{\partial X}=0,\quad 
{\partial G}/{\partial Y}=0 \quand 
{\partial G}/{\partial Z}=0.
$$
A general member $C$ of $|\ideal \sb{\ZZ{dG}} (5)|$ splits in $\Xpol$,
and the word $w(C)\in \gcodeXpol$ is equal to $\Ps\in\Pow (\Ps)$.
\end{proposition}
\begin{proposition}\label{prop:nodes}
Let $C$ be a reduced curve splitting in $\Xpol$,
and let $p$ be a point of $C$.
{\rm (1)}
If $p$ is an ordinary node of $C$, then $p\in \ZXpol$.
{\rm (2)}
If $p$ is an ordinary tacnode of $C$, then $p\notin \ZXpol$.
\end{proposition}
\begin{proposition}\label{prop:nodalsplit}
Let $C$ be a reduced curve of degree $6$ splitting in $\Xpol$,
and let $G\sprime =0$ be a defining equation of $C$.
If $C$ has only ordinary nodes as its singularities,
then the homogeneous polynomial $G\sprime$
is a point of  $\UUU$,
and the point $[G\sprime]\in \moduli$
corresponds to the isomorphism class of $\Xpol$.
\end{proposition}
\begin{proposition}\label{prop:linewt}
Let $L\st\Pt$ be a line.
The following conditions are equivalent:
{\rm (i)} $L$ splits in $\Xpol$, 
{\rm (ii)} $\card{L\cap\ZXpol}\ge 3$, and 
{\rm (iii)} $\card{L\cap\ZXpol}=5$.
\end{proposition}
\begin{proposition}\label{prop:conicwt}
Let $Q\st\Pt$ be a nonsingular conic curve.
The following conditions are equivalent:
{\rm (i)}  $Q$ splits in $\Xpol$, 
{\rm (ii)} $\card{Q\cap\ZXpol}\ge 6$, and 
{\rm (iii)} $\card{Q\cap\ZXpol}=8$.
\end{proposition}
\begin{corollary}
The word $w(L)=\gamma\inv (L\cap \ZXpol) $ of
a splitting line $L$ is of weight $5$,
and 
the word $w(Q)=\gamma\inv (Q\cap \ZXpol) $ of
a splitting nonsingular  conic curve $Q$ is of weight $8$.
\end{corollary}
\begin{definition}
A pencil $\EEE$ 
of cubic curves in $\Pt$  is called a {\it regular pencil}
if the following hold:
\begin{itemize}
\item the base locus $\Bs (\EEE)$  of $\EEE$  consists of distinct $9$ points, and 
\item  every singular member of $\EEE$  has only one ordinary node as its singularities.
\end{itemize}
We say that 
a regular  pencil $\EEE$ \emph{splits in $\Xpol$} if 
every member of $\EEE$ splits in $\Xpol$.
\end{definition}
\begin{proposition}
Let $\EEE$ be a regular pencil of cubic curves spanned by $E\sb 0$ and $E\sb\infty$.
Let $H\sb 0=0$ and $H\sb\infty=0$ be the defining equations of $E\sb 0$ and $E\sb\infty$, respectively.
Then $\EEE$ splits in $\Xpol$ if and only if 
there exist $c\in k\sptimes$ and $H\in \VVV$
such that
\begin{equation}\label{eq:EEEsplits}
G=c H\sb 0 H\sb{\infty} + H^2
\end{equation}
holds.
If $\EEE$ splits in $\Xpol$,
then $\Bs(\EEE)$ is contained in $\ZXpol$,
and 
$$
w (E\sb t) =\gamma\inv (\Bs (\EEE))
$$
holds for every member $E\sb t$ of $\EEE$.
In particular, the word $w(E\sb t)$ is of weight $9$.
\end{proposition}
\begin{remark}
The condition~\eqref{eq:EEEsplits} is equivalent to
$$
\ZZ{d (H\sb 0 H\sb{\infty})}=\ZZ{dG}=\ZXpol
$$
by Proposition~\ref{prop:equivG}.
\end{remark}
\begin{remark}
A regular pencil $\EEE$
has $12$  singular members  $E\spar{1}, \dots, E\spar{12}$.
We denote by $N\spar{i}$ the ordinary node of $E\spar{i}$.
 Suppose that $\EEE$
splits in $\Xpol$.
Then 
$\ZXpol$
is a disjoint union of $\Bs(\EEE)$ and $\{ N\spar{1}, \dots, N\spar{12}\}$.
\end{remark}
Let $L\sb 1$ and $L\sb 2$ be  distinct lines splitting in $\Xpol$.
Then the intersection point of $L\sb 1$ and $L\sb 2$ is in $\ZXpol$ by Proposition~\ref{prop:nodes},
and hence
$$
w(L\sb 1 \cup L\sb 2 )=w(L\sb 1) + w(L\sb 2)
$$
is a word of weight $8$.

Let $L\sb 1$,  $L\sb 2$ and $L\sb 3$ be  lines splitting in $\Xpol$
such that $L\sb 1\cap L\sb 2 \cap L\sb 3=\emptyset$.
Then the three ordinary nodes  of $L\sb 1 \cup L\sb 2 \cup L\sb 3$ 
are in $\ZXpol$ by Proposition~\ref{prop:nodes},
and hence
$$
w(L\sb 1 \cup L\sb 2 \cup L\sb 3 )=w(L\sb 1) + w(L\sb 2)+w(L\sb 3)
$$
is a word of weight $9$.

Let $Q$ be a nonsingular conic curve splitting in $\Xpol$,
and let $L$ be a line splitting in $\Xpol$.
Using Proposition~\ref{prop:nodes},
we see that 
$L$ intersects $Q$ transversely if and only if
$w(L\cup Q)=w(L)+w(Q)$ is of weight $9$.
We also see that 
$L$ is tangent to  $Q$  if and only if
$w(L)\cap w(Q)=\emptyset$.
\begin{definition}
Let $\abscode \st \Pow (\Ps)$ be a code satisfying the conditions in (ii) of Theorem~\ref{thm:codes},
and let $A$ be a word of $\abscode$ with $\card{A}\in \{ 5, 8, 9\}$.

(i) We say that $A$ is a \emph{linear word of $\abscode$} if $\card{A}=5$.

(ii) Suppose $\card{A}=8$.
If $A$ is \emph{not} a sum of two linear words  of $\abscode$,
then we say that $A$ is a \emph{quadratic word  of $\abscode$}.

(iii)
Suppose $\card{A}=9$.
If $A$ is neither  a sum of three linear words  of $\abscode$
nor a sum of a linear  and a quadratic word  of $\abscode$,
then we say that $A$ is  a \emph{cubic word  of $\abscode$}.
\end{definition}
\begin{proposition}~\label{prop:wordsbijection}
{\rm (1)}
The correspondence 
$L\mapsto w(L)$
yields a bijection
from the set of lines splitting in $\Xpol$ 
to the set of linear words in $\gcodeXpol$.

{\rm (2)}
The correspondence 
$Q\mapsto w(Q)$
yields a bijection
from the set of nonsingular conic curves splitting in $\Xpol$ 
to the set of quadratic words in $\gcodeXpol$.

{\rm (3)}
The correspondence 
$\EEE\mapsto \gamma\inv (\Bs(\EEE))$
yields a bijection
from the set of regular pencils of cubic curves  splitting in $\Xpol$ 
to the set of cubic  words in $\gcodeXpol$.
\end{proposition}
By Theorem~\ref{thm:codes},
the code  $\gcodeXpol$ is generated by the word $\Ps$ and by the linear, quadratic and cubic words in $\gcodeXpol$.
Combining this fact with Proposition~\ref{prop:wordsbijection},
we obtain the following:
\begin{corollary}\label{cor:projaut}
Let $g$ be an element of the group 
$$
\Aut\Xpol=\set{h\in \PGL (3, k)}{ h (\ZXpol)=\ZXpol}.
$$
Then we have $\gcode(X,\pol,\gamma)=\gcode(X,\pol, g\circ \gamma)$.
Hence there exists a unique element
$\sigma\sb g \in \Aut (\gcodeXpol)$
such that
$g\circ \gamma =\gamma \circ \sigma\sb{g}$
holds.
By $g\mapsto \sigma\sb g $,
we can embed $\Aut \Xpol$ into 
$\Aut (\gcodeXpol)$.
\end{corollary}
\subsection{Construction of $\moduli\sb{[\abscode]}$ from $\abscode$}
%
%
Let $[\abscode]$ be an isomorphism class of codes satisfying the conditions in (ii) of Theorem~\ref{thm:codes}.
We denote by
$$
\moduli\sb{[\abscode]} \;\;\st\;\;\moduli
$$
the locus of all isomorphism classes 
of polarized supersingular $K3$ surfaces $\Xpol$ of type $\typesharp$ such that
$\gcodeXpol$ is contained in $[\abscode]$
for a (and hence any) bijection 
$\gamma$ from $\Ps$ to $\ZXpol$.
We also denote by $\gs\sb{[\abscode]}$ the pull-back of $\moduli\sb{[\abscode]}$
by the quotient map
$$
\gs\;\;\maprightsp{}\;\;\moduli =\PGL (3, k) \backslash  \gs/\SSSS (\Ps).
$$
We will describe the locus $\gs\sb{[\abscode]}$.
\begin{definition}
For a point  $\gamma$ of $\gs$,
let $\gcode[\gamma]$ denote the code in $\Pow (\Ps)$ generated by the following words:
\begin{itemize}
\item $\Ps\in \Pow (\Ps)$,
\item words $A$ of weight $5$ such that the points  $\gamma (A)$ are collinear,
\item words $A$ of weight $8$ such that there exists a nonsingular conic curve containing $\gamma (A)$, and 
\item words $A$ of weight $9$ such that there exists a regular pencil $\EEE$ of cubic curves spanned by
$E\sb 0 =\{H\sb 0=0\}$ and $E\sb{\infty} =\{H\sb {\infty}=0\}$ such that $\Bs (\EEE)=\gamma(A)$ and 
$\ZZ{d(H\sb 0 H\sb{\infty})}=\gamma (\Ps)$ hold.
\end{itemize}
\end{definition}
From the results above,
we obtain the following:
\begin{corollary}
Suppose that $\gamma\in \gs$, and 
let $\Xpol$ be a polarized supersingular $K3$ surface of type $(\sharp)$ such that
$\gamma (\Ps)=\ZXpol$.
Then the code 
$\gcode [\gamma]$
coincides with 
the code $\gcodeXpol$.
\end{corollary}
By definition, we have
$$
\gcode[\gamma\circ\sigma]=\sigma\inv (\gcode [\gamma])
\qquad\hbox{for any $\sigma\in \SSSS (\Ps)$.}
$$
For each code $\abscode \in [\abscode]$,
we put
$$
\gs\sb{\abscode}:=\set{\gamma\in \gs}{\gcode [\gamma]=\abscode}.
$$
Then we have
$$
\gs\sb{\abscode}\sp{\sigma}=\gs\sb{\sigma\inv (\abscode)},
$$
where $\gs\sb{\abscode}\sp{\sigma}$ denotes the image of $\gs\sb{\abscode}$ by the action of $\sigma\in \SSSS(\Ps)$.
Therefore we obtain 
$$
\gs\sb{[\abscode]}=
\bigsqcup\sb{\abscode\sprime\in [\abscode]}\; \gs\sb{\abscode\sprime}=
\bigsqcup\sb{\sigma} \;\gs\sb{\abscode}\sp{\sigma}
\qquad \textrm{(disjoint union)},
$$
where $\sigma$ runs through the set of representatives for 
the right cosets in $\SSSS (\Ps)$ 
with respect to the subgroup $\Aut (\abscode)\st \SSSS(\Ps)$.
Hence we have
\begin{equation*}\label{eq:modulicode}
\moduli\sb{[\abscode]}=\PGL (3, k) \backslash  \gs\sb{\abscode}/\Aut(\abscode).
\end{equation*}
For $\gamma\in \gs\sb{\abscode}$,
let $[\gamma]\in \PGL (3, k)\backslash \gs\sb{\abscode}$ 
denote the projective equivalence class of $\gamma$.
From Corollary~\ref{cor:projaut}, we obtain the following:
\begin{corollary}\label{cor:projaut2}
Let $\Xpol$ be a polarized supersingular $K3$ surface of type $(\sharp)$
corresponding to the image of $[\gamma]\in \PGL (3, k)\backslash\gs\sb{\abscode}$
by the quotient map
$$
\PGL (3, k)\backslash\gs\sb{\abscode}\;\;\to\;\;
\moduli\sb{[\abscode]}=\PGL (3, k)\backslash\gs\sb{\abscode}/\Aut(\abscode).
$$
Via  the natural embedding of $\Aut\Xpol$
into 
$\Aut(\gcodeXpol)=\Aut(\abscode)$,
the automorphism group $\Aut\Xpol$ is equal to the stabilizer subgroup
of the point  $[\gamma]$.
\end{corollary}
\section{The isomorphism classes of codes with  Artin invariant $1$ and $2$}\label{sec:1and2}
We have classified all isomorphism classes of codes satisfying the conditions in (ii) of Theorem~\ref{thm:codes}.
The list is given in~\cite[\S8]{Shimada2004}.
Using the classification,
we have obtained the following~\cite[Corollary 1.11]{Shimada2004}:
\begin{theorem}\label{thm:DK}
There exists exactly one isomorphism class $[\abscode\sb{0}]$ of codes of dimension $10$
satisfying  the conditions in {\rm (ii)} of Theorem~\ref{thm:codes}.
The moduli space $\moduli\sb{[\abscode\sb{0}]}$ consists of a single point
corresponding to the Dolgachev-Kondo polynomial
$$
G\sb{\DK}:= XYZ(X^3+Y^3+Z^3).
$$
\end{theorem}
We call the point $[G\sb{\DK}]$  constituting $\moduli\sb{1}=\moduli\sb{[\abscode\sb 0]}$ the \emph{Dolgachev-Kondo point}.
We define the \emph{Dolgachev-Kondo code} 
$$
\abscode\sb{\DK}\;\;\st\;\; \Pow (\Pt (\F\sb 4))
$$
to be  the code generated by the words $\Lambda (\F\sb 4)$,
where $\Lambda$ are $\F\sb 4$-rational lines in $\Pt$.
The codes in the isomorphism class $[\abscode\sb 0]$
are precisely the codes $\gamma\inv  (\abscode\sb{\DK})$,
where 
$\gamma$ runs through the set of 
all bijections  from $\Ps$ to $\Pt (\F\sb 4)=\ZZ{dG\sb{\DK}}$.
The weight enumerator of any  code  in $[\abscode\sb 0]$  is
$$
1 + 21 z^5 + 210 z^8 + 280 z^9 + 280 z^{12} + 210 z^{13} + 21 z^{16} + z^{21}.
$$
There are no quadratic  nor cubic words in $\abscode\sb{0}$.
\par
\medskip
From the list in~\cite[\S8]{Shimada2004},
we obtain the following:
\begin{proposition}
There are exactly three isomorphism classes $[\CA]$, $[\CB]$, $[\CC]$
of codes of dimension $9$
satisfying the conditions in {\rm (ii)} of Theorem~\ref{thm:codes}.
\end{proposition}
As representatives of these isomorphism classes,
we can take codes $\CA$, $\CB$ and $\CC$
generated by vectors in Tables~\ref{table:genCA},~\ref{table:genCB} and~\ref{table:genCC}.
\begin{table}[t]
$$
\begin{array}{ccccccccccccccccccccccc}
\lb &1&1&1&1&1&1&1&1&1&1&1&1&1&1&1&1&1&1&1&1&1&\rb\\ 
\lb &0&0&0&0&0&0&0&0&0&0&0&0&0&0&0&0&1&1&1&1&1&\rb\\ 
\lb &0&0&0&0&0&0&0&0&0&0&0&0&1&1&1&1&0&0&0&0&1&\rb\\ 
\lb &0&0&0&0&0&0&0&0&0&1&1&1&0&0&0&1&0&0&0&1&0&\rb\\ 
\lb &0&0&0&0&0&0&0&1&1&0&0&1&0&0&1&0&0&0&1&0&0&\rb\\ 
\lb &0&0&0&0&0&0&1&0&1&0&1&0&0&1&0&0&0&1&0&0&0&\rb\\ 
\lb &0&0&0&0&1&1&0&0&0&0&0&1&0&1&0&0&0&1&1&1&1&\rb\\ 
\lb &0&0&0&1&0&1&0&0&0&0&1&0&0&1&1&1&0&0&1&0&1&\rb\\ 
\lb &0&0&1&0&0&1&0&0&1&0&0&1&0&0&1&1&0&1&1&0&0&\rb
\end{array}
$$
\caption{Generators of the code $\CA$}\label{table:genCA}
$$
\begin{array}{ccccccccccccccccccccccc}
\lb &1& 1& 1& 1& 1& 1& 1& 1& 1& 1& 1& 1& 1& 1& 1& 1& 1& 1& 1& 1& 1 &\rb\\   
\lb &0& 0& 0& 0& 0& 0& 0& 0& 0& 0& 0& 0& 0& 0& 0& 0& 1& 1& 1& 1& 1 &\rb\\   
\lb &0& 0& 0& 0& 0& 0& 0& 0& 0& 0& 0& 0& 1& 1& 1& 1& 0& 0& 0& 0& 1 &\rb\\   
\lb &0& 0& 0& 0& 0& 0& 0& 0& 0& 1& 1& 1& 0& 0& 0& 1& 0& 0& 0& 1& 0 &\rb\\     
\lb &0& 0& 0& 0& 0& 0& 0& 1& 1& 0& 0& 1& 0& 0& 1& 0& 0& 0& 1& 0& 0 &\rb\\    
\lb &0& 0& 0& 0& 0& 1& 1& 0& 0& 0& 0& 1& 0& 1& 0& 0& 0& 1& 0& 0& 0 &\rb\\     
\lb &0& 0& 0& 0& 1& 0& 1& 0& 1& 0& 1& 0& 0& 0& 0& 0& 0& 0& 0& 0& 1 &\rb\\   
\lb &0& 0& 1& 1& 0& 0& 0& 0& 0& 0& 1& 0& 0& 0& 1& 0& 0& 1& 0& 0& 0 &\rb\\   
\lb &0& 1& 0& 1& 0& 0& 0& 0& 1& 0& 0& 0& 0& 1& 0& 0& 0& 0& 0& 1& 0 &\rb
\end{array}
$$
\caption{Generators of the code $\CB$}\label{table:genCB}
$$
\begin{array}{ccccccccccccccccccccccc}
\lb& 1 &  1 &  1 &  1 &  1 &  1 &  1 &  1 &  1 &  1 &  1 &  1 &  1 &  1 &  1 &  1 &  1 &  1 &  1 &  1 &  1&\rb\\
\lb& 0 &  0 &  0 &  0 &  0 &  0 &  0 &  0 &  0 &  0 &  0 &  0 &  0 &  0 &  0 &  0 &  1 &  1 &  1 &  1 &  1&\rb\\ 
\lb& 0 &  0 &  0 &  0 &  0 &  0 &  0 &  0 &  0 &  0 &  0 &  0 &  1 &  1 &  1 &  1 &  0 &  0 &  0 &  0 &  1&\rb\\ 
\lb& 0 &  0 &  0 &  0 &  0 &  0 &  0 &  0 &  1 &  1 &  1 &  1 &  0 &  0 &  0 &  0 &  0 &  0 &  0 &  0 &  1&\rb\\
\lb& 0 &  0 &  0 &  0 &  0 &  0 &  1 &  1 &  0 &  0 &  1 &  1 &  0 &  0 &  1 &  1 &  0 &  0 &  1 &  1 &  0&\rb\\ 
\lb& 0 &  0 &  0 &  0 &  0 &  1 &  0 &  1 &  0 &  1 &  0 &  1 &  0 &  1 &  0 &  1 &  0 &  1 &  0 &  1 &  0&\rb\\ 
\lb& 0 &  0 &  0 &  0 &  1 &  0 &  0 &  1 &  0 &  1 &  1 &  0 &  0 &  1 &  1 &  0 &  0 &  1 &  1 &  0 &  1&\rb\\
\lb& 0 &  0 &  1 &  1 &  0 &  0 &  0 &  0 &  0 &  0 &  1 &  1 &  0 &  1 &  0 &  1 &  0 &  1 &  1 &  0 &  0&\rb\\
\lb& 0 &  1 &  0 &  1 &  0 &  0 &  0 &  0 &  0 &  1 &  0 &  1 &  0 &  1 &  1 &  0 &  0 &  0 &  1 &  1 &  0 &\rb
\end{array}
$$
\caption{Generators of the code $\CC$}\label{table:genCC}
\end{table}
The numbers  of linear, quadratic  and cubic words in these codes are given in the following table:
$$
\renewcommand{\arraystretch}{1.2}
\begin{array}{c | cccc}
 & \hbox{linear}  & \hbox{quadratic} & \hbox{cubic}&\\
\hline
\CA & 13 & 28 & 0& \\
\CB & 9 & 66 & 0& \\
\CC & 5 & 120 & 0 &.
\end{array}
$$
The weight enumerators of these codes are as follows:
\begin{eqnarray*}
\CA &:& 1+13 z^5 + 106 z^8 + 136 z^9 + 136 z^{12}+106 z^{13}+13 z^{16} +z^{21},\\
\CB &:& 1+9 z^5 + 102 z^8 + 144 z^9 + 144 z^{12}+102 z^{13}+9 z^{16} +z^{21},\\
\CC &:& 1+5 z^5 + 130 z^8 + 120 z^9 + 120 z^{12}+130 z^{13}+5 z^{16} +z^{21}.
\end{eqnarray*}
%
%
%
%
\begin{remark}
The Dolgachev-Kondo code $\abscode\sb{\DK}$
is related to the binary Golay code $\abscode_{24}$ in the following way.
Let $M:=\{\mu_1, \dots, \mu_{24}\}$
be the set of positions of the Miracle Octad Generator (MOG)
as is indicated in \cite[Table~6.1]{Shimada2001}.
The definition of $\abscode_{24}$ as a subcode of $\Pow (M)$ is described in \cite[Chapter~11]{CS}.
We put $N:=\{\mu_{22}, \mu_{23}, \mu_{24}\}\subset M$, and consider the $10$-dimensional subcode
$$
\abscode_{22}:=\set{w\in \abscode_{24}}{w\supset N\;\;\text{or}\;\;w\cap N=\emptyset}
$$
of $\abscode_{24}$.
We then define a map
$$
\Pt (\Ff)\;\;\longrightarrow\;\; M
$$
by \cite[Table~6.2]{Shimada2001}.
The pull-back of $\abscode_{22}$ by this map is just the Dolgachev-Kondo code $\abscode\sb{\DK}$.
\end{remark}
\begin{remark}
The codes $\CA$, $\CB$ and $\CC$ are isomorphic to linear subcodes of $\abscode\sb{\DK}$
defined as follows.
Let $F=\{Q_1, Q_2, Q_3, Q_4\}$ be a set of four points of $\Pt(\Ff)$,
and let $\abscode_{F}$ be the $9$-dimensional
linear subcode of  $\abscode\sb{\DK}$ defined by
$$
\abscode_{F}:=\set{w\in \abscode_{\DK}}{ \hbox{$|w\cap F|$ is even}}.
$$
If no three points of $F$ are collinear,
then $\abscode_{F}$ is isomorphic to $\CA$;
if exactly one triplet of the points of $F$ are collinear,
then $\abscode_{F}$ is isomorphic to $\CB$;
while if $F$ is on a line, then $\abscode_{F}$ is isomorphic to $\CC$.
\end{remark}
%
%
%
%
For $T=A, B$ and $C$,
we will write  $\moduli\sb{T}$ instead of $\moduli\sb{[\CT]}$,
and $\gs\sb{T}$ instead of $\gs\sb{\CT}$.
In the next three  sections, 
we will construct explicitly the space
$$
\moduli\sb{T}=\PGL (3, k) \backslash  \gs\sb{T}/\Aut(\CT)
$$
 for $T=A, B, C$,
and prove Theorems~\ref{thm:A},~\ref{thm:B} and~\ref{thm:C} stated in Introduction.
For this purpose,
we have to determine the group $\Aut(\abscode\sb T)$ and the space $\gs\sb{T}$.
Since $\CT$ is generated by $\Ps$ and the set of linear and quadratic words,
we obtain the following:
\begin{proposition}\label{prop:autcodeW}
Let $W\sb 1 (\CT)$ and $W\sb 2 (\CT)$ be the sets of linear and quadratic words
in $\CT$, respectively.
An element $\sigma$ of $\SSSS (\Ps)$ is contained in $\Aut (\CT)$
if and only if the following hold:
$$
\sigma(W\sb 1(\CT))=W\sb 1(\CT) \quand 
\sigma(W\sb 2(\CT))=W\sb 2(\CT).
$$
\end{proposition}
\begin{proposition}\label{prop:ings}
Suppose that a  map
$\gamma :\Ps\to \Pt$ is given.
Then $\gamma$ is contained in $ \gs\sb{T}=\shortset{\gamma\in\gs}{\CCC[\gamma]=\CT}$ if and only if the following hold:
\begin{itemize}
\item[(i)] $\gamma$ is injective, 
\item[(ii)] there exists a homogeneous polynomial $G$ of degree $6$ such that $\gamma (\Ps)=\ZZ{dG}$,
\item[(iii)] for every linear word $l$ of $\CT$,
there exists a line $L\st\Pt$ containing $\gamma (l)$, and 
\item[(iv)] for every quadratic word $q$ of $\CT$,
there exists a nonsingular conic curve  $Q\st\Pt$ containing $\gamma (q)$.
\end{itemize}
\end{proposition}
\begin{proof}
The ``only if " part is obvious from the definition of $\gs\sb{T}$.
Suppose that $\gamma$ satisfies (i)-(iv).
By (i) and (ii),
we have $\gamma\in \gs$.
Since $\CT$ is generated by the word  $\Ps$ and linear and quadratic words,
the properties (iii) and (iv) implies that
$\CT\subseteq \CCC[\gamma]$.
If $\CT\ne \CCC[\gamma]$,
then, by Theorem~\ref{thm:DK},
the code $\CCC[\gamma]\st\Pow(\Ps)$ is isomorphic to 
the  Dolgachev-Kondo code $\abscode\sb{\DK}\st\Pow(\Pt (\Ff))$ by some bijection from $\Ps$ to $\Pt (\Ff)$.
Hence  there exists $g\in \PGL(3, k)$ such that
$$
g(\gamma(\Ps))=\ZZ{dG\sb{\DK}}=\Pt (\Ff).
$$
However, there are no eight  points in $\Pt (\Ff)$ that are contained in a nonsingular conic curve.
\end{proof}
It will turn out that,
for $T=A,B$ and $C$, the following hold.
\par
\medskip
(1) The space $\PGL (3, k)\backslash \gs\sb T$ has exactly two connected components,
both of which are isomorphic to $\A\sp 1\sm\{0,1,\omega,\bar\omega\}$.
Let $N\sb T\st \Aut(\CT)$ be the subgroup consisting of the elements that do not interchange the two
connected components,
and let $\Gamma\sb T$ be the image of $N\sb T$ in $\Aut (\A\sp 1\sm \{0,1,\omega,\bar\omega\})$.
Then $N\sb T$ is of index $2$
in $\Aut(\CT)$.  
The moduli curve $\moduli\sb T$ is the quotient of $\A\sp 1\sm \{0,1,\omega,\bar\omega\}$ by $\Gamma\sb T$.
\par
\medskip
(2) 
The action of $\Gamma\sb T$ on the punctured affine line $\A\sp 1\sm \{0,1,\omega,\bar\omega\}$
is free.
Hence 
the order of the stabilizer subgroup
$\Stab([\gamma])\st\Aut(\CT)$ 
of a point $[\gamma]\in \PGL (3, k)\backslash \gs\sb T$
is constant on $\PGL (3, k)\backslash \gs\sb T$.
By Corollary~\ref{cor:projaut2}, $\Stab([\gamma])$ is equal to $\Aut\Xpol$,
where $\Xpol$ corresponds to the image of $[\gamma]$ in $\moduli\sb T$.
Hence 
we have an exact sequence
$$
1\;\;\to\;\; \Aut\Xpol\;\; \to\;\; N\sb T\;\;\to\;\; \Gamma\sb T\;\;\to\;\;1
$$
for any polarized supersingular $K3$ surface $\Xpol$ corresponding to a point of $\moduli\sb T$.
\par
\medskip
The orders of the groups above are given as follows.
$$
\renewcommand{\arraystretch}{1.2}
\begin{array}{|c|ccccccc|}
\hline
T & |\Aut(\CT)| &=& 2 &\times& | \Gamma\sb T| &\times& |\Aut\Xpol| \\
\hline
A & 1152 &=& 2  &\times& 6  &\times& 96\\
B & 432 &=& 2  &\times& 12  &\times&18\\
C & 23040 &=& 2  &\times& 12  &\times& 960\\
\hline
\end{array}
$$
\begin{remark}\label{rem:conicalg}
The following algorithm will be  used frequently.
Suppose that we are given eight points
$$
p\sb i=[\xi\sb i, \eta\sb i, \zeta\sb i]\qquad(i=1, \dots, 8)
$$
on $\Pt$. 
In order for them to be on a (possibly singular) conic  curve,
it is necessary and sufficient that the $8\times 6$ matrix 
$$
M:=
\left[
		\vcenter{
				\vbox{
						\hbox to 5cm{$
								\; \xi\sb 1^2,\; \eta\sb 1^2, \; \zeta\sb 1^2, \; \xi\sb 1\eta\sb 1, \; \eta\sb 1\zeta\sb 1, \;
								\zeta\sb 1 \xi\sb 1\;
						$}
						\hbox to 5cm{
								\hfil$\cdots\cdots$\hfil
						}
					\hbox to 5cm{
								\hfil$\cdots\cdots$\hfil
						}
						\hbox to 5cm{$
								\; \xi\sb 8^2,\; \eta\sb 8^2, \; \zeta\sb 8^2, \; \xi\sb 8\eta\sb 8, \; \eta\sb 8\zeta\sb 8, \;
								\zeta\sb 8 \xi\sb 8\;
						$}
				}
		}
\right]
$$
is of rank $<6$.
When  the rank of $M$ is $<6$,
a non-zero solution
$$
{}\sp{T} [A, B, C, D, E, F]
$$
of the linear equation $M\vx=0$ gives us a defining equation
\begin{equation}\label{eq:genconic}
AX^2+BY^2+CZ^2+DXY+EYZ+FZX=0
\end{equation}
of a conic curve containing $p\sb 1, \dots, p\sb 8$.
\end{remark}
The following are
 phenomena peculiar to projective geometry in characteristic $2$.
\begin{remark}\label{rem:nonsingconic}
The conic  curve defined by the equation~\eqref{eq:genconic} is singular
if and only if the following holds:
$$
AE^2+BF^2+CD^2 + DEF=0.
$$
\end{remark}
\begin{definition}\label{def:tangent}
Let $L\st \Pt$ be a line,
and let $Q\st\Pt$ be a (possibly singular) conic curve.
We say that $L$ and $Q$ are \emph{tangent}
if they fail to intersect at distinct two points.
\end{definition}
\begin{remark}\label{rem:tangent}
Let $L$ be a line.
Then the conic curves tangent to $L$ form a \emph{linear} system in $|\OOO\sb{\Pt} (2)|$.
If three distinct lines $L\sb 1$, $L\sb 2$ and $L\sb 3$ are concurrent,
then every conic curve that is tangent to $L\sb 1$ and $L\sb 2$ is tangent to $L\sb 3$.
\end{remark}
\begin{remark}\label{rem:tetrad}
Let $A, B, C, D\in \Pt$ be distinct points.
Suppose that no three of them are collinear.
Let 
$O$
(resp. $P$)
(resp. $Q$)
be the intersection point
of the lines $\overline{AB}$ and $\overline{CD}$
(resp. $\overline{AC}$ and $\overline{BD}$)
(resp. $\overline{AD}$ and $\overline{BC}$).
Then $O$, $P$ and $Q$ are collinear.
\end{remark}
\section{The moduli curve corresponding to the code $\CA$}\label{sec:A}
In this section,
we prove Theorem~\ref{thm:A}.
\par
\medskip
The linear words of $\CA$ are listed in Table~\ref{table:linesA}.
\begin{table}[b]
$$
\begin{array}{ccccccccc}
m & : & \{ &1, &2, &12, &13, &18&\},   \\ 
l\sb{12, 1} & : & \{ &10, &11, &12, &16, &20&\},  \\ 
l\sb{12, 2} & : & \{ &8, &9, &12, &15, &19&\},  \\ 
l\sb{12, 3} & : & \{ &5, &6, &12, &14, &17&\},  \\ 
l\sb{12, 4} & : & \{ &3, &4, &7, &12, &21&\},  \\ 
l\sb{13, 1} & : & \{ &13, &14, &15, &16, &21&\},  \\ 
l\sb{13, 2} & : & \{ &7, &8, &10, &13, &17&\},  \\ 
l\sb{13, 3} & : & \{ &4, &6, &11, &13, &19&\},  \\ 
l\sb{13, 4} & : & \{ &3, &5, &9, &13, &20&\},  \\ 
l\sb{18, 1} & : & \{ &17, &18, &19, &20, &21&\},  \\ 
l\sb{18, 2} & : & \{ &7, &9, &11, &14, &18&\},  \\ 
l\sb{18, 3} & : & \{ &4, &5, &10, &15, &18&\},  \\ 
l\sb{18, 4} & : & \{ &3, &6, &8, &16, &18&\}.  \\
\phantom{l\sb{18, 4}} & :& &&&&&&
\end{array}
$$
\caption{Linear words in $\CA$}\label{table:linesA}
\end{table}
From now on, we sometimes 
abbreviate, for example,  the set $\{ P\sb 8, P\sb 9, P\sb{12}, P\sb{15}, P\sb{19}\}$
to $\{ 8, 9, 12,15,19\}$.
The linear word $m$ stands out from the rest in that there are two points
$P\sb 1$ and $P\sb 2$ in $m$
through which no other linear words pass.
We call $m$ the \emph{special linear  word}.
The other linear words are divided into three groups
according to the intersection point with $m$.
For $\nu=12, 13, 18$ and $i=1,2,3,4$,
the non-special linear word  $l\sb{\nu, i}$ intersects $m$ at the  point $P\sb\nu$.
For each of $P\sb 1$ and $P\sb 2$,
there exists only one linear word  $m$ containing it.
For each of  $P\sb {12}$, $P\sb{13}$ and $P\sb {18}$,
there exist exactly five linear words containing it.
For each of the other $16$ points,
there exist exactly three linear words containing it.
For each $\alpha, \beta=1, \dots, 4$, there exists a unique
$\gamma=\gamma(\alpha, \beta)$ such that the three linear words 
$l\sb{12, \alpha}$,
$l\sb{13, \beta}$ and 
$l\sb{18, \gamma}$
have a point  in common.
\begin{table}[t]
\begin{tabular}{c|cccc}
$\beta\;\; \backslash\;\; \alpha$ & 1 & 2 & 3 & 4 \\
\hline
1 & 4 & 3 & 2 & 1 \\
2 & 3 & 4 & 1 & 2 \\
3 & 2 & 1 & 4 & 3 \\
4 & 1 & 2 & 3 & 4 \\
\end{tabular}
\vskip 3pt
\caption{Concurrent triples $(\alpha, \beta, \gamma (\alpha, \beta))$}\label{table:concurrent_triples}
\end{table}
We call such a triple $(\alpha, \beta, \gamma)$ a \emph{concurrent triple}.
The list of concurrent triples is given in Table~\ref{table:concurrent_triples}.
For a concurrent triple $(\alpha,\beta,\gamma)$,
we denote by $T\sb{\alpha\beta\gamma}$
the intersection point of  $l\sb{12, \alpha}$,
$l\sb{13, \beta}$ and 
$l\sb{18, \gamma}$.
\begin{table}[t]
$$
\renewcommand{\arraystretch}{1.1}
\begin{array}{|c|cccccccc|}
\hline 
\alpha \beta \gamma &114 &123 &132 &141 &213 &224 &231 &242 \\ 
T\sb{\alpha\beta\gamma} &P\sb{16}&P\sb{10}&P\sb{11}&P\sb{20}&P\sb{15}&P\sb{8}&P\sb{19}&P\sb{9}\phantom{\vrule depth 6pt}\\ 
\hline 
\alpha \beta \gamma &312 &321 &334 &343 &411 &422 &433 &444 \\ 
T\sb{\alpha\beta\gamma} &P\sb{14}&P\sb{17}&P\sb{6}&P\sb{5}&P\sb{21}&P\sb{7}&P\sb{4}&P\sb{3}\phantom{\vrule depth 6pt}\\ 
\hline
\end{array}
$$
\vskip 3pt
\caption{Points $T\sb{\alpha\beta\gamma}$}\label{table:center}
\end{table}
\par
\medskip
The $28$  quadratic words in $\CA$
 are divided into two groups.
The quadratic words $q\sb 1\sprime, \dots, q\sb{12}\sprime$ listed in Table~\ref{table:qsprimeA}
are disjoint from the special  linear word  $m$, and intersect each of 
the non-special linear words $l\sb{\nu, i}$ at distinct two points.
On the other hand,
for each concurrent triple  $(\alpha, \beta, \gamma)$,
there exists a unique quadratic word $q\sb{\alpha\beta\gamma}$
that is disjoint from  the three linear words 
$l\sb{12, \alpha}$, $l\sb{13, \beta}$, $l\sb{18, \gamma}$,
and intersects  other ten linear words at distinct two points.
The list of these quadratic words $q\sb{\alpha\beta\gamma}$ is given in Table~\ref{table:qabcA}.
\begin{table}[t]
$$
\begin{array}{cccccccccccc}
q\sprime\sb{1} & : & \{ &5, &6, &7, &9, &10, &16, &19, &21&\},  \\ 
q\sprime\sb{2} & : & \{ &5, &6, &7, &8, &11, &15, &20, &21&\},  \\ 
q\sprime\sb{3} & : & \{ &4, &6, &8, &9, &10, &14, &20, &21&\},  \\ 
q\sprime\sb{4} & : & \{ &4, &6, &7, &9, &15, &16, &17, &20&\},  \\ 
q\sprime\sb{5} & : & \{ &4, &5, &8, &9, &11, &16, &17, &21&\},  \\ 
q\sprime\sb{6} & : & \{ &4, &5, &7, &8, &14, &16, &19, &20&\},  \\ 
q\sprime\sb{7} & : & \{ &3, &6, &9, &10, &11, &15, &17, &21&\},  \\ 
q\sprime\sb{8} & : & \{ &3, &6, &7, &10, &14, &15, &19, &20&\},  \\ 
q\sprime\sb{9} & : & \{ &3, &5, &8, &10, &11, &14, &19, &21&\},  \\ 
q\sprime\sb{10} & : & \{ &3, &5, &7, &11, &15, &16, &17, &19&\},  \\ 
q\sprime\sb{11} & : & \{ &3, &4, &9, &10, &14, &16, &17, &19&\},  \\ 
q\sprime\sb{12} & : & \{ &3, &4, &8, &11, &14, &15, &17, &20&\}.  \\ 
\phantom{q\sprime\sb{12}} &  & &\phantom{20} &\phantom{20} &\phantom{20} &\phantom{20} &\phantom{20} &\phantom{20} &\phantom{20}
&\phantom{20}& 
\end{array}
$$
\caption{Quadratic words $q\sprime\sb\nu$ in $\CA$}\label{table:qsprimeA}
\end{table}
\begin{table}[t]
$$
\begin{array}{cccccccccccc}
q\sb{114} & : & \{ &1, &2, &4, &5, &7, &9, &17, &19&\},  \\ 
q\sb{123} & : & \{ &1, &2, &3, &6, &9, &14, &19, &21&\},  \\ 
q\sb{132} & : & \{ &1, &2, &3, &5, &8, &15, &17, &21&\},  \\ 
q\sb{141} & : & \{ &1, &2, &4, &6, &7, &8, &14, &15&\},  \\ 
q\sb{213} & : & \{ &1, &2, &3, &6, &7, &11, &17, &20&\},  \\ 
q\sb{224} & : & \{ &1, &2, &4, &5, &11, &14, &20, &21&\},  \\ 
q\sb{231} & : & \{ &1, &2, &3, &5, &7, &10, &14, &16&\},  \\ 
q\sb{242} & : & \{ &1, &2, &4, &6, &10, &16, &17, &21&\},  \\ 
q\sb{312} & : & \{ &1, &2, &3, &4, &8, &10, &19, &20&\},  \\ 
q\sb{321} & : & \{ &1, &2, &3, &4, &9, &11, &15, &16&\},  \\ 
q\sb{334} & : & \{ &1, &2, &7, &9, &10, &15, &20, &21&\},  \\ 
q\sb{343} & : & \{ &1, &2, &7, &8, &11, &16, &19, &21&\},  \\ 
q\sb{411} & : & \{ &1, &2, &5, &6, &8, &9, &10, &11&\},  \\ 
q\sb{422} & : & \{ &1, &2, &5, &6, &15, &16, &19, &20&\},  \\ 
q\sb{433} & : & \{ &1, &2, &8, &9, &14, &16, &17, &20&\},  \\ 
q\sb{444} & : & \{ &1, &2, &10, &11, &14, &15, &17, &19&\}.  \\ 
\phantom{q\sb{444}} &  & &\phantom{20} &\phantom{20} &\phantom{20} &\phantom{20} &\phantom{20} &\phantom{20} &\phantom{20}
&\phantom{20}& 
\end{array}
$$
\caption{Quadratic words $q\sb{\alpha\beta\gamma}$ in $\CA$}\label{table:qabcA}
\end{table}
\par
\medskip
In order to study $\Aut (\CA)$,
we embed $\CA$ into the Dolgachev-Kondo code
$\CDK\subset \Pow (\Pt(\Ff))$
by the bijection
$\phi : \Ps \isom \Pt(\F\sb 4)$ given in Table~\ref{table:phi}.
\begin{table}[t]
\par
\medskip
\noindent
{\phantom{a}
\hskip -6.5cm 
\hbox {
\hbox{
\hbox to 5cm{
\parbox{5cm}{
\begin{eqnarray*}
\phi(P\sb{1})&=& [1,\omega,0], \\
\phi(P\sb{2})&=& [1,\bar\omega,0], \\
\phi(P\sb{3})&=& [1,1,\omega], \\
\phi(P\sb{4})&=& [1,\bar\omega,\omega], \\
\phi(P\sb{5})&=& [1,1,\bar\omega], \\
\phi(P\sb{6})&=& [1,\omega,\bar\omega], \\
\phi(P\sb{7})&=& [1,\omega,\omega], \\
\phi(P\sb{8})&=& [1,\bar\omega,1], \\
\phi(P\sb{9})&=& [1,1,1], \\
\phi(P\sb{10})&=& [0,1,\bar\omega], \\
\phi(P\sb{11})&=& [0,1,1], 
\end{eqnarray*}
}
}
\hbox to 6cm{
\parbox{6cm}{
\begin{eqnarray*}
\phi(P\sb{12})&=& [0,1,0], \\
\phi(P\sb{13})&=& [1,1,0], \\
\phi(P\sb{14})&=& [1,\bar\omega,\bar\omega], \\
\phi(P\sb{15})&=& [1,\omega,1], \\
\phi(P\sb{16})&=& [0,1,\omega], \\
\phi(P\sb{17})&=& [1,0,\bar\omega], \\
\phi(P\sb{18})&=& [1,0,0], \\
\phi(P\sb{19})&=& [1,0,1], \\
\phi(P\sb{20})&=& [0,0,1], \\
\phi(P\sb{21})&=& [1,0,\omega]. \\
\phantom{\phi(P\sb{21})}&&
\end{eqnarray*}
}
}
}
}
}
\vskip 5pt
\caption{Bijection $\phi$ from $\Ps$ to $\Pt (\F\sb 4)$}\label{table:phi}
\end{table}
%
%
%
%
The following can be checked easily.
\begin{itemize}
\item[(1)]
If $l$ is a linear word of $\CA$, then the points in $\phi(l)$ are collinear.
The  linear words of $\CA$ coincide with  $\phi\inv (\Lambda (\Ff))$,
where $\Lambda $ are  $\Ff$-rational lines containing  at least one of
$\phi (\Psb{12})$, $\phi (\Psb{13})$, $\phi (\Psb{18})$.
\item[(2)]
The words  $q\sprime\sb 1, \dots, q\sprime\sb{12}$
coincide with the words written as
$$
\phi\inv (\Lambda\sb 1 (\Ff) + \Lambda\sb 2 (\Ff)),
$$
where $\Lambda\sb 1$ and $\Lambda\sb 2$ are distinct $\Ff$-rational lines 
such that both of $\Lambda\sb 1 (\Ff)$ and $\Lambda\sb 2 (\Ff)$
are disjoint from $\{ \phi (\Psb{12}), \phi (\Psb{13}), \phi (\Psb{18}) \}$,
and such that 
the  intersection point of $\Lambda\sb 1 (\Ff)$ and $\Lambda\sb 2 (\Ff)$ 
is either  $\phi (\Psb{1})$ or $\phi(\Psb{2})$.
\item[(3)]
For a concurrent triple $(\alpha, \beta, \gamma)$,
let 
$\Lambda\sb i$ be the $\Ff$-rational line passing through $\phi (T\sb{\alpha\beta\gamma})$ and $\phi (P\sb i)$
for $i=1, 2$.
Then we have
$q\sb{\alpha\beta\gamma}=\phi\inv  ( \Lambda\sb 1(\Ff) + \Lambda\sb 2 (\Ff) )$.
\end{itemize}
Let $\PG\sprime$ be the subgroup of $\PGL (3, \Ff)$
consisting of $g\in \PGL (3, \Ff)$
satisfying
$$
\{g(\phi (\Psb{12})), g(\phi(\Psb{13})), g(\phi(\Psb{18}))\}=\{ \phi (\Psb{12}), \phi(\Psb{13}), \phi(\Psb{18})\},
$$
and let $\PG$ be the subgroup $\phi\inv \circ \PG\sprime\circ \phi$ of $\SSSS (\Ps)$.
The order of $\PG$ is $288$.
Let $F\sprime \in \SSSS(\Pt(\Ff))$ be the element of order $2$ 
obtained by the conjugation $\omega\mapsto\bar\omega$ of $\Ff$ over $\F\sb 2$.
We then put
$$
F:=\phi\inv \circ F\sprime\circ \phi = (\Psb{1}\Psb{2})(\Psb{3}\Psb{5})(\Psb{4}\Psb{6})(\Psb{7}\Psb{14})(\Psb{8}\Psb{15})(\Psb{10}\Psb{16})(
\Psb{17}\Psb{21}) \in
\SSSS (\Ps).
$$
We also
put
$$
T:=(\Psb{1}\Psb{2}).
$$
\begin{proposition}
The group  $\Aut(\CA)$ is  of order $1152$,
and is 
generated by $\PG$, $F$ and $T$.
\end{proposition}
\begin{proof}
Since the actions of $\PG\sprime$ and $F\sprime$ on $\Pt(\Ff)$
leave the set 
$$
\{[0,1,0],[1,1,0],[1,0,0]\}=\{ \phi (\Psb{12}), \phi(\Psb{13}), \phi(\Psb{18}) \}
$$
invariant,
and 
preserve the line-point incidence configuration,
we see that 
$PG\subset \Aut(\CA)$ and $F\in \Aut(\CA)$.
It is obvious that $T\in \Aut (\CA)$.
By direct calculation,
we see that the subgroup of $\SSSS (\Ps)$
generated by  $\PG$, $F$ and  $T$ is  of order $1152$.
\par
Every automorphism of $\CA$ leaves each of the sets
$\{\Psb{1}, \Psb{2}\}$ and $\{\Psb{12}, \Psb{13}, \Psb{18}\}$
invariant.
Hence we have a homomorphism
\begin{equation}\label{eq:ASS}
\Aut (\CA) \;\to\;\SSSS(\{\Psb{1}, \Psb{2}\})\times\SSSS(\{\Psb{12}, \Psb{13}, \Psb{18}\}).
\end{equation}
Since $\PG$ acts on $\{\Psb{12}, \Psb{13}, \Psb{18}\}$
as the full-symmetric group,
and since $T$ is contained in $\Aut(\CA)$,
the homomorphism~\eqref{eq:ASS} is surjective.
Let $K$ denote the kernel of~\eqref{eq:ASS}.
We have a homomorphism 
\begin{equation}\label{eq:KSS}
K \;\to\;\SSSS\sb 4\times\SSSS\sb 4, 
\qquad
g\;\mapsto\; (\sigma, \sigma\sprime),
\end{equation}
where $\sigma$ and $\sigma\sprime$ are given by
$$
g(l\sb{12, \alpha})=l\sb{12, \sigma(\alpha)},
\quad
g(l\sb{13, \beta})=l\sb{13, \sigma\sprime(\beta)}.
$$
We also have a homomorphism
\begin{equation}\label{eq:SSS}
\SSSS\sb 4\times\SSSS\sb 4\;\to\; \SSSS(\Ps), 
\qquad
(\sigma, \sigma\sprime)\;\mapsto\;g\sb{\sigma, \sigma\sprime},
\end{equation}
where $g\sb{\sigma, \sigma\sprime} $ is given by
\begin{eqnarray*}
g\sb{\sigma, \sigma\sprime} (P\sb i) &=& P\sb i\qquad \hbox{if $P\sb i\in m$}, \\
g\sb{\sigma, \sigma\sprime} (T\sb{\alpha\beta\gamma}) &=& T\sb{\sigma(\alpha)\sigma\sprime(\beta)\gamma\sprime}\qquad \\
&& \hbox{where $(\alpha, \beta, \gamma)$ and $(\sigma(\alpha), \sigma\sprime(\beta), \gamma\sprime)$ are concurrent triples.} \\
\end{eqnarray*}
Since the composite of~\eqref{eq:KSS} and~\eqref{eq:SSS}
is the identity of $K$,
the homomorphism~\eqref{eq:KSS} is injective.
For each pair $(\sigma, \sigma\sprime)$ of $\SSSS\sb 4\times\SSSS\sb 4$,
we check whether $g\sb{\sigma,\sigma\sprime}$ is in $\Aut(\CA)$;
that is,
whether $g\sb{\sigma,\sigma\sprime}$ satisfies the following (see Proposition~\ref{prop:autcodeW}):
\begin{equation}\label{eq:condAutCA}
g\sb{\sigma,\sigma\sprime} (W\sb 1(\CA))=W\sb 1(\CA)
\quand
g\sb{\sigma,\sigma\sprime} (W\sb 2(\CA))=W\sb 2(\CA).
\end{equation}
Among $(4!)^2=576$ pairs,
exactly $96$ pairs satisfy~\eqref{eq:condAutCA}.
Hence $\Aut(\CA)$ is of order $|K||\SSSS\sb 2||\SSSS\sb 3|=96\cdot 12 =1152$,
and is generated by $\PG$, $F$ and $T$.
\end{proof}
For a parameter $\lambda$ of the affine line $\A\sp 1$,
we define a map
$$
\map{\gamma\sb{\lambda}}{\Ps}{\Pt}
$$
by Table~\ref{table:gammalambdaA}.
Note that $\gamma\sb{\omega}$ coincides with $\phi$ defined above.
\begin{table}
\par
\medskip
\noindent
{\phantom{a}
\hbox{
\hbox{
\hskip -6.5cm
\hbox to 5cm{
\parbox{5cm}{
\begin{eqnarray*}
\gamma\sb{\lambda}(P\sb{1}) &=& [1,  \omega,0], \\ 
\gamma\sb{\lambda}(P\sb{2}) &=& [1,\bar\omega,0], \\ 
\gamma\sb{\lambda}(P\sb{3}) &=& [1+\lambda,1+\lambda,1], \\ 
\gamma\sb{\lambda}(P\sb{4}) &=& [1+\lambda,\lambda,1], \\ 
\gamma\sb{\lambda}(P\sb{5}) &=& [\lambda,\lambda,1], \\ 
\gamma\sb{\lambda}(P\sb{6}) &=& [\lambda,1+\lambda,1], \\ 
\gamma\sb{\lambda}(P\sb{7}) &=& [1+\lambda,1,1], \\ 
\gamma\sb{\lambda}(P\sb{8}) &=& [1,1+\lambda,1], \\ 
\gamma\sb{\lambda}(P\sb{9}) &=& [1,1,1], \\ 
\gamma\sb{\lambda}(P\sb{10}) &=& [0,\lambda,1], \\ 
\gamma\sb{\lambda}(P\sb{11}) &=& [0,1,1], 
\end{eqnarray*}
}
}
\hbox to 6cm{
\parbox{6cm}{
\begin{eqnarray*}
\gamma\sb{\lambda}(P\sb{12}) &=& [0,1,0], \\ 
\gamma\sb{\lambda}(P\sb{13}) &=& [1,1,0], \\ 
\gamma\sb{\lambda}(P\sb{14}) &=& [\lambda,1,1], \\ 
\gamma\sb{\lambda}(P\sb{15}) &=& [1,\lambda,1], \\ 
\gamma\sb{\lambda}(P\sb{16}) &=& [0,1+\lambda,1], \\ 
\gamma\sb{\lambda}(P\sb{17}) &=& [\lambda,0,1], \\ 
\gamma\sb{\lambda}(P\sb{18}) &=& [1,0,0], \\ 
\gamma\sb{\lambda}(P\sb{19}) &=& [1,0,1], \\ 
\gamma\sb{\lambda}(P\sb{20}) &=& [0,0,1], \\ 
\gamma\sb{\lambda}(P\sb{21}) &=& [1+\lambda,0,1]. \\
\phantom{\gamma\sb{\lambda}(P\sb{21})}&&
\end{eqnarray*}
}
}
}
}
}
\vskip 5pt
\caption{Definition of $\gamma\sb{\lambda}$ for $\CA$}
\label{table:gammalambdaA}
\end{table}
We denote by $\widetilde{T}$ the subgroup $\{ 1, T\}$
of $\Aut (\CA)$.
\begin{proposition}\label{prop:Agammalambda}
The map
$\lambda\mapsto \gamma\sb{\lambda}$
induces an isomorphism  from
$\A\sp 1 \sm \{ 0, 1, \omega, \bar{\omega} \}$
to $\PGL(3, k)\backslash \GGG\sb{A} /\widetilde{T}$.
\end{proposition}
\begin{proof}
First note that $\gamma\sb{\lambda}$ is injective if and only if
\begin{equation}\label{eq:injective}
\lambda\ne 0
\quand
\lambda\ne 1.
\end{equation}
From now on, we assume~\eqref{eq:injective}.

We will show the following:
\begin{claim}\label{claim:tripleA}
Let $\gamma\sprime$ be an arbitrary element of $\gs\sb{A}$.
Then 
 there exists a unique triple
$$
(g, t, \lambda)\;\in\; 
\PGL(3, k) \times \widetilde{T}\times (k\sm \{ 0, 1, \omega,\bar{\omega} \})
$$
such that
\begin{equation*}
g\circ \gamma\sprime \circ t = \gamma\sb{\lambda}.
\end{equation*}
\end{claim}
Because the points $\gamma\sprime (\Psb{18})$, $\gamma\sprime (\Psb{19})$, $\gamma\sprime (\Psb{20})$
of $\gamma\sprime(l\sb{18, 1})$
are on a line
and the points 
$\gamma\sprime (\Psb{12})$, $\gamma\sprime (\Psb{13})$, $\gamma\sprime (\Psb{18})$
of $\gamma\sprime(m)$
are on another line,
there exists a unique element
$g\in \PGL (3, k)$
such that
$\gamma:=g\circ\gamma\sprime$ satisfies the following:
\begin{eqnarray}\label{eq:gammaPFIX}
&& \gamma(P\sb{18})=[1,0,0], \nonumber\\
&&\gamma (P\sb{12})=[0,1,0],\;
\gamma(P\sb{13})=[1,1,0],\;\\
&&\gamma(P\sb{20})=[0,0,1],\;
\gamma(P\sb{19})=[1,0,1].
\nonumber
\end{eqnarray}
Let $L\sb{12,\alpha}$, $L\sb{13, \beta}$ and $L\sb{18, \gamma}$ be the lines containing
$\gamma(l\sb{12,\alpha})$, $\gamma(l\sb{13,\beta})$ and $\gamma(l\sb{18,\gamma})$, respectively.
We put $x:=X/Z$, $y:=Y/Z$.
Then the defining equations of these lines can be written as follows:
\begin{eqnarray}\label{eq:lineeq}
L\sb{12, \alpha}&:& x+a\sb{\alpha}=0, \nonumber \\
L\sb{13, \beta }&:& x+y+b\sb{\beta}=0, \\
L\sb{18, \gamma}&:& y+c\sb{\gamma}=0. \nonumber
\end{eqnarray}
From~\eqref{eq:gammaPFIX},
we have
\begin{equation}\label{eq:rig}
a\sb 1=0, \quad a\sb 2 =1,
\quad b\sb 3 =1, \quad b\sb 4=0,\quad
c\sb 1 =0.
\end{equation}
The condition that $(\alpha,\beta,\gamma)$ is a concurrent  triple
is equivalent to
$$
a\sb \alpha + b\sb \beta +c\sb \gamma=0.
$$
Solving the linear equations corresponding to  the $16$ concurrent triples
and combining the result with~\eqref{eq:rig},
we obtain the following solutions:
\begin{eqnarray}\label{eq:abcsol}
(a_1, a_2, a_3, a_4)&=&(0,1,\lambda,1+\lambda), \nonumber\\ 
(b_1, b_2, b_3, b_4)&=&(1+\lambda,\lambda,1,0),\\
(c_1, c_2, c_3, c_4)&=&(0,1,\lambda,1+\lambda), \nonumber
\end{eqnarray}
where $\lambda$ is a parameter.
The coordinates of the   points 
$T\sb{\alpha\beta\gamma}$
are given by
$[a\sb\alpha, c\sb\gamma, 1]$.
Using Table~\ref{table:center},
we see that
$\gamma (P\sb i)=\gamma\sb{\lambda} (P\sb i)$
holds for every $i$ except for $i=1$ and $i=2$.
The line $M$ containing $\gamma (m)$ is defined by $Z=0$.
Hence we can put
$$
\gamma (P\sb 1)=[1,\tau\sb 1, 0], \quad \gamma (P\sb 2)=[1,\tau\sb 2, 0].
$$
By the algorithm in Remark~\ref{rem:conicalg},
we see that 
a conic curve
containing 
$\gamma (q\sb{114})$ exists if and only if the following hold:
\begin{eqnarray}\label{eq:ccd114}
&&\left (1+\tau\sb{2}+{\tau\sb{2}}^{2}\right ) \left (\lambda+1\right )^{2}{\lambda}^{2}=0, \nonumber\\
&&\left (\tau\sb{1} +\tau\sb{2}\right ) \left (1 +\tau\sb{2} + {\tau\sb{2}}^{2}\right )
\left (\lambda+1\right)\lambda=0,\\
&&\left (\tau\sb{1}+\tau\sb{2}\right ) \left (\tau\sb{1}+\tau\sb{2}+1\right ) \left (\lambda+1\right )\lambda=0. \nonumber
\end{eqnarray}
Here we have used the Buchberger algorithm to calculate
the Gr\"obner basis of the ideal in $k[\lambda, \tau\sb 1, \tau\sb 2]$
generated by $6\times 6$-minors of the $8\times 6$-matrix corresponding to the eight points in 
$\gamma (q\sb{114})$.
Replacing $\gamma$ by $\gamma\circ T$ if necessary,
we have
$$
\tau\sb 1=\omega \quand \tau\sb 2=\bar\omega
$$
by~\eqref{eq:injective}, \eqref{eq:ccd114} and $\tau\sb 1\ne\tau\sb 2$.
Then  the 
conic curve containing 
$\gamma (q\sb{114})$
is defined by 
$$
X^2+Y^2+\lambda Z^2+XY+(\lambda+1)ZX=0,
$$
which is nonsingular
if and only if $\lambda^2+\lambda+1\ne 0$.
(See Remark~\ref{rem:nonsingconic}.)
Thus we have proved the existence and the uniqueness of the triple
$(g, t, \lambda)$ satisfying $g\circ \gamma\sprime \circ t = \gamma\sb{\lambda}$.
In particular,
for  each double coset in  $\PGL(3, k)\backslash \GGG\sb{A}/\widetilde{T}$,  
there exists a unique $\lambda\in k\sm \{ 0, 1, \omega,\bar{\omega} \}$
such that  $\gamma\sb{\lambda}$ is contained in the coset.
\par
Conversely, let $\lambda$ be an element of $k\sm\{0,1,\omega, \bar\omega\}$.
We will show that $\gamma\sb\lambda$ is in $\GGG\sb{A}$.
The points $\gamma\sb{\lambda} (\Ps)$
coincides with $\ZZ{ d\GA [\lambda] }$,
where $\GA[\lambda]$ is given in Theorem~\ref{thm:A}.
Indeed, we can check that
$$
\frac{\partial \,\GA [\lambda]}{\partial X} (\gamma\sb{\lambda} (P\sb i))=
\frac{\partial \,\GA [\lambda]}{\partial Y}  (\gamma\sb{\lambda} (P\sb i)) =
\frac{\partial \,\GA [\lambda]}{\partial Z}  (\gamma\sb{\lambda} (P\sb i))=
0
$$
holds for $i=1,\dots, 21$.
For each linear word $l$ of $\CA$,
there exists a line containing $\gamma\sb{\lambda} (l)$.
The defining equations of them are given by~\eqref{eq:lineeq} and~\eqref{eq:abcsol}. 
(The line $M$ containing $\gamma\sb{\lambda} (m)$ is defined by $Z=0$.)
For each quadratic word $q\sprime\sb i$  of $\CA$ (resp.~$q\sb{\alpha\beta\gamma}$),
there exists a nonsingular conic curve
$Q\sprime\sb i$ (resp.~$Q\sb{\alpha\beta\gamma}$)
containing 
$\gamma\sb{\lambda}(q\sprime\sb i)$ (resp. $\gamma\sb{\lambda} (q\sb{\alpha\beta\gamma})$).
The defining equations of them are given in Tables~\ref{table:AQ1} and~\ref{table:AQ2}.
Hence $\gamma\sb \lambda \in \gs\sb A$
by Proposition~\ref{prop:ings}.
\end{proof}
\begin{table}
{\small
\begin{eqnarray*}
Q\sprime\sb{1} &:&\lambda\,{X}^{2}+{Y}^{2}+\left ({\lambda}^{2}+\lambda\right ){Z}^{2
}+YZ+{\lambda}^{2}ZX
=0, \\ 
Q\sprime\sb{2} &:&(\lambda+1) X^2+Y^2+YZ+(\lambda^2+1) ZX
=0, \\ 
Q\sprime\sb{3} &:&(\lambda+1) X^2+ \lambda Y^2+\lambda ^2 YZ+(\lambda^2+1) ZX
=0, \\ 
Q\sprime\sb{4} &:&\lambda\,{X}^{2}+\left (\lambda+1\right ){Y}^{2}+\left ({\lambda}^{
2}+1\right )YZ+{\lambda}^{2}ZX
=0, \\ 
Q\sprime\sb{5} &:&{X}^{2}+\lambda\,{Y}^{2}+\left ({\lambda}^{2}+\lambda\right ){Z}^{2
}+{\lambda}^{2}YZ+ZX
=0, \\ 
Q\sprime\sb{6} &:&{X}^{2}+\left (\lambda+1\right ){Y}^{2}+\left ({\lambda}^{2}+1
\right )YZ+ZX
=0, \\ 
Q\sprime\sb{7} &:&{X}^{2}+\left (\lambda+1\right ){Y}^{2}+\left ({\lambda}^{2}+\lambda\right ){Z}^{2}+\left
({\lambda}^{2}+1\right )YZ+ZX =0, \\ 
Q\sprime\sb{8} &:&{X}^{2}+\lambda\,{Y}^{2}+{\lambda}^{2}YZ+ZX
=0, \\ 
Q\sprime\sb{9} &:&\lambda\,{X}^{2}+\left (\lambda+1\right ){Y}^{2}+\left ({\lambda}^{
2}+\lambda\right ){Z}^{2}+\left ({\lambda}^{2}+1\right )YZ+{\lambda
}^{2}ZX
=0, \\ 
Q\sprime\sb{10} &:&(\lambda+1) X^2 +\lambda Y^2 +(\lambda^2+\lambda) Z^2+\lambda^2 YZ+(\lambda^2+1)ZX
=0, \\ 
Q\sprime\sb{11} &:&(\lambda+1) X^2 +							 Y^2 +(\lambda^2+\lambda) Z^2+									 YZ+(\lambda^2+1)ZX
=0, \\ 
Q\sprime\sb{12} &:&\lambda\,{X}^{2}+{Y}^{2}+YZ+{\lambda}^{2}ZX
=0. 
\end{eqnarray*}
}
\caption{Defining equations of the conic curves $Q\sb{i}\sprime$}\label{table:AQ1}
\end{table}
\begin{table}
{\small
\begin{eqnarray*}
Q\sb{114} &:& {X}^{2}+{Y}^{2}+XY+\lambda\,{Z}^{2} +\left (\lambda+1\right )ZX
=0, \\ 
Q\sb{123} &:& {X}^{2}+{Y}^{2}+XY+\left (\lambda+1\right ){Z}^{2}+\lambda\,ZX
=0, \\ 
Q\sb{132} &:& {X}^{2}+{Y}^{2}+XY+\left ({\lambda}^{2}+\lambda\right ){Z}^{2}+ZX
=0, \\ 
Q\sb{141} &:& {X}^{2}+{Y}^{2}+XY+\left ({\lambda}^{2}+\lambda+1\right ){Z}^{2}
=0, \\ 
Q\sb{213} &:&{X}^{2}+{Y}^{2}+XY+ YZ+\lambda\,ZX
=0, \\ 
Q\sb{224} &:&{X}^{2}+{Y}^{2}+XY+ YZ+\left (\lambda+1\right )ZX
=0, \\ 
Q\sb{231} &:& {X}^{2}+{Y}^{2}+XY+\left ({\lambda}^{2}+\lambda\right ){Z}^{2} +YZ
=0, \\ 
Q\sb{242} &:& {X}^{2}+{Y}^{2}+XY+\left ({\lambda}^{2}+\lambda\right ){Z}^{2} +YZ+ZX
=0, \\ 
Q\sb{312} &:&{X}^{2}+{Y}^{2}+XY+ \lambda\,YZ+ZX
=0, \\ 
Q\sb{321} &:& {X}^{2}+{Y}^{2}+XY+\left (\lambda+1\right ){Z}^{2} +\lambda\,YZ
=0, \\ 
Q\sb{334} &:&{X}^{2}+{Y}^{2}+XY+ \lambda\,YZ+\left (\lambda+1\right )ZX
=0, \\ 
Q\sb{343} &:& {X}^{2}+{Y}^{2}+XY+\left (\lambda+1\right ){Z}^{2} +\lambda\,YZ+\lambda\,ZX
=0, \\ 
Q\sb{411} &:& {X}^{2}+{Y}^{2}+XY+\lambda\,{Z}^{2} +\left (\lambda+1\right )YZ
=0, \\ 
Q\sb{422} &:&{X}^{2}+{Y}^{2}+XY+ \left (\lambda+1\right )YZ+ZX
=0, \\ 
Q\sb{433} &:&{X}^{2}+{Y}^{2}+XY+ \left (\lambda+1\right )YZ+\lambda\,ZX
=0, \\ 
Q\sb{444} &:& {X}^{2}+{Y}^{2}+XY+\lambda\,{Z}^{2} +\left (\lambda+1\right )YZ+
\left (\lambda+1\right )ZX
=0. 
\end{eqnarray*}
}
\caption{Defining equations of the conic curves $Q\sb{\alpha\beta\gamma}$}
\label{table:AQ2}
\end{table}
\begin{remark}
The polynomial  $\GA [\lambda]$   defines 
the nodal splitting curve
$$
M \cup L\sb{12, 1} \cup L\sb{18, 1} \cup L\sb{13, 3} \cup Q\sb{242}.
$$
See Proposition~\ref{prop:nodalsplit}.
\end{remark}
\begin{remark}\label{rem:limitDKA}
When $\lambda\in \{\omega, \bar\omega\}$,
the set $\gamma\sb\lambda (\Ps)$ coincides with $\Pt (\Ff)$,
and the point $[\GA[\lambda]]\in \moduli $ is the Dolgachev-Kondo point.
\end{remark}
Let $k(\lambda)$ be the rational function field with variable $\lambda$.
For each $\sigma\in \Aut (\CA)$,
we calculate the unique triple
$$
(g\sb{\sigma}, t\sb{\sigma}, \lambda\sp{\sigma})\;\in\; \PGL (3, k(\lambda)) \times \widetilde{T}\times k(\lambda)
$$
such that 
$$
g\sb{\sigma} \circ (\gamma\sb{\lambda}\circ \sigma)\circ t\sb{\sigma}=\gamma\sb{\lambda\sp{\sigma}}
$$
holds (see~Claim~\ref{claim:tripleA}.)
The calculation is done as follows: 
$g\sb\sigma$ is the unique linear automorphism of $\Pt$ characterized by
\begin{eqnarray*}
&& g\sb{\sigma} (\gamma\sb{\lambda}(\sigma (P\sb {18})))=[1,0,0], \nonumber\\
&&g\sb{\sigma} (\gamma\sb{\lambda}(\sigma (P\sb {12})))=[0,1,0],\;\;\;
g\sb{\sigma} (\gamma\sb{\lambda}(\sigma (P\sb {13})))=[1,1,0]\;\;\;\quand\\
&&g\sb{\sigma} (\gamma\sb{\lambda}(\sigma (P\sb {20})))=[0,0,1],\;\;\;
g\sb{\sigma} (\gamma\sb{\lambda}(\sigma (P\sb {19})))=[1,0,1];
\nonumber
\end{eqnarray*}
$t\sb \sigma\in  \widetilde T$ is given by
$$
t\sb \sigma=
\begin{cases}
\id & \textrm{if $g\sb{\sigma} (\gamma\sb{\lambda}(\sigma (P\sb 1)))=[1, \omega, 0]$},\\
T & \textrm{if $g\sb{\sigma} (\gamma\sb{\lambda}(\sigma (P\sb 1)))=[1, \bar\omega, 0]$};
\end{cases}
$$
and $\lambda\sp{\sigma}$ is the rational
function of the parameter $\lambda$ satisfying
$$
g\sb{\sigma} (\gamma\sb{\lambda}(\sigma (P\sb {10})))=[0, \lambda\sp{\sigma}, 1].
$$
The map $\sigma\mapsto t\sb{\sigma}$ is a homomorphism from $\Aut (\CA)$ to $\widetilde T$.
We  put
$$
N\sb A:=\Ker (\Aut (\CA)\to \widetilde{T} ).
$$
From the proof of Proposition~\ref{prop:Agammalambda},
we obtain the following:
\begin{corollary}\label{cor:connectedcompsA}
The space $\PGL(3, k)\backslash \gs\sb A$ has exactly two connected components,
each of which is isomorphic to $\A\sp 1\sm \{0,1,\omega,\bar\omega\}$.
Set-theoretically, they are given by
\begin{eqnarray*}
(\PGL(3, k)\backslash \gs\sb A)\sp + &:=&\set{[\gamma\sb\alpha]}{\alpha\in k\sm \{0,1,\omega,\bar\omega\}},
\quand\\
(\PGL(3, k)\backslash \gs\sb A)\sp-&:=&\set{[\gamma\sb\alpha\circ T]}{\alpha\in k\sm \{0,1,\omega,\bar\omega\}}.
\end{eqnarray*}
The group $N\sb A$ acts on $(\PGL(3, k)\backslash \gs\sb A)\sp +$, and
the moduli curve $\moduli\sb A$ is equal to the quotient space $(\PGL(3, k)\backslash \gs\sb A)\sp +/N\sb A$.
\end{corollary}
Let
$$
\map{p\sb A}{ \A\sp 1\sm \{0,1,\omega, \bar\omega\} \cong (\PGL (3, k) \backslash \GGG\sb{A})\sp + }{%
\moduli\sb{A} = (\PGL (3, k) \backslash \GGG\sb{A})\sp + / N\sb A}
$$
denote  the natural projection.
For $\alpha\in k\sm\{0, 1, \omega, \bar\omega\}$,
let $P [\alpha]$ be the  point of $\A\sp 1\sm \{0,1,\omega, \bar\omega\}$
given by $\lambda=\alpha$.
Then $p\sb A (P[\alpha])\in\moduli\sb A$ corresponds to the isomorphism class 
of the polarized supersingular $K3$ surface $( X\sb{\GA [\alpha]}, \pol\sb{\GA [\alpha]})$.
\begin{proposition}\label{prop:groupA}
The set 
$p\sb A\inv (p\sb A (P[\alpha]))$ is equal to 
\begin{equation}\label{eq:fiberA}
\{ \;
P[\alpha], \;
P[1/\alpha],  \;
P[\alpha+1],  \;
P[1/(\alpha+1)],  \;
P[\alpha/(\alpha+1)],  \;
P[(\alpha+1)/\alpha]  \;
\},
\end{equation}
and  $\Aut( X\sb{\GA [\alpha]}, \pol\sb{\GA [\alpha]})$
is equal to the group~\eqref{eq:groupA}.
\end{proposition}
\begin{proof}
The set 
$\shortset{\lambda\sp{\sigma}}{\sigma\in N\sb A}\st k(\lambda)$
coincides with the group $\Gamma\sb A$ given in Theorem~\ref{thm:A}.
The fiber $p\sb A\inv (p\sb A (P[\alpha]))$ is therefore 
equal to~\eqref{eq:fiberA}.
Note that the fiber $p\sb A\inv (p\sb A (P[\alpha]))$ 
 consists of six distinct points for any $\alpha\in k\sm\{0,1,\omega, \bar\omega\}$;
 that is, the action of $\Gamma\sb A$ on $(\PGL (3, k) \backslash \GGG\sb{A})\sp +$
 is free.
Hence,
for any $\alpha\in k\sm\{0,1,\omega,\bar\omega\}$ and any $\sigma\in \Aut (\CA)$, 
the projective equivalence classes 
$[\gamma\sb{\alpha}]$ and 
$$
[\gamma\sb{\alpha}\circ \sigma]=[\gamma\sb{\alpha\sp\sigma}\circ t\sb{\sigma}]
\;\in\; \PGL(3, k)\backslash\gs\sb C
$$
coincide if and only if $t\sb{\sigma}=\id$ and $\lambda\sp{\sigma}=\lambda$ hold.
Therefore, using Corollary~\ref{cor:projaut2},  we can obtain  
$\Aut(X\sb{\GA [\alpha]}, \pol\sb{\GA [\alpha]})$ from the subgroup 
$$
\set{g\sb\sigma}{t\sb{\sigma}=\id\quand \lambda\sp{\sigma}=\lambda}\;\;\st\;\; \PGL(3, k(\lambda))
$$
by substituting   $\alpha$ for  $\lambda$.
\end{proof}
\begin{corollary}
We have
$\moduli\sb{A}=\Spec k[J\sb A, 1/J\sb A]$,
where  $J\sb A=(\lambda^2+\lambda+1)/\lambda^2(\lambda+1)^2$.
The morphism $p\sb A$ is an \'etale  Galois covering
with Galois group $\Gamma\sb A$,
which is isomorphic to $\SSSS\sb 3$.
\end{corollary}
\section{The moduli curve corresponding to the code $\CB$}\label{sec:B}
In this section,
we prove Theorem~\ref{thm:B}.
\par
\medskip
Let $\AF$ be  the affine plane over $\F\sb 3$,
$P(\AF)$ the set of rational points of $\AF$, 
and   $L(\AF)$  the set of rational  affine lines of $\AF$.
Each element of $P(AF)$ is expressed  by a pair $aa\sprime$ of elements of $\F\sb 3$, and 
each element of $L(\AF)$ is expressed as a subset $\{aa\sprime, bb\sprime, cc\sprime\}$
of $P(\AF)$ with cardinality $3$.
We have
$$
| P(\AF) |=9 \quand | L(\AF)| =12.
$$
The incidence relation
$$
\set{(p, \ell)\in P(\AF) \times L(\AF)}{p\in \ell}
$$
is called the \emph{Hesse configuration}
(\cite{DolCon}).
The automorphism group
$$
\GH:=\set{\sigma\in \SSSS (P(\AF))}{\sigma (\ell)\in L(\AF)\;\;\textrm{for all $\ell \in L(\AF)$}\;\;}
$$
of this configuration is isomorphic to the group of affine transformations of $\AF$ 
defined over $\F\sb 3$.
In particular, the order of $\GH$ is $432$.
\par
\medskip
We define injective maps
$$
C: P(\AF) \to \Ps  \quand T: L(\AF) \to \Ps
$$
by Table~\ref{table:CpointsBTpointsB}.
Then $\Ps$ is a disjoint union of $C(P(\AF))$ and $T(L(\AF))$.
A point $P\in \Ps$ is called a \emph{$C$-point}  or a \emph{$T$-point}
according to whether 
$P\in C(P(\AF))$  or $P\in T(L(\AF))$.
\begin{table}
$$
\renewcommand{\arraystretch}{1.2}
\begin{array}{|c|ccccccccc|}
\hline
 aa\sprime & 00 & 01 & 02 & 10 & 11 & 12 & 20 & 21 & 22 \\
\hline 
C(aa\sprime) & P\sb{17} & P\sb{13} & P\sb {5}  & P\sb{10} & P\sb {8} & P\sb {6} & P\sb {2} & P\sb{3} &P\sb 1 \\
\hline 
\end{array}
$$
\vskip 5pt
$$
\renewcommand{\arraystretch}{1.2}
\begin{array}{|ccc|c |}
\hline
aa\sprime & bb\sprime & cc\sprime  & T(\ell)  \\
\hline
00 & 01 & 02 & P\sb{21}  \\
00 & 10 & 20 & P\sb{20}  \\
00 & 11 & 22 & P\sb{19} \\
00 & 12 & 21 & P\sb{18} \\
01 & 10 & 22 & P\sb{16} \\
01 & 11 & 21 & P\sb{15}  \\
\hline
\end{array}
\qquad
\begin{array}{|ccc|c |}
\hline
aa\sprime & bb\sprime & cc\sprime  & T(\ell)  \\
\hline
01 & 12 & 20 & P\sb{14} \\
02 & 10 & 21 & P\sb{11} \\
02 & 11 & 20 & P\sb{9}  \\
02 & 12 & 22 & P\sb{7}  \\
10 & 11 & 12 & P\sb{12}  \\
20 & 21 & 22 & P\sb{4} \\
\hline
\end{array}
$$
\vskip 5pt
\caption{$C$-points $C(aa\sprime)$ and $T$-points $T(\ell)$ for $\ell=\{aa\sprime, bb\sprime,
cc\sprime\}$}\label{table:CpointsBTpointsB}
\end{table}
The code $\CB$ is described as follows.
\par\smallskip
The linear words of $\CB$ are precisely the words
$$
l\sb{aa\sprime}:=\{ C(aa\sprime), T(\ell\sb 1), T(\ell\sb 2), T(\ell\sb 3), T(\ell\sb 4) \}\qquad (aa\sprime \in P(\AF)),
$$
where $\ell\sb 1, \dots, \ell\sb 4 \in L (\AF)$ are the four affine lines passing through
the point $aa\sprime \in P(\AF)$.
\par\smallskip
There are two types of quadratic words.
\par\smallskip
(I) 
Let $\ell=\{aa\sprime, bb\sprime, cc\sprime\}$ be an element of $L(\AF)$.
There exists a unique pair of  distinct  affine lines
$$
\ell\sb 1=\{a\sb 1 a\sprime\sb 1, b\sb 1 b\sprime\sb 1, c\sb 1 c\sprime\sb 1\}\ne\ell,
\qquad
\ell\sb 2=\{a\sb 2 a\sprime\sb 2, b\sb 2 b\sprime\sb 2, c\sb 2 c\sprime\sb 2\}\ne\ell
$$
that are parallel to $\ell$.
Then the word
$$
q\sb{\ell}:=\{
C(a\sb 1 a\sb 1\sprime),
C(b\sb 1 b\sb 1\sprime),
C(c\sb 1 c\sb 1\sprime),
C(a\sb 2 a\sb 2\sprime),
C(b\sb 2 b\sb 2\sprime),
C(c\sb 2 c\sb 2\sprime),
T(\ell\sb 1),
T(\ell\sb 2)
\}
$$
is a quadratic word of $\CB$.
\par
\begin{figure}
\includegraphics[height=5.5cm]{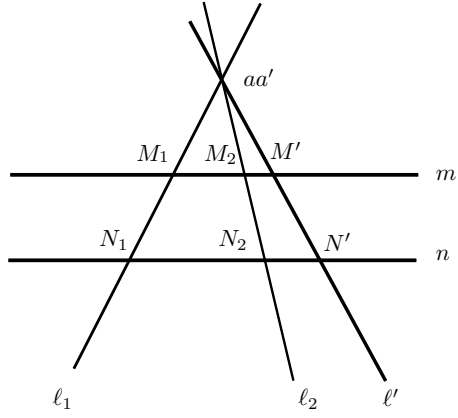}
\caption{Intersection points}\label{fig:MNs}
\end{figure}

(II) 
Let $\ell\sb 1$ and $\ell\sb 2$ be two distinct elements of $L(\AF)$
that are not parallel,
and let $aa\sprime\in P(\AF)$ be the intersection point of $\ell\sb 1$ and $\ell\sb 2$.
Then there exists a  pair $\{ m, n \}$ of elements of $L(\AF)$ 
with the following properties:
\begin{itemize}
\item[(i)] $m$ and $n$ are parallel, 
\item[(ii)] $aa\sprime \notin m$, $aa\sprime \notin n$, and
\item[(iii)] none of the pairs $(\ell\sb 1, m)$, $(\ell\sb 2, m)$, $(\ell\sb 1, n)$, $(\ell\sb 2, n)$  are parallel.
\end{itemize}
For such a pair $\{ m, n \}$,  there exists a unique line
$\ell\sprime\in L(\AF)$ such that 
\begin{itemize}
\item[(a)] $aa\sprime\in \ell\sprime$,
\item[(b)] is distinct from $\ell\sb 1$ and $\ell\sb 2$, and 
\item[(c)] intersects both of  $m$ and $n$.
\end{itemize}
We denote the intersection points of these affine lines as in Figure~\ref{fig:MNs}.
Then the word
$$
q\sprime \sb{\ell\sb 1, \ell\sb 2}:=\{
C(M\sb 1), C(M\sb 2), C(N\sb 1), C(N\sb 2), T (M\sb 1 N\sprime ), T(M\sb 2 N\sprime ), T(N\sb 1 M\sprime ), T(N\sb 2 M\sprime )
\}
$$
is a quadratic word of $\CB$,
where $MN \in L(\AF)$ denotes the affine line containing the points $M$ and $N$.
For each $(\ell\sb 1, \ell\sb 2)$, there exist exactly two pairs  satisfying (i), (ii) and (iii).
However, 
the word $q\sprime \sb{\ell\sb 1, \ell\sb 2}$ is independent of the choice of the pair.
\par\smallskip
There exist $12$ quadratic words of type I,
and $54$ quadratic words of type II.
The quadratic words of $\CB$ are precisely these $66$ words.
The linear and quadratic words of $\CB$ are explicitly presented in
Tables~\ref{table:linearB},~\ref{table:QSB} and~\ref{table:QFB}.
\begin{table}[tb]
$$
\begin{array}{ccccccccc}
l\sb{00} & : & \{ &17, &18, &19, &20, &21&\},  \\ 
l\sb{01} & : & \{ &13, &14, &15, &16, &21&\},  \\ 
l\sb{02} & : & \{ &5, &7, &9, &11, &21&\},  \\ 
l\sb{10} & : & \{ &10, &11, &12, &16, &20&\},  \\ 
l\sb{11} & : & \{ &8, &9, &12, &15, &19&\},  \\ 
l\sb{12} & : & \{ &6, &7, &12, &14, &18&\},  \\ 
l\sb{20} & : & \{ &2, &4, &9, &14, &20&\},  \\ 
l\sb{21} & : & \{ &3, &4, &11, &15, &18&\},  \\ 
l\sb{22} & : & \{ &1, &4, &7, &16, &19&\}.\\
\phantom{l\sb{}} &&&&&&&&
\end{array}
$$
\caption{Linear words of $\CB$}\label{table:linearB}
\end{table}
\begin{table}[tb]
$$
\begin{array}{cccccccccccc}
q \sb{00,01,02}& : & \{ &1, &2, &3, &4, &6, &8, &10, &12&\},  \\ 
q \sb{00,10,20}& : & \{ &1, &3, &5, &6, &7, &8, &13, &15&\},  \\ 
q \sb{00,11,22}& : & \{ &2, &3, &5, &6, &10, &11, &13, &14&\},  \\ 
q \sb{00,12,21}& : & \{ &1, &2, &5, &8, &9, &10, &13, &16&\},  \\ 
q \sb{01,10,22}& : & \{ &2, &3, &5, &6, &8, &9, &17, &18&\},  \\ 
q \sb{01,11,21}& : & \{ &1, &2, &5, &6, &7, &10, &17, &20&\},  \\
q \sb{01,12,20}& : & \{ &1, &3, &5, &8, &10, &11, &17, &19&\},  \\  
q \sb{02,10,21}& : & \{ &1, &2, &6, &8, &13, &14, &17, &19&\},  \\  
q \sb{02,11,20}& : & \{ &1, &3, &6, &10, &13, &16, &17, &18&\},  \\  
q \sb{02,12,22}& : & \{ &2, &3, &8, &10, &13, &15, &17, &20&\},  \\  
q \sb{10,11,12}& : & \{ &1, &2, &3, &4, &5, &13, &17, &21&\},  \\  
q \sb{20,21,22}& : & \{ &5, &6, &8, &10, &12, &13, &17, &21&\}.  \\  
\phantom{Q\sb{00,10,20}} &  & &\phantom{20} &\phantom{20} &\phantom{20} &\phantom{20} &\phantom{20} &\phantom{20} &\phantom{20}
&\phantom{20}& 
\end{array}
$$
\caption{Quadratic words of type I in $\CB$}\label{table:QSB}
\end{table}
{\tiny
\begin{table}[tb]
$$
\renewcommand{\arraystretch}{1.0}
\begin{array}{cc|ccccccccccc}
 T(\ell) &  T(\ell\sprime) &  &&&&& \hskip - 5mm\hbox to 1mm{$q\sprime \sb{\ell, \ell\sprime}$} &&\\
\hline
 18 & 19 & \{ & 1, &3, &6, &8, &9, &11, &14, &16  &\} \\ 
 18 & 20 & \{ & 2, &3, &6, &7, &9, &10, &15, &16  &\} \\ 
 18 & 21 & \{ & 3, &4, &5, &6, &9, &12, &13, &16  &\} \\ 
 19 & 20 & \{ & 1, &2, &7, &8, &10, &11, &14, &15  &\} \\ 
 19 & 21 & \{ & 1, &4, &5, &8, &11, &12, &13, &14  &\} \\ 
 20 & 21 & \{ & 2, &4, &5, &7, &10, &12, &13, &15  &\} \\ 
 11 & 12 & \{ & 3, &4, &5, &6, &8, &14, &19, &21  &\} \\ 
 11 & 16 & \{ & 1, &3, &5, &9, &13, &14, &18, &19  &\} \\ 
 11 & 20 & \{ & 2, &3, &5, &7, &14, &15, &17, &19  &\} \\ 
 12 & 16 & \{ & 1, &4, &6, &8, &9, &13, &18, &21  &\} \\ 
 12 & 20 & \{ & 2, &4, &6, &7, &8, &15, &17, &21  &\} \\ 
 16 & 20 & \{ & 1, &2, &7, &9, &13, &15, &17, &18  &\} \\ 
 4 & 9 & \{ & 1, &3, &5, &8, &12, &16, &18, &21  &\} \\ 
 4 & 14 & \{ & 1, &3, &6, &11, &12, &13, &19, &21  &\} \\ 
 4 & 20 & \{ & 1, &3, &7, &10, &12, &15, &17, &21  &\} \\ 
 9 & 14 & \{ & 5, &6, &8, &11, &13, &16, &18, &19  &\} \\ 
 9 & 20 & \{ & 5, &7, &8, &10, &15, &16, &17, &18  &\} \\ 
 14 & 20 & \{ & 6, &7, &10, &11, &13, &15, &17, &19  &\} \\ 
 14 & 15 & \{ & 2, &3, &6, &7, &8, &11, &19, &20  &\} \\ 
 14 & 16 & \{ & 1, &2, &6, &9, &10, &11, &18, &19  &\} \\ 
 14 & 21 & \{ & 2, &4, &5, &6, &11, &12, &17, &19  &\} \\ 
 15 & 16 & \{ & 1, &3, &7, &8, &9, &10, &18, &20  &\} \\ 
 15 & 21 & \{ & 3, &4, &5, &7, &8, &12, &17, &20  &\} \\ 
 16 & 21 & \{ & 1, &4, &5, &9, &10, &12, &17, &18  &\} \\ 
 9 & 12 & \{ & 2, &4, &5, &6, &10, &16, &18, &21  &\} \\ 
 9 & 15 & \{ & 2, &3, &5, &7, &13, &16, &18, &20  &\} \\ 
 9 & 19 & \{ & 1, &2, &5, &11, &14, &16, &17, &18  &\} \\ 
 12 & 15 & \{ & 3, &4, &6, &7, &10, &13, &20, &21  &\} \\ 
 12 & 19 & \{ & 1, &4, &6, &10, &11, &14, &17, &21  &\} \\ 
 15 & 19 & \{ & 1, &3, &7, &11, &13, &14, &17, &20  &\} \\ 
 4 & 11 & \{ & 1, &2, &5, &10, &12, &14, &19, &21  &\} \\
 4 & 15 & \{ & 1, &2, &7, &8, &12, &13, &20, &21  &\} \\ 
 4 & 18 & \{ & 1, &2, &6, &9, &12, &16, &17, &21  &\} \\ 
 11 & 15 & \{ & 5, &7, &8, &10, &13, &14, &19, &20  &\} \\ 
 11 & 18 & \{ & 5, &6, &9, &10, &14, &16, &17, &19  &\} \\ 
 15 & 18 & \{ & 6, &7, &8, &9, &13, &16, &17, &20  &\} \\ 
 7 & 9 & \{ & 1, &2, &6, &8, &15, &16, &18, &20  &\} \\ 
 7 & 11 & \{ & 1, &3, &6, &10, &14, &15, &19, &20  &\} \\ 
 7 & 21 & \{ & 1, &4, &6, &12, &13, &15, &17, &20  &\} \\ 
 9 & 11 & \{ & 2, &3, &8, &10, &14, &16, &18, &19  &\} \\ 
 9 & 21 & \{ & 2, &4, &8, &12, &13, &16, &17, &18  &\} \\ 
 11 & 21 & \{ & 3, &4, &10, &12, &13, &14, &17, &19  &\} \\ 
 7 & 12 & \{ & 1, &4, &5, &8, &10, &15, &20, &21  &\} \\ 
 7 & 14 & \{ & 1, &2, &5, &11, &13, &15, &19, &20  &\} \\ 
 7 & 18 & \{ & 1, &3, &5, &9, &15, &16, &17, &20  &\} \\ 
 12 & 14 & \{ & 2, &4, &8, &10, &11, &13, &19, &21  &\} \\ 
 12 & 18 & \{ & 3, &4, &8, &9, &10, &16, &17, &21  &\} \\ 
 14 & 18 & \{ & 2, &3, &9, &11, &13, &16, &17, &19  &\} \\ 
 4 & 7 & \{ & 2, &3, &5, &6, &12, &15, &20, &21  &\} \\ 
 4 & 16 & \{ & 2, &3, &9, &10, &12, &13, &18, &21  &\} \\ 
 4 & 19 & \{ & 2, &3, &8, &11, &12, &14, &17, &21  &\} \\ 
 7 & 16 & \{ & 5, &6, &9, &10, &13, &15, &18, &20  &\} \\ 
 7 & 19 & \{ & 5, &6, &8, &11, &14, &15, &17, &20  &\} \\ 
 16 & 19 & \{ & 8, &9, &10, &11, &13, &14, &17, &18  &\} \\ 
 &&& &\phantom{20} &\phantom{20} &\phantom{20} &\phantom{20} &\phantom{20} &\phantom{20} &\phantom{20}
&\phantom{20}& 
\end{array}
$$
\caption{Quadratic words of type II  in $\CB$}\label{table:QFB}
\end{table}
}
\par
\medskip
The following proposition can be checked easily:
\begin{proposition}\label{prop:B1}
Let $\ell=\{aa\sprime, bb\sprime, cc\sprime\}$ be an element of $ L(\AF)$.
Then the  quadratic word $q\sb{\ell}$
of type {\rm I} is disjoint from the three linear words
$l\sb{aa\sprime}, l\sb{bb\sprime}, l\sb{cc\sprime}$ 
containing $T(\ell)\in \Ps$.
\end{proposition}
We define a homomorphism
$$
\map{\Psi}{\GH}{\SSSS (\Ps)}
$$
by
$$
\Psi(g) (C(aa\sprime)):= C(g(aa\sprime))
\qquad
\Psi(g) (T(\ell)):= T(g(\ell)).
$$
It is obvious that $\Psi$ is injective.
\begin{proposition}\label{prop:autoCB}
The automorphism group $\Aut (\CB)$ of the code $\CB$ coincides with the image of $\Psi$.
\end{proposition}
\begin{proof}
The above description of the linear and quadratic words in $\CB$
shows that every element in the image of $\Psi$ preserves
the sets of these words.
Since $\CB$ is generated by the word $\Ps\in \Pow(\Ps)$ and these words,
the image of $\Psi$ is contained in $\Aut (\CB)$.
\par
Suppose that $\sigma\in \Aut (\CB)$ is given.
A point $P\in \Ps$ is a $C$-point  if and only if there exists
exactly one linear word in $\CB$ that contains $P$.
Hence $\sigma$ preserves the set of $C$-points.
Via the bijection $C: P(\AF)\cong \Im C$,
we obtain a unique element $\tilde \sigma \in \SSSS (P(\AF))$ such that $\sigma \circ C =C\circ \tilde \sigma$
holds.
When $P=C(aa\sprime)$,  the unique linear word in $\CB$ containing $P$ is just $l\sb{aa\sprime}$.
The Hesse configuration on $\AF$ is recovered from $\CB$ as follows;
a set  $\{aa\sprime, bb\sprime, cc\sprime\}$ of cardinality $3$ 
is  an element  of $L(\AF)$
if and only if the three linear words
$l\sb{aa\sprime}$, $l\sb{bb\sprime}$, $l\sb{cc\sprime}$
have a point in common.
In this case, 
the common point of $l\sb{aa\sprime}$, $l\sb{bb\sprime}$, $l\sb{cc\sprime}$
is just $T(\{aa\sprime, bb\sprime, cc\sprime\})$.
Therefore
we see that  $\tilde \sigma \in \GH$,
and  that $\sigma \circ T =T\circ \tilde \sigma$ holds.
 Thus $\sigma  = \Psi (\tilde \sigma)$.
\end{proof}
Let $\lambda$ be a parameter of the affine line $\A\sp 1$.
We define
$\gamma\sb\lambda : \Ps \to \Pt$ by Table~\ref{table:gammalambdaB}.
We also denote by $\widetilde T=\langle T\rangle$
the subgroup of $\Aut (\CB)$
of order $2$ generated by 
$$
T:=(P\sb 2P\sb 5)(P\sb 3 P\sb 6) (P\sb 4 P\sb 7) (P\sb {10} P\sb{13}) (P\sb{11} P\sb{14}) (P\sb{12} P\sb{15}) (P\sb{20}
P\sb{21}), 
$$
which corresponds to 
the automorphism of the Hesse configuration given by
$aa\sprime\mapsto a\sprime a$.
\begin{table}
\begin{center}
\begin{tabular}{l c|c c  }
$P\sb i$ &\phantom{aa}&\phantom{aa}&  $\gamma\sb{\lambda} (P\sb i)$ \\
\hline
$P\sb 1=C(22)$ &&&  $[\lambda+1, \bar\omega\lambda+\omega, \bar\omega\lambda+\omega]$ \\
$P\sb 2=C(20)$ &&&  $[1, \omega \lambda+\omega, \omega]$ \\
$P\sb 3=C(21)$ &&&  $[\lambda+\bar\omega, 1, \lambda+1]$ \\
$P\sb 4=T(20,21,22)$ &&&  $[1, \omega,\omega]$ \\
$P\sb 5=C(02)$ &&&  $[\lambda,\bar\omega\lambda,\bar\omega\lambda+\bar\omega]$ \\
$P\sb 6=C(12)$ &&&  $[\lambda+\bar\omega, \bar\omega\lambda+\bar\omega, \bar\omega\lambda]$ \\
$P\sb 7=T(02,12,22)$ &&&  $[1, \bar\omega,\bar\omega]$ \\
$P\sb 8=C(11)$ &&&  $[\lambda+1, 1, \lambda]$ \\
$P\sb 9=T(02,11,20)$ &&&  $[1, \bar\omega, \omega]$ \\
$P\sb {10}=C(10)$ &&&  $[1,\omega\lambda+1, 0 ]$ \\
$P\sb {11}=T(02,10,21)$ &&&  $[1, \bar\omega, 0]$ \\
$P\sb {12}=T(10,11,12)$ &&&  $[1,1,0]$ \\
$P\sb {13}=C(01)$ &&&  $[\lambda, 0, \lambda+\bar\omega]$ \\
$P\sb {14}=T(01,12,20)$ &&&  $[1,0,\omega]$ \\
$P\sb {15}=T(01,11,21)$ &&&  $[1,0,1]$ \\
$P\sb {16}=T(01,10,22)$ &&&  $[1,0,0]$ \\
$P\sb {17}=C(00)$ &&&  $[0, \lambda, 1]$ \\
$P\sb {18}=T(00,12,21)$ &&&  $[0, 1, \omega]$ \\
$P\sb {19}=T(00,11,22)$ &&&  $[0,1,1]$ \\
$P\sb {20}=T(00,10,20)$ &&&  $[0,1,0]$ \\
$P\sb {21}=T(00,01,02)$ &&&  $[0,0,1]$ \\
\end{tabular}
\end{center}
\vskip 5pt
\caption{Definition of $\gamma\sb{\lambda}$ for $\CB$}
\label{table:gammalambdaB}
\end{table}
\begin{proposition}
The map $\lambda \mapsto \gamma\sb{\lambda}$ induces an isomorphism   
from $\A\sp 1\sm \{ 0,1, \omega, \bar\omega\}$ to 
$\PGL (3, k) \backslash \GGG\sb{B} / \widetilde T$.
\end{proposition}
\begin{proof}
First note that
$\gamma\sb{\lambda}$ is injective
if and only if 
$$
\lambda\ne 0, \quad
\lambda\ne 1
\quand
\lambda\ne\bar\omega
$$
hold.
\par
\medskip
Suppose that $\lambda\ne 0, 1, \omega$ and $\bar\omega$.
Then $\gamma\sb \lambda$ is injective,
and the image $\gamma\sb{\lambda} (\Ps)$ coincides with $ \ZZ {d\GB[\lambda]}$,
where $\GB[\lambda]$ is given in Theorem~\ref{thm:B}.
Moreover,
for each linear word $l\sb{aa\sprime}$
(resp. each quadratic word  $q\sb{\ell}$ of type I)
(resp. each quadratic word  $q\sprime\sb{\ell, \ell\sprime}$ of type II)
of the code $\CB$,
there exists a line $L\sb{aa\sprime}$ containing $\gamma\sb{\lambda} (\ell\sb{aa\sprime})$
(resp. a  conic curve $Q\sb{\ell}$ containing $\gamma\sb{\lambda} (q\sb{\ell})$)
(resp. a  conic curve $Q\sprime\sb{\ell, \ell\sprime}$ containing $\gamma\sb{\lambda} (q\sprime\sb{\ell, \ell\sprime})$)
given in Tables~\ref{table:linesB},~\ref{table:coniccurvesIB},~\ref{table:coniccurvesIIB}.
The conic curves
in Tables~~\ref{table:coniccurvesIB} and~\ref{table:coniccurvesIIB}
are nonsingular because  $\lambda\notin\{0,1,\omega,\bar\omega\}$.
Hence $\gamma\sb{\lambda}$ is in $\GGG\sb{B}$ by Proposition~\ref{prop:ings}.
\begin{table}
$$
\begin{array}{c|l}
aa\sprime &  \hbox{The defining equation of $L\sb{aa\sprime}$}\\
\hline
00 & X=0 \\
01 & Y=0 \\
02 & X+\omega Y=0 \\
10 & Z=0 \\
11 & X+Y+Z=0 \\
12 & \omega X+\omega Y+Z =0\\
20 & \omega X+Z =0\\
21 & X+\omega Y +Z =0\\
22 & Y+Z =0\\
\end{array}
$$
\caption{Defining equations of the  lines $L\sb{aa\sprime}$}\label{table:linesB}
\end{table}
\begin{table}
$$
\renewcommand{\arraystretch}{1.3}
\begin{array}{ccc|c}
aa\sprime & bb\sprime & cc\sprime & \hbox{The defining equation of $Q\sb{aa\sprime, bb\sprime, cc\sprime}$} \\
\hline 
00 & 01& 02&\left (\lambda+\bar\omega\right ){X}^{2}+\bar\omega {Y}^
{2}+\left (\lambda+\omega\right ){Z}^{2}+\lambda XY
\\ 
00 & 10& 20&\left (\omega \lambda+1\right ){X}^{2}+\left (\bar\omega \lambda+1\right ){Y}^{2}+\omega \lambda {Z}^{2}+ZX
\\ 
00 & 11& 22&\left (\bar\omega \lambda+\omega\right ){X}^{2}+\bar\omega {Y}^{2}+\omega \lambda {Z}^{2}+\left
(\lambda+1\right)XY+\left (\lambda+1\right )ZX
\\ 
00 & 12& 21&{Y}^{2}+\lambda {Z}^{2}+\left (\omega \lambda+1\right )XY+ \left (\lambda+\bar\omega\right )ZX
\\ 
01 & 10& 22&\left (\bar\omega \lambda+1\right ){X}^{2}+{Y}^{2}+ \bar\omega \lambda {Z}^{2}+\left (\lambda+\bar\omega\right )Y Z
\\ 
01 & 11& 21&\left (\omega \lambda+1\right ){X}^{2}+\lambda {Z}^{2}+XY+YZ
\\ 
01 & 12& 20&\left (\lambda+\bar\omega\right ){X}^{2}+{Y}^{2}+\lambda {Z}^{2}+
\left (\omega \lambda+\omega\right )XY+\left (\lambda+1\right )YZ
\\ 
02 & 11& 20&\omega {Y}^{2}+\lambda {Z}^{2}+\left (\bar\omega \lambda+\omega\right )XY+\left (\omega \lambda+1\right )YZ+\left (\lambda+\bar\omega\right )ZX
\\ 
02 & 12& 22&\left (\omega \lambda+1\right ){X}^{2}+\omega \lambda {Z}^{2}+X
Y+\omega YZ+ZX
\\ 
10 & 11& 12&\left (\lambda+\bar\omega\right ){X}^{2}+{Y}^{2}+\lambda YZ+\lambda ZX
\\ 
20 & 21& 22&\left (\lambda+\bar\omega\right ){X}^{2}+\bar\omega {Y}^
{2}+\lambda XY+ \bar\omega \lambda YZ+\lambda ZX
\\ 
\end{array}
$$
\caption{Defining equations of the conic curves of type I}\label{table:coniccurvesIB}
\end{table}
{\tiny
\begin{table}
$$
\renewcommand{\arraystretch}{1.3}
\begin{array}{cc|c}
 T(\ell) & T(\ell\sprime)   & \hbox{The defining equation of $Q\sprime\sb{\ell, \ell\sprime}$} \\
\hline 
  18& 19&  \bar\omega\lambda{Y}^{2}+ \bar\omega{Z}^{2}+\omega \lambda XY+ZX
\\ 
  18& 20&{Y}^{2}+ (\lambda+1 ){Z}^{2}+ (\omega \lambda+1 )XY+ (\lambda+1 )ZX
\\ 
  18& 21& (\lambda+1 ){Y}^{2}+\lambda {Z}^{2}+ (\lambda+1 )XY+ ( \lambda+\bar\omega)ZX
\\ 
  19& 20& ( \bar\omega\lambda+\omega ){X}^{2}+ \bar\omega{Y}^{2}+ (\omega \lambda+1 ){Z}^{2}+ (\lambda+1 )XY+ ( \lambda+\bar\omega)ZX
\\ 
  19& 21& ( \bar\omega\lambda+\omega ){X}^{2}+ (\lambda+\bar\omega ){Y}^{2}+\omega \lambda {Z}^{2}+ (\omega \lambda+1 )XY+ (\lambda+1 )ZX
\\ 
  20& 21& (\omega \lambda+1 ){X}^{2}+{Y}^{2}+\omega \lambda {Z}^{2}+\omega \lambda XY+ZX
\\ 
 11& 12&\omega {X}^{2}+ ( \bar\omega\lambda+\omega ){Y}^{2}+ ( \bar\omega\lambda+\omega )YZ+ZX
\\ 
 11& 16& ( \bar\omega\lambda+\omega ){X}^{2}+  \bar\omega\lambda{Y}^{2}+\omega \lambda {Z}^{2}+\lambda YZ+ (\lambda+1 )ZX
\\ 
 11& 20&  \bar\omega\lambda{X}^{2}+\omega {Y}^{2}+\omega \lambda {Z}^{2}+ \omega (\lambda+1) YZ+\lambda ZX
\\ 
 12& 16& ( \lambda+\bar\omega){X}^{2}+ \bar\omega (\lambda+1 ){Y}^{2}+ \omega (\lambda+1) YZ+\lambda ZX
\\ 
 12& 20& (\lambda+1 ){X}^{2}+{Y}^{2}+\lambda YZ+ (\lambda+1)ZX
\\ 
 16& 20& (\omega \lambda+1 ){X}^{2}+ \bar\omega{Y}^{2}+\omega \lambda {Z}^{2}+ ( \bar\omega\lambda+\omega )YZ+ZX
\\ 
  4& 9& \bar\omega (\lambda+1 ){Y}^{2}+ \bar\omega (\lambda+1 )XY+ \omega (\lambda+1) YZ+ZX
\\ 
  4& 14& ( \lambda+\bar\omega){X}^{2}+ (\omega \lambda+1 ){Y}^{2}+ ( \bar\omega\lambda+\omega )XY+ (\omega \lambda+1 )YZ+\lambda ZX
\\ 
  4& 20& ( \lambda+\bar\omega){X}^{2}+ \bar\omega{Y}^{2}+\lambda XY+  \bar\omega\lambda YZ+ (\lambda+\bar\omega )ZX
\\ 
  9& 14&\omega \lambda {Y}^{2}+\lambda {Z}^{2}+\lambda XY+  \bar\omega\lambda YZ+ ( \lambda+\bar\omega)ZX
\\ 
  9& 20&\omega {Y}^{2}+\lambda {Z}^{2}+ ( \bar\omega\lambda+\omega )XY+ (\omega \lambda+1 )YZ+\lambda ZX
\\ 
  14& 20& (\omega \lambda+1 ){X}^{2}+\omega {Y}^{2}+\omega \lambda {Z}^{2}+ \bar\omega (\lambda+1 )XY+ \omega (\lambda+1) YZ+ZX
\\ 
  14& 15& \lambda {X}^{2}+ ( \lambda+\bar\omega){Z}^{2}+\omega 
\lambda XY+ ( \lambda+\bar\omega)YZ
\\ 
  14& 16& ( \bar\omega\lambda+\omega ){X}^{2}+ \bar\omega{Y}^{2}+\omega {Z}^{2}+ (\lambda+1 )XY+YZ
\\ 
  14& 21& \bar\omega{X}^{2}+{Y}^{2}+\lambda {Z}^{2}+\omega XY+ (\lambda+1 )YZ
\\ 
  15& 16& ( \bar\omega\lambda+\omega ){X}^{2}+ \bar\omega (\lambda+1 ){Z}^{2}+\omega XY+ (\lambda+1 )YZ
\\ 
  15& 21& (\lambda+1 ){X}^{2}+\lambda {Z}^{2}+ (\lambda+1)XY+YZ
\\ 
  16& 21& (\omega \lambda+1 ){X}^{2}+{Y}^{2}+  \bar\omega\lambda{Z}^{2}+\omega \lambda XY+ (\lambda+\bar\omega )YZ
\\ 
 9& 12& \bar\omega{Y}^{2}+ ( \lambda+\bar\omega)XY+\omega YZ+ (\lambda+1 )ZX
\\ 
 9& 15&\lambda {Z}^{2}+ ( \bar\omega\lambda+\bar\omega)XY+\omega \lambda YZ+ ( \lambda+\bar\omega)ZX
\\ 
  9& 19&{Y}^{2}+  \bar\omega\lambda{Z}^{2}+ \bar\omega XY+ ( \lambda+\bar\omega)YZ+\lambda ZX
\\ 
  12& 15& ( \lambda+\bar\omega){X}^{2}+ \bar\omega XY+ ( \lambda+\bar\omega)YZ+\lambda ZX
\\ 
  12& 19& (\omega \lambda+1 ){X}^{2}+\omega {Y}^{2}+ \bar\omega (\lambda+1 ) XY+\omega \lambda YZ+
 ( \lambda+\bar\omega)ZX
\\ 
  15& 19& ( \bar\omega\lambda+\omega ){X}^{2}+\omega \lambda {Z}^{2}+ ( \lambda+\bar\omega)XY+\omega YZ+ (\lambda+1 )ZX
\\ 
  4& 11& (\omega \lambda+1 ){X}^{2}+{Y}^{2}+\omega \lambda XY+YZ+ ( \lambda+\bar\omega)ZX
\\ 
  4& 15& ( \lambda+\bar\omega){X}^{2}+ ( \lambda+\bar\omega)XY+ \bar\omega (\lambda+1 )YZ+\lambda ZX
\\ 
  4& 18&\omega {Y}^{2}+\omega XY+\omega \lambda YZ+ (\lambda+1 )ZX
\\ 
  11& 15& ( \bar\omega\lambda+\omega ){X}^{2}+\omega \lambda {Z}^{2}+\omega XY+\omega \lambda YZ+ (\lambda+1)ZX
\\ 
  11& 18& \bar\omega{Y}^{2}+  \bar\omega\lambda{Z}^{2}+ ( \lambda+\bar\omega)XY+ \bar\omega (\lambda+1 )YZ+\lambda ZX
\\ 
  15& 18&\lambda {Z}^{2}+\omega \lambda XY+YZ+ ( \lambda+\bar\omega)ZX
\\ 
  7& 9& (\lambda+1 ){Z}^{2}+  \bar\omega\lambda XY+ \omega (\lambda+1) YZ+ (\lambda+1 )ZX
\\ 
  7& 11& ( \bar\omega\lambda+\omega ){X}^{2}+ (\omega \lambda+1 ){Z}^{2}+\omega XY+ (\omega \lambda+1 )YZ+ ( \lambda+\bar\omega)ZX
\\ 
  7& 21& (\omega \lambda+1 ){X}^{2}+\omega \lambda {Z}^{2}+ (\omega \lambda+1 )XY+\omega YZ+ZX
\\ 
  9& 11&{Y}^{2}+ \bar\omega{Z}^{2}+ (\omega \lambda+1 )XY+\omega YZ+ZX
\\ 
  9& 21&\omega {Y}^{2}+\lambda {Z}^{2}+\omega XY+ (\omega \lambda+1 )YZ+ ( \lambda+\bar\omega)ZX
\\ 
  11& 21& ( \bar\omega\lambda+\omega ){X}^{2}+\omega {Y}^{2}+\omega \lambda {Z}^{2}+  \bar\omega\lambda XY+ \omega (\lambda+1) YZ+
(\lambda+ 1 )ZX
\\ 
  7& 12& ( \lambda+\bar\omega){X}^{2}+ \bar\omega XY+ (\lambda+1 )YZ+ ( \lambda+\bar\omega)ZX
\\ 
  7& 14& (\omega \lambda+1 ){X}^{2}+\omega \lambda {Z}^{2}+ ( \bar\omega\lambda+\omega )XY+\omega \lambda YZ+ZX
\\ 
  7& 18&\lambda {Z}^{2}+ \omega (\lambda+1) XY+YZ+ \lambda ZX
\\ 
  12& 14& ( \lambda+\bar\omega){X}^{2}+{Y}^{2}+ \omega (\lambda+1) XY+YZ+\lambda ZX
\\ 
  12& 18&\omega {Y}^{2}+ ( \bar\omega\lambda+\omega )XY+\omega \lambda YZ+ZX
\\ 
  14& 18&{Y}^{2}+\lambda {Z}^{2}+ \bar\omega XY+ (\lambda+1)YZ+ ( \lambda+\bar\omega)ZX
\\ 
  4& 7& (\lambda+1 ){X}^{2}+ (\lambda+1 )XY+ (\bar\omega\lambda+\omega )YZ+ (\lambda+1 )ZX
\\ 
  4& 16& ( \lambda+\bar\omega){X}^{2}+ \bar\omega{Y}^{2}+\lambda XY+\omega YZ+\lambda ZX
\\ 
  4& 19&\omega {X}^{2}+ \bar\omega{Y}^{2}+XY+  \bar\omega\lambda YZ+ZX
\\ 
  7& 16& (\omega \lambda+1 ){X}^{2}+\omega \lambda {Z}^{2}+XY+  \bar\omega\lambda YZ+ZX
\\ 
  7& 19&  \bar\omega\lambda{X}^{2}+\omega \lambda {Z}^{2}+\lambda XY+\omega YZ+\lambda ZX
\\ 
  16& 19& ( \bar\omega\lambda+\omega ){X}^{2}+ \bar\omega{Y}^{2}+\omega \lambda {Z}^{2}+ (\lambda+1)XY+
 ( \bar\omega\lambda+\omega
)YZ+(\lambda+1 )ZX
\\ 
\end{array}
$$
\caption{Defining equations of  the conic curves  of type II}
\label{table:coniccurvesIIB}
\end{table}
}
\par
\medskip
Conversely,  let $\gamma\sprime$ be an arbitrary element of $\GGG\sb{B}$.
We will show the following:
\begin{claim}\label{claim:tripleB}
There exists a unique triple
$$
(g, t, \lambda)\;\in\; 
\PGL(3, k) \times \widetilde{T}\times (k\sm \{ 0, 1, \omega,\bar{\omega} \})
$$
such that $g\circ \gamma\sprime \circ t = \gamma\sb{\lambda}$ holds.
\end{claim}
The points 
$\gamma\sprime (P\sb{15})$, $\gamma\sprime (P\sb{16})$, $\gamma\sprime (P\sb{21})$
of $\gamma\sprime (\ell\sb{01})$ 
are on a line,
and the points
$\gamma\sprime (P\sb{12})$, $\gamma\sprime (P\sb{16})$, $\gamma \sprime (P\sb{20})$
of $\gamma\sprime (\ell\sb{10})$ 
are on another line.
Hence there exists a unique element $g\in \PGL (3, k)$ such that
$\gamma :=g\circ \gamma\sprime$ satisfies the following:
\begin{eqnarray}\label{eq:gammaPFIXB}
&& \gamma(P\sb{16})=[1,0,0],\nonumber\\
&&\gamma (P\sb{12})=[1,1,0],\;\;
\gamma(P\sb{20})=[0,1,0],\\
&&\gamma(P\sb{15})=[1,0,1],\;\;
\gamma(P\sb{21})=[0,0,1]. \nonumber
\end{eqnarray}
For $aa\sprime\in P(\AF)$,
let $L\sb{aa\sprime}\st \Pt$ be the line containing $\gamma (l\sb{aa\sprime})$,
and let
$$
\xi \sb{aa\sprime} X + \eta \sb{aa\sprime} Y +\zeta \sb{aa\sprime} Z =0 
$$
be the defining equation of $L\sb{aa\sprime}$.
By~\eqref{eq:gammaPFIXB},
we can put
\begin{eqnarray*}
 \xi\sb{00}=1, & \eta\sb{00}=0, & \zeta\sb{00}=0,\;\\
 \xi\sb{01}=0, & \eta\sb{01}=1, & \zeta\sb{01}=0,\;\\
 \xi\sb{02}=1, & \phantom{\eta\sb{02}=1,  } & \zeta\sb{02}=0,\;\\
 \xi\sb{10}=0, & \eta\sb{10}=0,& \zeta\sb{10}=1,\; \\
 \xi\sb{11}=1, & \eta\sb{11}=1, & \zeta\sb{11}=1,\; \\
 \xi\sb{12}=1, & \eta\sb{12}=1, & \phantom{\zeta\sb{12}=1,}\; \\
 \xi\sb{20}=1, & \eta\sb{20}=0, & \phantom{\zeta\sb{20}=1,}\; \\
 \xi\sb{21}=1, & \phantom{\eta\sb{21}=0, } & \zeta\sb{21}=1,\; \\
 \xi\sb{22}=0, & \eta\sb{22}=1.  & \phantom{\zeta\sb{22}=1.}\; 
\end{eqnarray*}
The three lines $L\sb{aa\sprime}$, $L\sb{bb\sprime}$, $L\sb{cc\sprime}$ 
are concurrent 
if  $\{aa\sprime, bb\sprime, cc\sprime\}\in L(\AF)$.
Hence we obtain a system of equations
\begin{equation}\label{eq:deteq}
\det\left[
\begin{array}{ccc}
\xi\sb{aa\sprime} & \eta\sb{aa\sprime}  &\zeta\sb{aa\sprime} \\
\xi\sb{bb\sprime} & \eta\sb{bb\sprime}  & \zeta\sb{bb\sprime} \\
\xi\sb{cc\sprime} & \eta\sb{cc\sprime}  & \zeta\sb{cc\sprime} 
\end{array}
\right]=0
\qquad\textrm{for every  $\{aa\sprime, bb\sprime, cc\sprime\}\in L(\AF)$.}
\end{equation}
A Gr\"obner basis of the ideal generated by the left hand side of~\eqref{eq:deteq}
in the polynomial ring 
$k[\eta\sb{02},  \eta\sb{21}, \zeta\sb{12}, \zeta\sb{20}, \zeta\sb{22}]$
is calculated as follows:
$$
\langle\;\; 
1+\zeta\sb{22}, 
\;\;
1+\zeta\sb{20}+\eta\sb{21},
\;\;
1+\zeta\sb{12}+\eta\sb{21},
\;\;
\eta\sb{02}+\eta\sb{21},
\;\;
1+\eta\sb{21}+\eta\sb{21}^2
\;\;\rangle.
$$
Hence there are two solutions of this system of equations,
\begin{eqnarray*}
&&\eta\sb{21}=\eta\sb{02}=\omega, \quad \zeta\sb{12}=\zeta\sb{20}=\bar\omega, \quad \zeta\sb{22}=1, 
\quad\textrm{or}\\
&&\eta\sb{21}=\eta\sb{02}=\bar\omega, \quad \zeta\sb{12}=\zeta\sb{20}=\omega, \quad \zeta\sb{22}=1, 
\end{eqnarray*}
which are conjugate over $\F\sb 2$.
If the latter holds, then 
we replace   $\gamma$ with $g\sb 0 \circ \gamma\circ T$,
where
$$
g\sb 0:=\left[
\begin{array}{ccc}
1 & 0  &0 \\
0 & 0 & 1 \\
0 & 1  & 0
\end{array}
\right],
$$
so that we can assume the former always holds.
The image of the $T$-points by $\gamma$ is therefore equal to 
the ones given in  Table~\ref{table:gammalambdaB},
and the lines $L\sb{aa\sprime}$ are given by equations in Table~\ref{table:linesB}.
\par
We next determine the coordinates  of the image of $C$-points by $\gamma$.
The point $\gamma (C(00))= \gamma (P\sb{17})$
is on the line $L\sb{00}=\{X=0\}$
and is different from $\gamma(P\sb{20})=[0,1,0]$ and $\gamma(P\sb{21})=[0,0,1]$.
Hence we can put
\begin{equation}\label{eq:17B}
\gamma (P\sb{17})=[0, \lambda, 1],
\end{equation}
where $\lambda$ is a parameter $\ne 0$.
Let $\ell, \ell\sb 1, \ell\sb 2$ be three distinct elements of $L(\AF)$
that are parallel to each other.
The conic curves $Q$ satisfying the following conditions form a pencil $PQ\sb{\ell}$:
\begin{itemize}
\item[(i)] $Q$ contains  $\gamma (T(\ell\sb 1))$ and $\gamma (T(\ell\sb 2))$, and
\item[(ii)] $Q$ is tangent to the lines $L\sb{aa\sprime}$, $L\sb{bb\sprime}$, $L\sb{cc\sprime}$,
where $\ell=\{aa\sprime, bb\sprime, cc\sprime\}$.
(Recall Definition~\ref{def:tangent} and Remark~\ref{rem:tangent}.)
\end{itemize}
Using the coordinates of the points $\gamma (T(\ell))$ and the defining equations of the nine lines  $L\sb{aa\sprime}$
determined so far,
we can calculate this pencil explicitly.
By Proposition~\ref{prop:B1},
the conic curve $Q\sb{\ell}$ containing  $\gamma(q\sb{\ell})$
is a nonsingular member of the pencil $PQ\sb{\ell}$.
Starting from~\eqref{eq:17B},
we can determine the coordinates of $\gamma (C(aa\sprime))$, 
and see that they coincide with Table~\ref{table:gammalambdaB}.
For example,
consider  $\ell=\{01,10,22\}\in L(\AF)$.
We have 
$$
\ell\sb 1=\{ 02,11,20\},
\quad
\ell\sb 2=\{00, 12, 21\}.
$$
The pencil of conic curves passing through the points 
$$
\gamma( T(\ell\sb 1) )=\gamma(P\sb 9)=[1, \bar\omega,\omega],
\quad
\gamma(T(\ell\sb 2))=\gamma(P\sb{18})=[0, 1,\omega],
$$
and tangent to the lines
$$
L\sb{01}=\{Y=0\}, 
\quad
L\sb{10}=\{Z=0\},
\quad  
L\sb{22}=\{ Y+Z=0\}
$$
is spanned by the two conic curves 
defined by
$$
\omega X^2+ \bar\omega Y^2 +Z^2=0
\quand
\omega X^2+\omega Y^2 +YZ=0.
$$
Because the conic curve  $Q\sb\ell$ passes through $\gamma (P\sb{17})=[0, \lambda, 1]$, it is defined by
$$
\lambda (\omega X^2+ \bar\omega Y^2 +Z^2)
+( \omega \lambda+1)(\omega X^2+\omega Y^2 +YZ)=0.
$$
The intersection points of $Q\sb\ell$  with
the line $ L \sb{12}=\{\omega X +\omega Y +Z=0\}$
are $\gamma\ (T(\ell\sb 2))=\gamma(P\sb{18})=[0, 1, \omega]$
and $\gamma (C(12))=\gamma(P\sb 6)$.
Hence we obtain 
$$
\gamma (C(12))=\gamma(P\sb 6)= [\omega \lambda +1, \lambda+1, \lambda].
$$
See Table~\ref{table:PQ}
for the detail of the calculation.
\begin{table}
$$
\renewcommand{\arraystretch}{1.3}
\begin{array}{c|c|c}
\ell & F\sb 1,  F\sb 2 & \beta\sb{\ell} \\
\hline
00, 01, 02 & 
\omega X^2+\omega Y^2 +Z^2, \;\;\; \bar\omega X^2 +\omega Y^2 +XY &
(\lambda^2+\bar\omega)/\lambda^2
\\
00, 10, 20 &
X^2+\omega Y^2 +Z^2, \;\;\;  X^2 +Y^2 +ZX &
(\omega\lambda^2+1)/\lambda^2 
\\
00, 11, 22 &
\bar\omega X^2+\omega Y^2 +Z^2, \;\;\;  \omega X^2 +\bar\omega Y^2 +XY +ZX & 
(\bar\omega \lambda^2+\omega )/\lambda^2 
\\
00, 12, 21 &
\omega Y^2 + Z^2, \;\;\;  \omega Y^2 +\omega XY +ZX &
(\lambda^2+\bar\omega)/\lambda^2
\\
01, 10, 22 &
\omega X^2 +\bar\omega Y^2 +Z^2,\;\;\; \omega X^2 +\omega Y^2 + YZ &
(\omega\lambda +1)/\lambda 
\\
01, 11, 21 &
\omega X^2 + Z^2,\;\;\;  X^2 +XY +YZ &
1/\lambda 
\\
01, 12, 20 &
\omega X^2 + Y^2 + Z^2,\;\;\;\bar\omega X^2 +Y^2 +\omega X Y + YZ &
(\lambda+1)/\lambda
\\
02, 10, 21 &
\bar\omega X^2 + Y^2 +Z^2,\;\;\; \omega X^2 +\omega Y^2 +\omega YZ +ZX &
(\bar\omega\lambda +\bar\omega)/\lambda
\\
02, 11, 20 &
\bar\omega Y^2 +Z^2, \;\;\; \bar\omega Y^2 +\bar\omega XY +\omega YZ +ZX &
(\lambda+\bar\omega)/\lambda
\\
02, 12, 22 &
X^2+Z^2, \;\;\; X^2 +XY +\omega YZ +ZX &
\bar\omega /\lambda
\\
10, 11, 12 &
\bar\omega X^2 +Y^2, \;\;\; X^2 +YZ+ZX &
\lambda
\\
20, 21, 22 &
X^2+Y^2, \;\;\; X^2 +XY +\bar\omega YZ +ZX &
\omega\lambda
\end{array}
$$
\vskip 5pt
\caption{Basis $\{F\sb 1=0,  F\sb 2=0\}$ of the pencil $PQ\sb{\ell}$ and the member
$Q\sb{\ell}=\{ F\sb 1 +\beta\sb{\ell} F\sb 2=0\}$}\label{table:PQ}
\end{table}
Thus we have proved that $\gamma$ is equal to $\gamma\sb{\lambda}$.
Because $\gamma\sb{\lambda}$ is injective,
$\lambda$ is not among $\{0,1, \bar\omega\}$.

There exists a unique conic curve containing $\gamma\sb{\lambda} (q)$
for each quadratic word $q$ of $\CB$,
and the defining equations of those conic curves  are given in Tables~\ref{table:coniccurvesIB} and~\ref{table:coniccurvesIIB}.
The smoothness of these curves implies that $\lambda\ne \omega$.
Thus we have proved Claim~\ref{claim:tripleB}.
\end{proof}
\begin{remark}
The polynomial $\GB[\lambda]$ is the defining equation of the nodal splitting curve
$$
L\sb{00} \cup L\sb{01} \cup L\sb{10} \cup L\sb{11} \cup Q\sprime \sb{\ell, \ell\sprime},
$$
where $T(\ell)=P\sb{16}$ and $T(\ell\sprime)=P\sb{19}$.
See Proposition~\ref{prop:nodalsplit}.
\end{remark}
\begin{remark}
Consider the projective plane $(\Pt)\dual$ of lines on $\Pt$.
Let $[U, V, W]$ be the homogeneous coordinates of $(\Pt)\dual$ 
dual to the homogeneous coordinates $[X,Y,Z]$ of $\Pt$.
Let $E\sb{\lambda}$ be the cubic curve in $(\Pt)\dual$ defined by
\begin{multline*}
\bar\omega\, U^2 W + U W^2+\omega\lambda\, UV^2+(\omega\lambda+1) \,V^2 W+\\
+(\omega \lambda +1) \,V W^2 +\bar\omega \lambda\, U^2 V +(\omega+\lambda ) \,UVW=0.
\end{multline*}
Then the points $\gamma\sb{\lambda} (C(P(\AF)))$ correspond 
to the nine flex tangents to $E\sb{\lambda}$,
and 
the points $\gamma\sb{\lambda} (T(L(\AF)))$ correspond to the 
twelve lines 
containing three flex points of $E\sb{\lambda}$.
\end{remark}
\begin{remark}\label{rem:limitDKB}
When $\lambda=\omega$,
the set $\gamma\sb\lambda (\Ps)$ coincides with $\Pt (\Ff)$,
and the point $[\GB[\lambda]]\in\moduli$ is equal to the Dolgachev-Kondo point.
\end{remark}
For each $\sigma\in \Aut (\CB)$,
we calculate the unique triple
$$
(g\sb{\sigma}, t\sb{\sigma}, \lambda\sp{\sigma})\;\in\; \PGL (3, k(\lambda)) \times \widetilde{T}\times k(\lambda)
$$
such that 
$g\sb{\sigma} \circ (\gamma\sb{\lambda}\circ \sigma)\circ t\sb{\sigma}=\gamma\sb{\lambda\sp{\sigma}}$
holds.
The map $\sigma\mapsto t\sb{\sigma}$ is a homomorphism from $\Aut (\CB)$ to $\widetilde T$.
We  put
$N\sb B:=\Ker (\Aut (\CB)\to \widetilde{T} )$.
\begin{corollary}\label{cor:connectedcompsB}
The space $\PGL(3, k)\backslash \gs\sb B$ has exactly two connected components,
each of which is isomorphic to $\A\sp 1\sm \{0,1,\omega,\bar\omega\}$.
One of them is given, set-theoretically, by
$$
(\PGL(3, k)\backslash \gs\sb B)\sp + :=\set{[\gamma\sb\alpha]}{\alpha\in k\sm \{0,1,\omega,\bar\omega\}},
$$
and the other one is equal to $((\PGL(3, k)\backslash \gs\sb B)\sp +)\cdot T$.
The group $N\sb B$ acts on $(\PGL(3, k)\backslash \gs\sb B)\sp +$, and
the moduli curve $\moduli\sb B$ is equal to the quotient space $(\PGL(3, k)\backslash \gs\sb B)\sp +/N\sb B$.
\end{corollary}
Consider the natural projection
$$
\map{p\sb B}{ \A\sp 1\sm \{0,1,\omega, \bar\omega\} \cong (\PGL (3, k) \backslash \GGG\sb{B})\sp + }{%
\moduli\sb{B} = (\PGL (3, k ) \backslash \GGG\sb{B})\sp + / N\sb B}.
$$
For $\alpha\in k\sm\{0,1,\omega,\bar\omega\}$,
let $P[\alpha]$ denote the point of $\A\sp 1\sm \{0,1,\omega, \bar\omega\}$
given by $\lambda=\alpha$.
Then $p\sb B (P[\alpha])\in\moduli\sb B$ corresponds to the isomorphism class 
of the polarized supersingular $K3$ surface $( X\sb{\GB [\alpha]}, \pol\sb{\GB [\alpha]})$.
The following is proved in the same way as Proposition~\ref{prop:groupA}.
\begin{proposition}\label{prop:groupB}
The fiber
$p\sb B\inv (p\sb B (P[\alpha]))$ is equal 
$\{ P[\varphi]\}$, where $\varphi$ runs through 
the set  $\Gamma\sb B$ in Theorem~\ref{thm:B}  with $\lambda$ replaced by $\alpha$.
The group  $\Aut( X\sb{\GB [\alpha]}, \pol\sb{\GB [\alpha]})$
is equal to the subgroup of $\PGL (3, k)$  generated by the elements in~\eqref{eq:groupB}.
\end{proposition}
\begin{corollary}
We have
$\moduli\sb{B}=\Spec k[J\sb B, 1/J\sb B]$,
where 
$$
J\sb B={(\lambda+\omega)^{12}}/{\lambda^3(\lambda+1)^3(\lambda+{\bar\omega})^3}.
$$
The morphism $p\sb B$ is an \'etale  Galois covering
with Galois group $\Gamma\sb B$,
which is isomorphic to the alternating group $\AAAA\sb 4$.
\end{corollary}
Indeed the group $\Gamma\sb B$ acts on the set $\{0, 1,  \bar\omega, \infty\}$ as $\AAAA\sb 4$.
\section{The moduli curve corresponding to the code $\CC$}\label{sec:C}
In this section,
we prove Theorem~\ref{thm:C}.
\par
\medskip
The linear words of $\CC$ are listed in Table~\ref{table:linesC}.
\begin{table}[tb]
$$
\begin{array}{ccccccccc}
l\sb{1} & := & \{ &17, &18, &19, &20, &21&\},  \\ 
l\sb{2} & := & \{ &13, &14, &15, &16, &21&\},  \\ 
l\sb{3} & := & \{ &9, &10, &11, &12, &21&\},  \\ 
l\sb{4} & := & \{ &5, &6, &7, &8, &21&\},  \\ 
l\sb{5} & := & \{ &1, &2, &3, &4, &21&\}.
\end{array}
$$
\caption{Linear words of $\CC$}\label{table:linesC}
\end{table}
The list of quadratic words in $\CC$ is omitted.
The point $P\sb{21}$ is special because every linear word contains it.
The following can be checked directly by computer:
\begin{proposition}
Let $\phi: \PPP \isom \P^2 (\F\sb 4)$ be the bijection given in Table~\ref{table:embedC}.

{\rm (1)} The linear words of $\CC$ are precisely the words $\phi\inv (\Lambda (\Ff))$,
where $\Lambda$ are  $\F\sb 4$-rational lines passing through 
$$
O:=[0,0,1]=\phi(P\sb{21}).
$$

{\rm (2)}
The  quadratic words of $\CC$ are  precisely the words
$\phi\inv (\Lambda (\F\sb 4) + \Lambda\sprime (\F\sb 4))$, 
where $\Lambda$ and $\Lambda\sprime$ are distinct $\F\sb 4$-rational lines 
that do not pass through  $O$.
\end{proposition}
\begin{table}
\par
\medskip
\noindent
{\phantom{a}
\hskip -6.5cm 
\hbox {
\hbox{
\hbox to 5cm{
\parbox{5cm}{
\begin{eqnarray*}
\phi(P\sb{1})&=& [1,1,0], \\
\phi(P\sb{2})&=& [1,1,1], \\
\phi(P\sb{3})&=& [1,1,\oo], \\
\phi(P\sb{4})&=& [1,1,\bo], \\
\phi(P\sb{5})&=& [1,\bo,\oo], \\
\phi(P\sb{6})&=& [1,\bo,0], \\
\phi(P\sb{7})&=& [1,\bo,1], \\
\phi(P\sb{8})&=& [1,\bo,\bo], \\
\phi(P\sb{9})&=& [1,\oo,\bo], \\
\phi(P\sb{10})&=& [1,\oo,0], \\
\phi(P\sb{11})&=& [1,\oo,\oo], \\
\phi(P\sb{12})&=& [1,\oo,1], 
\end{eqnarray*}
}
}
\hbox to 6cm{
\parbox{6cm}{
\begin{eqnarray*}
\phi(P\sb{13})&=& [1,0,1], \\
\phi(P\sb{14})&=& [1,0,0], \\
\phi(P\sb{15})&=& [1,0,\bo], \\
\phi(P\sb{16})&=& [1,0,\oo], \\
\phi(P\sb{17})&=& [0,1,1], \\
\phi(P\sb{18})&=& [0,1,0], \\
\phi(P\sb{19})&=& [0,1,\bo], \\
\phi(P\sb{20})&=& [0,1,\oo], \\
\phi(P\sb{21})&=& [0,0,1]. \\
\phantom{\phi(P\sb{21})}&& \\
\phantom{\phi(P\sb{21})}&& \\
\phantom{\phi(P\sb{21})}&& 
\end{eqnarray*}
}
}
}
}
}
\vskip 5pt
\caption{Bijection  $\phi$ from  $\PPP$ to $\P^2 (\F\sb 4)$}\label{table:embedC}.
\end{table}
Note that $\phi$ embeds $\CC$ into the Dolgachev-Kondo code $\abscode\sb{\DK}$.
\begin{corollary}\label{cor:unquedisjoint}
For each quadratic word $q$ in $\CC$,
there exists a unique linear word $l$ in $\CC$ such that $q\cap l=\emptyset$.
\end{corollary}
From Remark~\ref{rem:tetrad}, we obtain the following:
\begin{corollary}\label{cor:unquedisjoint2}
Let $l$ and $l\sprime$ be distinct linear words of  $\CC$,
and let $A\sb 1, A\sb 2\in l$
{\rm(}resp.~$B\sb 1, B\sb 2\in l\sprime${\rm)}
be distinct points
not equal to $P\sb{21}$.
Then there are exactly two quadratic words
$q$ and $q\sprime$ in $\CC$ containing the points $\{ A\sb 1, A\sb 2, B\sb 1, B\sb 2\}$.
Moreover,
if a linear word $l\spprime \in \CC$ is disjoint from $q$,
then $l\spprime$ is also disjoint from $q\sprime$.
\end{corollary}
For $\alpha\sb 1, \alpha\sb 2, \alpha\sb 3\in \Ff$,
we denote by $\Lambda [\alpha\sb 1\alpha\sb 2\alpha\sb 3]$ the $\Ff$-rational line defined by
$$
\alpha\sb 1 X +\alpha\sb 2 Y +\alpha\sb 3 Z=0,
$$
and by $q[\alpha\sb 1\alpha\sb 2\alpha\sb 3, \beta\sb1 \beta\sb 2 \beta\sb 3]\in\CC$
the quadratic word
$$
\phi\inv ( \Lambda [\alpha\sb 1\alpha\sb 2\alpha\sb 3] (\F\sb 4) + \Lambda [\beta\sb 1\beta\sb 2\beta\sb 3] (\F\sb 4)).
$$
We put
\begin{equation*}
LG\sprime:=\set{g \in \PGL (3, \F\sb 4)}{ g(O)=O}.
\end{equation*}
The automorphism group $\Aut (\CC)$
of the code  $\CC$ contains a subgroup 
$$
LG:=\phi\inv \circ LG\sprime\circ \phi.
$$
The order of $LG$ is $2880$.
The  group $\Aut (\CC)$ also contains 
the permutation
\begin{equation*}
T:= (P\sb{3} P\sb {4}) (P\sb{5} P\sb {9}) (P\sb{6} P\sb{10})(P\sb{7} P\sb {12})(P\sb{8} P\sb {11})(P\sb{15} P\sb {16})
(P\sb{19} P\sb {20})
\end{equation*}
of $\PPP$ that corresponds, via the bijection $\phi$,  
to the action of  the conjugation $\oo\mapsto \bo$ over $\F\sb 2$
on $\P^2 (\F\sb 4)$.
It can be checked easily by computer that the following  permutation is also 
contained in $\Aut(\CC)$:
$$
S := (P\sb{1}P\sb{3})(P\sb{2}P\sb{4})(P\sb{5}P\sb{7})(P\sb{6}P\sb{8})(P\sb{9}P\sb{11})(P\sb{10}P\sb{12}).
$$
The automorphisms $T$ and $S$ of $\CC$ generate a subgroup isomorphic to the dihedral group of order $8$
in $\Aut(\CC)$. 
An ordered quartet 
$$
(R\sb 1, R\sb 2, R\sprime\sb 1, R\sprime\sb 2)
$$
of points in $\PPP\sm \{P\sb{21}\}$ is called a \emph{marking quartet} 
if $P\sb{21}$, $R\sb 1$, $R\sb 2$ are in a linear word, and 
$P\sb{21}$, $R\sb 1\sprime$, $R\sb 2\sprime$ are in another  linear word.
There are $2880$ marking quartets,
and the action  of  $LG$  on the set of marking quartets is 
simply transitive.
\begin{proposition}
The group 
$\Aut (\CC)$ is generated by $LG$, $T$ and $S$,
and 
the order of $\Aut(\CC)$ is $23040$.
\end{proposition}
\begin{proof}
Let $\sigma$ be an arbitrary element of $\Aut (\CC)$.
Because   $(P\sb{17}, P\sb{18}, P\sb {13}, P\sb{14})$ 
and 
$(\sigma (P\sb{17}), \sigma(P\sb{18}), \sigma(P\sb{13}),\sigma(P\sb{14}))$   
are marking quartets,  there exists an element
$\tau\in LG$ such that
$\tau \sigma(P\sb{i})=P\sb{i}$ for $i=13, 14, 17, 18, 21$.
Because $\tau\sigma (l\sb 1)=l\sb 1$ and $\tau\sigma(l\sb 2)=l\sb 2$, 
we have 
$$
\{\tau \sigma(P\sb{19}), \tau \sigma(P\sb{20})\}=\{P\sb{19}, P\sb{20}\}
\quand
\{\tau \sigma(P\sb{15}), \tau \sigma(P\sb{16})\}=\{P\sb{15}, P\sb{16}\}.
$$
If $\tau \sigma(P\sb{19})=P\sb{20}$,
then we replace $\tau$ by $T \tau$.
Therefore,
modulo the subgroup generated by $LG$ and $T$,
we can assume that $\sigma$ has the following properties:
\begin{itemize}
\item[($\sigma$-i)] $\sigma$ fixes each of the seven points 
$P\sb {13}, P\sb{14}, P\sb{17}, P\sb{18}, P\sb{19}, P\sb{20}, P\sb{21}$, 
\item[($\sigma$-ii)] $\{ \sigma(P\sb{15}),  \sigma(P\sb{16})\}=\{P\sb{15}, P\sb{16}\}$
\end{itemize}
Consider, for example,
a set of four points $\{ P\sb{13}, P\sb{14}, P\sb{17}, P\sb{18}\}$,
each of which is fixed by $\sigma$.
The two quadratic words containing them are 
\begin{eqnarray*}
q[101, 011] &=&\{2, 7, 12, 13, 18\}+\{2, 8, 11, 14, 17\}=\{7,8,11,12,13,14,17,18\}, \\
q[111, 001] &=&\{1, 5, 9, 13, 17\}+ \{1, 6, 10, 14, 18\}=\{5,6,9,10,13,14,17,18\}.
\end{eqnarray*}
Both of $q[101, 011]$ and $q[111, 001]$ are disjoint from $l\sb 5$.
By Corollary~\ref{cor:unquedisjoint},
we have $\sigma(l\sb 5)=l\sb 5$.
Considering other sets of four points fixed by $\sigma$,
we can show that $\sigma(l\sb 4)=l\sb 4$ and $\sigma(l\sb 3)=l\sb 3$.
In Table~\ref{table:fqql},
we list the triples $(f, \{q, q\sprime\}, l\sb{\nu})$,
where  $f$ is a set of four points pointwise fixed by $\sigma$,
$\{q, q\sprime\}$ is the pair of quadratic words containing $f$, 
and $l\sb\nu$ is the linear word  disjoint from both of $q$ and $q\sprime$.
\begin{table}
$$
\renewcommand{\arraystretch}{1.3}
\begin{array}{c|c|c}
f & [q, q\sprime] & l\sb{\nu} \\
\hline
\{13,14,17,18\} & q[101,011]=\{ {7, 8, 11, 12, 13, 14, 17, 18}     \} &l\sb 5 \\
															 & q[111,001]=\{   {5, 6, 9, 10, 13, 14, 17, 18}   \} & \\
\hline
\{13,14,17,19\} & q[1\bar\omega 1,011]=\{  {2, 3, 10, 11, 13, 14, 17, 19}    \} &l\sb 4 \\
															 & q[111,0\bar\omega 1]=\{  {1, 4, 9, 12, 13, 14, 17, 19}    \} & \\
\hline
\{13,14,17,20\} & q[1\omega 1,011]=\{   {2, 4, 6, 8, 13, 14, 17, 20}   \} &l\sb 3 \\
															 & q[111,0\omega 1]=\{    {1, 3, 5, 7, 13, 14, 17, 20}  \} & \\
\hline
\{13,14,18,19\} & q[101,0\bar\omega 1]=\{   {2, 4, 5, 7, 13, 14, 18, 19}   \} &l\sb 3 \\
															 & q[1\bar\omega 1,001]=\{  {1, 3, 6, 8, 13, 14, 18, 19}    \} & \\
\hline
\{13,14,18,20\} & q[101,0\omega 1]=\{ {2, 3, 9, 12, 13, 14, 18, 20}     \} &l\sb 4 \\
															 & q[1\omega 1,001]=\{  {1, 4, 10, 11, 13, 14, 18, 20}   \} & \\
\hline
\{13,14,19,20\} & q[1\bar\omega 1,0\omega 1]=\{  {7, 8, 9, 10, 13, 14, 19, 20}    \} &l\sb 5 \\
															 & q[1\omega 1,0\bar\omega 1]=\{   {5, 6, 11, 12, 13, 14, 19, 20}   \} & \\

\end{array}
$$
\vskip 5pt
\caption{List of the triples $(f, \{q, q\sprime\}, l\sb{\nu})$}\label{table:fqql}
\end{table}
Therefore we have the following:
\begin{itemize}
\item[($\sigma$-iii)] $\sigma$ leaves each of the sets
$$
\{ P\sb{1}, P\sb{2}, P\sb{3}, P\sb{4} \},
\;\;
\{  P\sb{5}, P\sb{6}, P\sb{7}, P\sb{8} \}, 
\;\;
\{  P\sb{9}, P\sb{10}, P\sb{11}, P\sb{12}\}
$$
invariant.
\end{itemize}
Let us consider the quadratic words  $q\sb 1:=q[101, 011]$ and $q\sb 2:=q[111, 001] $ again.
Since
$$
\{ \sigma (q\sb 1 \cap l\sb 4), \sigma (q\sb 2 \cap l\sb 4)\}
=\{q\sb 1 \cap l\sb 4, q\sb 2 \cap l\sb 4\}, 
$$
the action of $\sigma$ on $\{P\sb 5, P\sb 6, P\sb 7, P\sb 8\}$ preserves the decomposition
$$
\{P\sb 5, P\sb 6, P\sb 7, P\sb 8\}=\{P\sb 5, P\sb 6\}\cup \{P\sb 7, P\sb 8\};
$$
that is, $\{\sigma(P\sb 5), \sigma (P\sb 6)\}$ is either $\{P\sb 5, P\sb 6\}$ or $\{P\sb 7, P\sb 8\}$.
By the same argument applied to  the pairs $\{q, q\sprime\}$ of quadratic words in Table~\ref{table:fqql},
we see  the following:
\begin{itemize}
\item[($\sigma$-iv)] $\sigma$ preserves the decompositions
\begin{eqnarray*}
\{P\sb 1, P\sb 2, P\sb 3, P\sb 4\}&=&\{P\sb 1, P\sb 4\}\cup \{P\sb 2, P\sb 3\}=\{P\sb 1, P\sb 3\}\cup \{P\sb 2, P\sb 4\}, \\
\{P\sb 5, P\sb 6, P\sb 7, P\sb 8\}&=&\{P\sb 5, P\sb 6\}\cup \{P\sb 7, P\sb 8\}=\{P\sb 5, P\sb 7\}\cup \{P\sb 6, P\sb 8\}, \quand
\\
\{P\sb 9, P\sb {10}, P\sb {11}, P\sb {12}\}&=&\{P\sb {9}, P\sb {10}\}\cup \{P\sb {11}, P\sb {12}\}=
\{P\sb {9}, P\sb {12}\}\cup\{P\sb {10}, P\sb {11}\}.
\end{eqnarray*}
\end{itemize}
The two quadratic words
containing $\{P\sb{13}, P\sb{16}, P\sb{17}, P\sb{18}\}$ are 
\begin{eqnarray*}
q[\omega 11, 101]&=&\{ {2, 4, 10, 12, 13, 16, 17, 18} \}\quand \\
q[111, \omega 01]&=&\{ {1, 3, 9, 11, 13, 16, 17, 18} \},
\end{eqnarray*}
both of which are disjoint from $l\sb 4$.
On the other hand, 
the two quadratic words
containing $\{P\sb{13}, P\sb{15}, P\sb{17}, P\sb{18}\}$ are
\begin{eqnarray*}
q[\bar\omega 11, 101]&=&\{ {2, 3, 6, 7, 13, 15, 17, 18} \} \quand \\
q[111, \bar\omega 01]&=&\{ {1, 4, 5, 8, 13, 15, 17, 18} \},
\end{eqnarray*}
both of which 
are disjoint from $l\sb 3$.
Since $\sigma$ fixes each of $l\sb 4$ and $l\sb 3$,
we see that 
the property ($\sigma$-ii) of $\sigma$ can be  strengthened to the following:
\begin{itemize}
\item[($\sigma$-ii)${}\sprime$] $\sigma(P\sb{15})=P\sb{15}$, $\sigma(P\sb{16})=P\sb{16}$.
\end{itemize}
Using computer,
we can easily list  all $4^3=64$ permutations  $\sigma$
satisfying ($\sigma$-i), ($\sigma$-ii)${}\sprime$, ($\sigma$-iii) and ($\sigma$-iv).
We can  check that exactly four  of them $\id$, $S$,
\begin{eqnarray*}
(ST)^2 &=&(P\sb{1}P\sb{2})(P\sb{3}P\sb{4})(P\sb{5}P\sb{6})(P\sb{7}P\sb{8})(P\sb{9}P\sb{10})(P\sb{11}P\sb{12}) \quand \\
(ST)^2 S &=& (P\sb{1}P\sb{4})(P\sb{2}P\sb{3})(P\sb{5}P\sb{8})(P\sb{6}P\sb{7})(P\sb{9}P\sb{12})(P\sb{10}P\sb{11})
\end{eqnarray*}
preserve the set of quadratic words in $\CC$.
Hence, by Proposition~\ref{prop:autcodeW},  $\Aut (\CC)$ is generated by $LG$, $S$ and $T$.
It can be checked by computer that
the order of $\Aut (\CC)$ is $23040$.
\end{proof}
\begin{table}
\par
\medskip
\noindent
{\phantom{a}
\hskip -6.5cm 
\hbox {
\hbox{
\hbox to 5cm{
\parbox{5cm}{
\begin{eqnarray*}
\gamma\sb{\lambda  } (P\sb{1}) &=&  [1,1,\lambda  ],\\
\gamma\sb{\lambda  } (P\sb{2}) &=&  [1,1, \lambda  +1], \\
\gamma\sb{\lambda  } (P\sb{3}) &=&  [1,1,\lambda  +\oo],\\
\gamma\sb{\lambda  } (P\sb{4}) &=&  [1,1,\lambda  +\bo],\\
\gamma\sb{\lambda  } (P\sb{5}) &=&  [1, \bo, \oo\lambda  +\oo],\\
\gamma\sb{\lambda  } (P\sb{6}) &=&  [1, \bo, \oo\lambda  ], \\
\gamma\sb{\lambda  } (P\sb{7}) &=&  [1, \bo, \oo\lambda  +1],\\
\gamma\sb{\lambda  } (P\sb{8}) &=&  [1, \bo, \oo\lambda  +\bo],\\
\gamma\sb{\lambda  } (P\sb{9}) &=&  [1, \oo, \bo\lambda  +\bo],\\
\gamma\sb{\lambda  } (P\sb{10}) &=&  [1, \oo,\bo\lambda  ] , \\
\gamma\sb{\lambda  } (P\sb{11}) &=&  [1, \oo, \bo\lambda  +\oo],\\
\gamma\sb{\lambda  } (P\sb{12}) &=&  [1, \oo, \bo\lambda  +1],
\end{eqnarray*}
}
}
\hbox to 6cm{
\parbox{6cm}{
\begin{eqnarray*}
\gamma\sb{\lambda  } (P\sb{13})&=&[1,0,1],\\
\gamma\sb{\lambda  } (P\sb{14}) &=& [1,0,0],\\
\gamma\sb{\lambda  } (P\sb{15}) &=&  [1, 0, \bo],\\
\gamma\sb{\lambda  } (P\sb{16}) &=&  [1, 0, \oo] , \\
\gamma\sb{\lambda  } (P\sb{17}) &=&[0,1,1],\\
\gamma\sb{\lambda  } (P\sb{18})&=&[0,1,0], \\
\gamma\sb{\lambda  } (P\sb{19}) &=& [0, 1, \bo], \\
\gamma\sb{\lambda  } (P\sb{20}) &=& [0, 1, \oo], \\
\gamma\sb{\lambda  } (P\sb{21})&=&[0,0,1]. \\
\phantom{\gamma\sb{\lambda  } (\phi(P\sb{21}))}&&\\
\phantom{\gamma\sb{\lambda  } (\phi(P\sb{21}))}&&\\
\phantom{\gamma\sb{\lambda  } (\phi(P\sb{21}))}&&
\end{eqnarray*}
}
}
}
}
}
\vskip 5pt
\caption{Definition of  $\gamma\sb{\lambda  }$ for $\CC$}\label{table:gammalambdaC}
\end{table}
Let $\lambda$ be a parameter of the affine line $\A\sp 1$.
We define 
$\gamma\sb{\lambda} : \Ps\to \Pt$
by Table~\ref{table:gammalambdaC}.
When $\lambda=0$,
the map $\gamma\sb\lambda$ is equal to $\phi$.
Let $\widetilde T$ denote the subgroup of $\Aut (\CC)$ generated by the involution $T$.
\begin{proposition}
The map $\lambda\mapsto \gamma\sb{\lambda}$ induces an isomorphism   
from $\A\sp 1\sm \{0,1, \omega, \bar\omega\}$ to $\PGL(3, k)\backslash \GGG\sb{C}/ \widetilde T $.
\end{proposition}
\begin{proof}
First note that $\gamma\sb{\lambda}$ is injective for every $\lambda$.
\begin{claim}\label{claim:tripleC}
Let $\gamma\sprime$ be an arbitrary element of $\gs\sb{C}$.
Then 
 there exists a unique triple
$$
(g, t, \lambda)\;\in\; 
\PGL(3, k) \times \widetilde{T}\times (k\sm \{ 0, 1, \omega,\bar{\omega} \})
$$
such that
$g\circ \gamma\sprime \circ t = \gamma\sb{\lambda}$.
\end{claim}
Since $\gamma\sprime (P\sb{21})$,  $\gamma\sprime (P\sb{13})$, $\gamma\sprime (P\sb{14})$
are on a line,
and 
$\gamma\sprime (P\sb{21})$, $\gamma\sprime (P\sb{17})$, $\gamma\sprime (P\sb{18})$ are on another line,
there exists a unique $g\in \PGL(3, k)$
such that $\gamma:=g\circ \gamma\sprime$ satisfies 
\begin{eqnarray*}
&& \gamma (P\sb{21})=[0,0,1]=O,\\
&& \gamma (P\sb{17})=[0,1,1],\quad
\gamma (P\sb{18})=[0,1,0],\quad \\
&&\gamma (P\sb{13})=[1,0,1],\quad
\gamma (P\sb{14})=[1,0,0].
\end{eqnarray*}
The $X$-coordinate of $\gamma (P\sb i)$ is not  $0$
for $i=1, \dots, 16$,
because otherwise
$\gamma (P\sb i)$, $\gamma (P\sb{17})$ and $\gamma(P\sb{21})$ would be collinear,
and hence there would exist a linear word of $\CC$ containing $\{P\sb i, P\sb{17}, P\sb{21}\}$ by 
Proposition~\ref{prop:linewt}.
Therefore
there exist  parameters
$\alpha\sb 1, \alpha\sb 2$, $\beta\sb 1, \beta\sb 2$,
$t\sb i$, $s\sb{ij}$ ($i=3,4,5$, $j=1, \dots, 4$) such that 
$\gamma$ is given by Table~\ref{table:parametricembedC}.
\begin{table}
\par
\medskip
\noindent
{\phantom{a}
\hskip -6.5cm 
\hbox {
\hbox{
\hbox to 5cm{
\parbox{5cm}{
\begin{eqnarray*}
\gamma (P\sb{1}) &=&  [1, t\sb{5},  s\sb{5,1}],\\
\gamma (P\sb{2}) &=&  [1, t\sb{5},  s\sb{5,2}], \\
\gamma (P\sb{3}) &=&  [1, t\sb{5},  s\sb{5,3}],\\
\gamma (P\sb{4}) &=&  [1, t\sb{5},  s\sb{5,4}],\\
\gamma (P\sb{5}) &=&  [1, t\sb{4},  s\sb{4,1}],\\
\gamma (P\sb{6}) &=&  [1, t\sb{4},  s\sb{4,2}], \\
\gamma (P\sb{7}) &=&  [1, t\sb{4},  s\sb{4,3}],\\
\gamma (P\sb{8}) &=&  [1, t\sb{4},  s\sb{4,4}],\\
\gamma (P\sb{9}) &=&  [1, t\sb{3},  s\sb{3,1}],\\
\gamma (P\sb{10}) &=&  [1, t\sb{3},  s\sb{3,2}] , \\
\gamma (P\sb{11}) &=&  [1, t\sb{3},  s\sb{3,3}],\\
\gamma (P\sb{12}) &=&  [1, t\sb{3},  s\sb{3,4}],
\end{eqnarray*}
}
}
\hbox to 6cm{
\parbox{6cm}{
\begin{eqnarray*}
\gamma (P\sb{13})&=&[1,0,1],\\
\gamma (P\sb{14}) &=& [1,0,0],\\
\gamma (P\sb{15}) &=&  [\alpha\sb 2, 0, 1],\\
\gamma (P\sb{16}) &=&  [\beta\sb 2, 0, 1] , \\
\gamma (P\sb{17}) &=&[0,1,1],\\
\gamma (P\sb{18})&=&[0,1,0], \\
\gamma (P\sb{19}) &=& [0, \alpha\sb 1, 1], \\
\gamma (P\sb{20}) &=& [0, \beta\sb 1, 1], \\
\gamma (P\sb{21})&=&[0,0,1]. \\
\phantom{\gamma (\phi(P\sb{21}))}&&\\
\phantom{\gamma (\phi(P\sb{21}))}&&\\
\phantom{\gamma (\phi(P\sb{21}))}&&
\end{eqnarray*}
}
}
}
}
}
\vskip 5pt
\caption{Parametric presentation of  $\gamma$}\label{table:parametricembedC}
\end{table}
The lines $L\sb {\nu}$ containing  the points $\gamma(l\sb{\nu})$ are defined by
\begin{multline*}
L \sb 1 =\{X=0\}, \quad 
L \sb 2 =\{Y=0\},\\
L \sb 3 =\{Y=t\sb 3 X \}, \quad 
L \sb 4 =\{Y=t\sb 4 X \}, \quad 
L \sb 5 =\{Y=t\sb 5 X \}.
\end{multline*}
%
%
%
%
\begin{claim}\label{claim:C1} $t\sb 5=1$.\end{claim}
Consider the quadratic word
$$
q\sb 1:=\{7, 8, 11, 12, 13, 14, 17, 18\}= \{2, 7, 12, 13, 18\}+ \{2, 8, 11, 14, 17\},
$$
%
%
which passes through the four points $P\sb{13}, P\sb{14}, P\sb{17}, P\sb{18}$,
and is disjoint from the linear word $l\sb 5$.
The conic curves containing the points
$\gamma (P\sb{13})=[1,0,1]$,
$\gamma (P\sb{14}) = [1,0,0]$,
$\gamma (P\sb{17}) =[0,1,1]$ and
$\gamma (P\sb{18})=[0,1,0]$
form a pencil 
$$
\sigma Z(X+Y+Z) +XY=0 \qquad(\sigma\in \P^1).
$$
The conic curve $Q \sb 1 \st\P^2$ containing $\gamma (q\sb 1)$ is a member of this pencil.
Since $Q \sb 1 $ is nonsingular,
the  value of the  parameter  $\sigma$ corresponding to $Q\sb 1$ is not $ 0$  nor $ \infty$.
Since $Q \sb 1 $  is tangent to the line $L \sb 5 =\{Y=t\sb 5 X \}$,
we have $\sigma  (1+t\sb 5)=0$.
Hence $t\sb 5=1$.
\par
\medskip
From the quadratic words that
\begin{itemize}
\item
contain exactly one of $\{P\sb{17}, P\sb{18}\}$,
\item
contain exactly one of $\{P\sb{13}, P\sb{14}\}$, and 
\item
are disjoint from $l\sb 5$,
\end{itemize}
we obtain the following relations:
\begin{claim}\label{claim:C2} 
$\alpha\sb 1=\alpha\sb 2(=:\alpha)$,  $\beta\sb 1=\beta\sb 2(=:\beta)$,
$\alpha+\beta=\alpha\beta$.\end{claim}
Consider, for example,
the quadratic word
$$
q\sb 2:=\{6, 8, 9, 11, 14, 16, 17, 19\}=\{2, 6, 9, 16, 19\}+ \{2, 8, 11, 14, 17\}.  
$$
Since the conic curve $Q \sb 2$ containing $\gamma (q\sb 2)$ passes through the  points 
$\gamma ( P\sb{14})=[1,0,0]$ and  $\gamma ( P\sb{17})=[0,1,1]$
and is tangent to $L \sb 5=\{ X=Y\}$,
it   is a member of the web
$$
\sigma\sb 1 (Y^2+ Z^2)+ \sigma\sb 2  X Y + \sigma \sb 3 (Y^2+YZ+ZX)=0
\qquad ([\sigma\sb 1, \sigma\sb 2 , \sigma \sb 3 ]\in \P^2)
$$
of conic curves.
Since $\gamma(P\sb{16})=[\beta\sb 2, 0,  1] \in Q \sb 2$,
we have $\beta\sb 2=\sigma\sb 1 /\sigma\sb 3 $.
Since $\gamma (P\sb{19})=[0, \alpha\sb 1,  1] \in Q \sb 2$
and $\alpha\sb 1\ne 1$, 
we have $\alpha\sb 1=\sigma\sb 1 /(\sigma\sb 1+\sigma\sb 3) $.
Therefore we obtain a relation $\alpha\sb 1 + \beta\sb 2 +\alpha\sb 1 \beta\sb 2 =0$.
\par
\medskip
From the quadratic words that
contain exactly  three of  $P\sb{17}, P\sb{18}, P\sb{13}, P\sb{14}$, 
we obtain the following relations:
\begin{claim}\label{claim:C3}

$\alpha\sb 1  +t\sb 3 =0, \quad \beta\sb 1 +t\sb 4=0$.

$1+\alpha\sb 2 t\sb 4 =0,\quad 1 + \beta\sb 2 t\sb 3=0$.

$1+\alpha\sb 2 + \alpha\sb 2 t\sb 3 =0, \quad 1 + \beta\sb 2 + \beta\sb 2 t\sb 4 =0$.

$\alpha\sb 1 +t\sb 4 + \alpha\sb 1 t\sb 4 =0,\quad \beta\sb 1 +t\sb 3 +\beta\sb 1 t\sb 3=0$.
\end{claim}
Consider, for example,
the quadratic word
$$
q\sb 3:=\{2, 4, 10, 12, 13, 16, 17, 18\}=\{4, 7, 10, 16, 17\}+ \{2, 7, 12, 13, 18\}.  
$$
Since the conic curve  $Q \sb 3$ containing $\gamma (q\sb 3)$  passes through the points 
$\gamma ( P\sb{13})=[1,0,1]$, $\gamma ( P\sb{17})=[0,1,1]$ and 
 $\gamma ( P\sb{18})=[0,1,0]$,
it   is a member of the web
$$
\sigma\sb 1 (X^2+Z^2+YZ)+ \sigma\sb 2 (X^2+ZX) + \sigma \sb 3 XY=0
\qquad ([\sigma\sb 1, \sigma\sb 2 , \sigma \sb 3 ]\in \P^2)
$$
of conic curves.
Since $\gamma (P\sb{16})=[\beta\sb 2, 0,  1]$ is contained in $ Q \sb 3$,
we obtain ${\beta\sb 2} \sp 2 (\sigma\sb 1 + \sigma \sb 2)+\beta\sb 2 \sigma\sb 2 +\sigma\sb 1=0$.
Since  $Q \sb 3$ is tangent to the line $L \sb 4 =\{Y=t\sb 4 X\}$,
we have $t\sb 4 \sigma\sb 1 +\sigma\sb 2=0$.
Combining these two relations and $\beta\sb 2 \ne 1$, 
we obtain a relation $1+\beta\sb 2+\beta\sb 2 t\sb 4=0$.
\par
\medskip
Combining Claims~\ref{claim:C1}-\ref{claim:C3},
we obtain the following two possibilities for the parameters;
\begin{eqnarray*}\label{eq:twopos}
&& \alpha\sb 1=\alpha\sb 2=\omega, \;\; \beta\sb 1=\beta\sb 2=\bar\omega, \;\; 
t\sb 3=\omega, \;\; t\sb 4 =\bar\omega, \;\; t\sb 5=1,
\quad\textrm{or}\\ 
&& \alpha\sb 1=\alpha\sb 2=\bar\omega, \;\; \beta\sb 1=\beta\sb 2=\omega, \;\; t\sb 3=\bar\omega, \;\; t\sb 4
=\omega, \;\;  t\sb5=1.
\end{eqnarray*}
If the latter holds,
then we replace $\gamma$ by $\gamma \circ T$
so that we  assume that the former  always holds.
\par
\medskip
Next we put
$$
P\sb 1=[1,1,\lambda  ],
$$
where $\lambda  =s\sb{5,1}$ is a parameter.
Using quadratic words that
\begin{itemize}
\item
contain exactly four points among $l\sb 1 \cup l\sb 2$, and 
\item
are not disjoint from $l\sb 5$,
\end{itemize}
we obtain the following:
\begin{claim}\label{claim:C4}
\begin{eqnarray*}
&& 
s\sb{5, 1}=\lambda  , \quad 
s\sb{5, 2}=\lambda  +1, \quad 
s\sb{5, 3}=\lambda  +\oo, \quad  
s\sb{5, 4}=\lambda  +\bo, \\
&& 
s\sb{4, 1}=\oo \lambda   +\oo, \quad 
s\sb{4, 2}=\oo\lambda  , \quad 
s\sb{4, 3}=\oo\lambda   +1, \quad  
s\sb{4, 4}=\oo\lambda   +\bo, \\
&& 
s\sb{3, 1}=\bo\lambda   +\bo, \quad 
s\sb{3, 2}=\bo\lambda  , \quad 
s\sb{3, 3}=\bo\lambda   +\oo, \quad  
s\sb{3, 4}=\bo\lambda  +1. \\
\end{eqnarray*}
\end{claim}
Consider, for example,
the quadratic word
$$
q\sb 4:=\{1, 2, 11, 12, 14, 16, 17, 20\}=\{1, 8, 12, 16, 20\}+ \{2, 8, 11, 14, 17\},  
$$
which is disjoint from $l\sb 4$.
Because  there exists a conic curve $Q \sb 4$ that contains    $\gamma (q\sb 4)$
and is tangent to the line $L \sb 4=\{ Y=\bo X\}$,
the following matrix $\wtil{M}$  is of rank $< 6$.
$$
\wtil{M}:=\left [\begin {array}{cccccc} 
1&1&{\lambda  }^{2}&1&\lambda  &\lambda  \\
\noalign{\medskip}1&1&{s_{{5,2}}}^{2}&1&s_{{5,2}}&s_{{5,2}}\\
\noalign{\medskip}1&\bo&{s_{{3,3}}}^{2}&\oo&\oo\,s_{{3,3}}&s_{{3,3}}\\\noalign{\medskip}1&\bo&{s_{{3,4}}}^{2}&\oo&\oo\,s_{{3,4}}&s_{{3,4}}\\
\noalign{\medskip}1&0&0&0&0&0\\
\noalign{\medskip}\oo&0&1&0&0&\bo\\
\noalign{\medskip}0&1&1&0&1&0\\
\noalign{\medskip}0&\oo&1&0&\bo&0\\
\noalign{\medskip}0&0&0&0&\bo&1\end {array}\right ].
$$
Indeed, if 
 the equation 
$$
a\sb 1 X^2 + a\sb 2 Y^2 + a\sb 3 Z^2 + a\sb 4 XY + a\sb 5 YZ + a\sb 6 ZX=0
$$
defines 
a conic  curve containing  $\gamma (q\sb 4)$
and tangent to $L \sb 4$,
then
$\va=\spT [a\sb 1, a\sb 2, \dots, a\sb 6]$ is a non-zero  solution of
$\wtil{M}\vx=\mathord{\bf 0}$.
(The condition $\bo a\sb 5 +a\sb 6=0$ is equivalent to the condition that
the conic curve is tangent to  $L \sb 4$.)
Let $\wtil{M} [i\sb 1, \dots, i\sb 6]$ denote the submatrix of $\wtil{M}$
consisting of  $i\sb j$-th rows of $\wtil{M}$.
Because
$$
\det \wtil{M}[1, 2, 5, 7, 8, 9] =
\left (s_{{5,2}}+\lambda  \right )\left (s_{{5,2}}+\lambda  +1\right )
$$
and $s_{{5,2}}\ne s_{{5,1}}=\lambda  $,
we obtain $s_{{5,2}}=\lambda  +1$.
Because
$$
\det \wtil{M}[1, 3, 5, 7, 8, 9] =
\left (s_{{3,3}}+\bo \lambda  +1\right )
\left (s_{{3,3}}+\bo\lambda  +\oo\right ),
$$
we
obtain
$$
s\sb{3,3}= \bo \lambda  +1\quad\textrm{or}\quad \bo\lambda  +\oo.
$$
Continuing the same calculations,
we get the relations in Claim~\ref{claim:C4}.
\par
\medskip
Thus we have proved Claim~\ref{claim:tripleC}.
\par
\medskip
Conversely,
suppose that $\lambda\in k\sm \{0,1,\oo, \bo\}$ is given.
Then $\gamma\sb{\lambda} (\Ps)$ is equal to $\ZZ{d\GC[\lambda]}$,
where $\GC[\lambda]$ is given in Theorem~\ref{thm:C}.
Moreover,
for every linear word $l$ of $\CC$,
there exists a line containing $\gamma\sb{\lambda} (l)$,
and 
for every quadratic word $q$ of $\CC$,
there exists a unique conic curve  containing $\gamma\sb{\lambda} (q)$.
The defining equations of the $120$ conic curves are omitted.
These conic curves are nonsingular because  $\lambda\notin \{0,1,\oo, \bo\}$.
Hence
$\gamma\sb\lambda$ is in $\gs\sb B$ by Proposition~\ref{prop:ings}.
\end{proof}
\begin{remark}\label{rem:limitDKC}
When $\lambda\in \{0,1, \omega, \bar\omega\}$,
the set $\gamma\sb\lambda (\Ps)$ coincides with $\Pt (\Ff)$,
and the point $[\GC[\lambda]]\in \moduli$ is equal to the Dolgachev-Kondo point.
\end{remark}
For each $\sigma\in \Aut (\CC)$,
we calculate the unique triple
$$
(g\sb{\sigma}, t\sb{\sigma}, \lambda\sp{\sigma})\;\in\; \PGL (3, k(\lambda)) \times \widetilde{T}\times k(\lambda)
$$
such that 
$g\sb{\sigma} \circ (\gamma\sb{\lambda}\circ \sigma)\circ t\sb{\sigma}=\gamma\sb{\lambda\sp{\sigma}}$
holds.
The map $\sigma\mapsto t\sb{\sigma}$ is a homomorphism from $\Aut (\CC)$
to $\widetilde T$. We  put
$N\sb C:=\Ker (\Aut (\CC)\to \widetilde{T} )$.
\begin{corollary}\label{cor:connectedcompsC}
The space $\PGL(3, k)\backslash \gs\sb C$ has exactly two connected
components
$$
(\PGL(3, k)\backslash \gs\sb C)\sp +
:=\set{[\gamma\sb\alpha]}{\alpha\in k\sm \{0,1,\omega,\bar\omega\}}
$$
and
$((\PGL(3, k)\backslash \gs\sb C)\sp +)\cdot T$, each of which is isomorphic to $\A\sp 1\sm
\{0,1,\omega,\bar\omega\}$. 
The group $N\sb C$ acts on $(\PGL(3, k)\backslash \gs\sb C)\sp +$, and
the moduli curve $\moduli\sb C$ is equal to the quotient space $(\PGL(3, k)\backslash\gs\sb C)\sp +/N\sb C$.
\end{corollary}
Consider the natural projection
$$
\map{p\sb C}{ \A\sp 1\sm \{0,1,\omega, \bar\omega\} \cong (\PGL (3, k)
\backslash \GGG\sb{C})\sp + }{%
\moduli\sb{C} = (\PGL (3, k) \backslash \GGG\sb{C})\sp + / N\sb C}.
$$
For $\alpha\in k\sm\{0,1,\omega,\bar\omega\}$, 
let $P[\alpha]$ denote the point of $\A\sp 1\sm \{0,1,\omega, \bar\omega\}$
given by $\lambda=\alpha$.
Then $p\sb C (P[\alpha])\in\moduli\sb C$ corresponds to the isomorphism
class  of the polarized supersingular $K3$ surface 
$( X\sb{\GC [\alpha]},\pol\sb{\GC [\alpha]})$.
\begin{proposition}\label{prop:groupC}
We have
$$
p\sb C\inv (p\sb C (P[\alpha]))=\set{P[u\alpha +v]}{u\in \F\sb 4\sptimes, v\in \F\sb 4}.
$$
The group $\Aut( X\sb{\GC [\alpha]},
\pol\sb{\GC [\alpha]})$ is equal to the subgroup~\eqref{eq:groupC} of $\PGL (3, k)$.
\end{proposition}
\begin{corollary}
We have
$\moduli\sb{C}=\Spec k[J\sb C, 1/J\sb C]$,
where 
$J\sb C:=(\lambda^4+\lambda)^3$.
The morphism $p\sb C$ is an \'etale  Galois covering with Galois group
$\Gamma\sb C$.
\end{corollary}
\section{Cremona transformations by quintic curves}\label{sec:cremona}
\subsection{Preliminaries}
Let $\Sigma\sb 1$ and $\Sigma\sb 2$ be disjoint  sets  of reduced points of $\Pt$
with $|\Sigma\sb 1|=n\sb 1$ and $|\Sigma\sb 2|=n\sb 2$,
and 
let $\ideal\sb{\Sigma\sb 1}\st\OOO\sb{\Pt}$ and $\ideal\sb{\Sigma\sb 2}\st\OOO\sb{\Pt}$ be the ideal sheaves defining 
$\Sigma\sb 1$ and
$\Sigma\sb 2$.
We define $\wtil{\Sigma}$ to be the $0$-dimensional subscheme of $\Pt$ defined by the ideal sheaf
$$
\ideal \sb{\wtil{\Sigma}} :=\ideal \sb{\Sigma\sb 1}\sp{\phantom{2}} \ideal \sb{\Sigma \sb 2}\sp 2.
$$
The length of $\OOO\sb{\wtil{\Sigma}}$ is $n\sb 1+3 n\sb 2$.
Let $d$ be a positive integer.
The linear system $|\ideal \sb{\wtil{\Sigma}} (d)|$
consists of plane curves of degree $d$ that pass through the points of $\Sigma\sb 1 \cup \Sigma\sb 2$
and are singular at each point of $\Sigma\sb 2$.
\begin{proposition}\label{prop:fcf}
Suppose that the linear system
$|\ideal \sb{\wtil{\Sigma}} (d)|$ is of dimension $\ge 1$ and has no fixed components.
If
\begin{equation}\label{eq:dimassump}
\dim |\ideal \sb{\wtil{\Sigma}} (d)| > \frac{(d+2)(d+1)}{2} -(n\sb 1+3 n\sb 2) -1,
\end{equation}
then there exists a projective plane curve of degree $d-3$ that passes through 
all the points of $\Sigma\sb 2$.
\end{proposition}
\begin{corollary}\label{cor:fcf}
Suppose that the linear system
$|\ideal \sb{\wtil{\Sigma}} (d)|$ is of dimension $\ge 1$ and  has no fixed components.
If $d\le 3$ and $n\sb 2 >0$, then the dimension of the linear system
$|\ideal \sb{\wtil{\Sigma}} (d)|$ is equal to $(d+2)(d+1)/2 -(n\sb 1+3 n\sb 2) -1$.
\end{corollary}
\begin{proof}
We follow the argument in~\cite[pp.712-714]{GH}.
From the exact sequence
$$
0 \;\to\; \ideal \sb{\wtil{\Sigma}} (d) \;\to\; \OOO \sb{\Pt} (d) \;\to\; \OOO\sb{\wtil{\Sigma}} (d)\;\to\; 0,
$$
we obtain
\begin{equation}\label{eq1}
h\sp 0 (\Pt, \ideal \sb{\wtil{\Sigma}} (d))= (d+2)(d+1)/2 -(n\sb 1+3n\sb 2) + h\sp 1 (\Pt, \ideal \sb{\wtil{\Sigma}} (d)).
\end{equation}
Let $\blu: S\to\Pt$ denote the blowing up of $\Pt$ at the points of $\Sigma\sb 1 \cup \Sigma\sb 2$.
We put
$$
\Delta\sb 1:= \blu\inv (\Sigma\sb 1), \quad \Delta\sb 2 :=\blu\inv (\Sigma\sb 2),
$$
both of which are considered to be reduced divisors.
Let $\Line\st S$ be the pull-back of a general line on $\Pt$.
We  put
$$
L:=\blu\sp * \OOO\sb{\Pt} (d) \otimes \OOO\sb S ( -\Delta\sb 1 -2 \Delta\sb 2)=\OOO\sb{S}(dH-\Delta\sb 1-2\Delta\sb 2).
$$
Because $K\sb S =-3\Line +\Delta\sb 1 +\Delta\sb 2$,
we have
$$
L^2=d^2-n\sb 1 - 4 n\sb 2, \qquad L\intsct K\sb S = -3d +n\sb 1 + 2 n\sb 2.
$$
The complete linear system $|\Line|=|\blu\sp * \OOO\sb{\Pt} (1) |$ on $S$ is fixed component free.
Since $\Line\intsct (K\sb S-L)=-d-3<0$, we have
$$
h\sp 2 (S, L)= h\sp 0 (S, K\sb S-L)=0.
$$
By the Riemann-Roch theorem, we obtain
\begin{equation}\label{eq2}
h\sp 0 (S, L)= (d+2)(d+1)/2 -(n\sb 1+3n\sb 2) + h\sp 1 (S, L).
\end{equation}
There exists a canonical isomorphism 
\begin{equation}\label{eq3}
|\ideal \sb{\wtil{\Sigma}} (d)|\;\cong\; |L|
\end{equation}
that maps a member $C$ of $|\ideal \sb{\wtil{\Sigma}} (d)|$ to the member
$\blu\sp * C-\Delta\sb 1 -2\Delta\sb 2$ of $|L|$.
From~\eqref{eq1}, \eqref{eq2} and~\eqref{eq3}, we obtain
\begin{equation}\label{eq4}
h\sp 1 (\Pt, \ideal \sb{\wtil{\Sigma}} (d))=h\sp 1 (S, L).
\end{equation}
Using the assumption~\eqref{eq:dimassump} and the equalities~\eqref{eq1} and~\eqref{eq4}, we obtain
\begin{equation}\label{eq5}
h\sp 1 (S, L)>0.
\end{equation}
Since $|\ideal \sb{\wtil{\Sigma}} (d)|$ is of dimension $\ge 1$ and has no fixed components,
we obtain by the isomorphism~\eqref{eq3} global sections $s$ and $s\sprime$ of $L$
such that the subscheme $R=\{s=s\sprime=0\}$ of $S$ 
is of dimension $0$.
Let $\ideal \sb{R}\st\OOO\sb{S}$ be the ideal sheaf defining $R$.
From the Koszul complex
$$
0
\;\maprightsp{}\; 
\OOO\sb S (K\sb S -L)
\; \maprightsp{(s, s\sprime)}\;
\OOO\sb S (K\sb S) \oplus \OOO\sb S (K\sb S)
\; \maprightsp{{}\sp T \hskip -1pt (-s\sprime, s)}\;
\ideal \sb R (K\sb S+L)
\; \maprightsp{}\;
0
$$
and $h\sp 0 (S, \OOO\sb S (K\sb S))=h\sp 1 (S, \OOO\sb S (K\sb S))=0$,
we obtain 
$$
h\sp 1 (S, L)= h\sp 1 (S, \OOO\sb S (K\sb S-L)) =h\sp 0 (\ideal \sb R (K\sb S +L)).
$$
From~\eqref{eq5}, we see that the linear system $|\ideal \sb R (K\sb S +L)|$ is non-empty.
Since $K\sb S+L= \blu\sp * \OOO\sb{\Pt} (d-3) \otimes\OOO\sb S (-\Delta\sb 2)$,
a member of $|\ideal \sb R (K\sb S +L)|$ is mapped by $\blu$ to a projective plane curve of degree $d-3$ that passes through
the points of $\Sigma\sb 2$.
\end{proof}
\begin{definition}
Let $F$ be an effective divisor of $\Pt$.
We put
$$
\Sigma\sb 1 \sprime:= (\Sigma\sb 1\sm (\Sigma\sb 1 \cap F)) \cup (\Sigma\sb 2 \cap F\sp 0),
\qquad
\Sigma\sb 2 \sprime:= \Sigma\sb 2\sm (\Sigma\sb 2 \cap F),
$$
where $F\sp 0$ is the locus of all $p\in \Supp (F)$ at which $F$ is reduced and nonsingular.
We then define  $\wtil{\Sigma}\sm F$ to be the  $0$-dimensional subscheme of $\Pt$
defined by the ideal sheaf
$$
\ideal \sb{\wtil{\Sigma} \sm F}:=\ideal \sp{\phantom{2}}\sb{\Sigma\sb 1\sprime}\ideal \sp 2 \sb{\Sigma\sb 2\sprime}.
$$
\end{definition}
If $F$ is a fixed component of $|\ideal \sb{\wtil{\Sigma}} (d)|$,
then $C\mapsto C-F$ gives an isomorphism
$$
|\ideal \sb{\wtil{\Sigma}} (d)|\;\cong\; |\ideal \sb{\wtil{\Sigma}\sm F} (d-\deg F)|
$$
of linear systems.
By the definition,
we have
\begin{equation}\label{eq:F1F2}
\wtil{\Sigma} \sm (F\sb 1+F\sb 2)=(\wtil{\Sigma} \sm F\sb 1)\sm F\sb 2
\end{equation}
for any effective (not necessarily distinct) divisors $F\sb 1$ and $F\sb 2$ of $\Pt$.
\subsection{A homaloidal system of quintic curves}
Let $\Sigma=\{p\sb 1, \dots, p\sb 6\}$
be a set of distinct six points of $\Pt$ satisfying the following:
\begin{itemize}
\item[($\Sigma$1)] no three points of $\Sigma$ are collinear, and 
\item[($\Sigma$2)] there are no conic curves containing $\Sigma$.
\end{itemize}
These conditions are equivalent to the following:
\begin{itemize}
\item[($\Sigma$3)] for each $p\sb i\in \Sigma$, there exists a nonsingular conic curve $\EED\sprime\sb i\st \Pt$
that contains $\Sigma\sm \{p\sb i\}$ and does not contain $p\sb i$.
\end{itemize}
\begin{proposition}\label{prop:Lfcf}
The  linear system $|\ideal \sb{\Sigma}\sp 2 (5)|$ of quintic curves that
pass through the points of $\Sigma$ and are singular at each point of $\Sigma$
is of dimension $2$, and has no fixed components.
\end{proposition}
\begin{proof}
Because each point of $\Sigma$ imposes three linear conditions on $|\OOO \sb{\Pt} (5)|$,
we have $\dim |\ideal \sb{\wtil{\Sigma}} (5)|\ge 2$.
\par
Suppose that $|\ideal \sb{\Sigma}\sp 2 (5)|$ has a fixed component.
Let $F$ be the fixed component, and let 
$$
F=F\sb 1 + \cdots + F\sb N
$$
be the decomposition into the reduced irreducible components of $F$,
where non-reduced components are expressed by repetition of summation.
We have
\begin{equation*}\label{eqF:deg}
\deg F=\deg F\sb 1 +\cdots + \deg F\sb N\;\;>\;\; 0.
\end{equation*}
As in the previous subsection, we denote by $\wtil{\Sigma}$ the $0$-dimensional subscheme of $\Pt$
defined by the ideal sheaf $\ideal \sb{\Sigma}^2$.
We will consider the linear system $|\ideal \sb{\wtil{\Sigma} \sm F} (5-\deg F)|$,
which has no fixed components and is of dimension equal to $\dim |\ideal \sb{\wtil{\Sigma}} (5)|\ge 2$.
For $\nu=0, \dots, N$, we define reduced $0$-dimensional
subschemes $\Sigma\sp{(\nu)}\sb 1$ and  $\Sigma\sp{(\nu)}\sb 2$ of $\Pt$  by
$$
\ideal \sb{\wtil{\Sigma} \sm (F\sb 1 + \cdots + F\sb{\nu})}=
\ideal \sb{\Sigma\sp{(\nu)}\sb 1}\sp{\phantom{2}} \ideal \sb{\Sigma\sp{(\nu)}\sb 2}\sp 2.
$$
Then
$$
|\Sigma\sp{(\nu+1)}\sb 1|= |\Sigma\sp{(\nu)}\sb 1|-i+j-k, \quad
|\Sigma\sp{(\nu+1)}\sb 2|= |\Sigma\sp{(\nu)}\sb 2|-j,
$$
where
$$
i:=|\Sigma\sp{(\nu)}\sb 1 \cap F\sb{\nu+1}|,
\quad
j:=|\Sigma\sp{(\nu)}\sb 2 \cap F\sb{\nu+1}|,
\quad
k:=|\Sigma\sp{(\nu)}\sb 2 \cap \Sing F\sb{\nu+1}|.
$$
The integers $i$, $j$ and $k$ are subject to the following conditions:
\begin{itemize}
\item $i+j\le 2$ and $k=0$ if $\deg F\sb{\nu+1}=1$, because of ($\Sigma$1), 
\item $i+j\le 5$ and $k=0$ if $\deg F\sb{\nu+1}=2$, because of ($\Sigma$2),
\item $k\le 1$ if $\deg F\sb{\nu+1}=3$, because an irreducible cubic curve has at most one singular point, and
\item $k\le 4$ if $\deg F\sb{\nu+1}=4$, because if $k\ge 5$, there would exist a conic curve $C$ with $C\intsct F\sb{\nu+1}\ge
10$.
\end{itemize}
Since $|\ideal\sb{\Sigma}\sp 2 (5)|\cong |\ideal \sb{\Sigma\sb 1\sp{(N)}} \ideal\sb{\Sigma\sb 2 \sp{(N)}}^2 (5-\deg F)|$
is of dimension $\ge 2$ and fixed component free, we have
\begin{eqnarray*}
\deg F =4 &\Longrightarrow & |\Sigma\sb 1 \sp{(N)}|\le 1\quand |\Sigma\sb 2 \sp{(N)}|=0,\\
\deg F =3 &\Longrightarrow & |\Sigma\sb 1 \sp{(N)}|\le 4\quand |\Sigma\sb 2 \sp{(N)}|=0.
\end{eqnarray*}
We put
$$
\delta :=(6-\deg F)(7-\deg F)/2-(|\Sigma\sb 1 \sp{(N)}|+3 |\Sigma\sb 2 \sp{(N)}|)-1.
$$
From Corollary~\ref{cor:fcf}, we also have
$$
\deg F\ge 2 \;\;\textrm{ and }\;\; 2>\delta \;\;\Longrightarrow\;\;
|\Sigma\sb 2 \sp{(N)}|=0.
$$
Using these considerations,
we see that the triple $(|\Sigma\sb 1 \sp{(N)}|, |\Sigma\sb 2 \sp{(N)}|, \deg F)$ is one of the following:
$$
(0, 6, 1), \;\;  (1, 5, 1), \;\; (2, 4, 1).
$$
For these triples, however, we have $ |\ideal \sb{\wtil{\Sigma}\sm F} (5-\deg F)|=\emptyset$,
because otherwise, there would exist an irreducible quartic curve $C\sb 4$ and a conic curve $C\sb 2$
such that $C\sb 4\intsct C\sb 2>8$.
Thus we have proved that $|\ideal \sb{\wtil{\Sigma}} (5)|$ is fixed component free.
\par
If $\dim |\ideal \sb{\wtil{\Sigma}} (5)|>2 $, then, by Proposition~\ref{prop:fcf},
there would exists a conic curve that contains $\Sigma$,
which contradicts ($\Sigma$2).
\end{proof}
\begin{remark}\label{rem:multQ}
Recall from ($\Sigma$3) that $\EED\sprime\sb i\st\Pt$ is the conic curve such that  
$\EED\sprime\sb i\cap\Sigma=\Sigma\sm \{p\sb
i\}$. Let $Q$ be a general member of $|\ideal\sp 2\sb{\Sigma} (5)|$.
Since $\EED\sprime\sb i Q=10$ for each $i$,
the multiplicity of $Q$ at each point of $\Sigma$ is $2$.
\end{remark}
Let $\blu:S \to \Pt$ be the blowing up of $\Pt$ at the points of $\Sigma$,
and 
let $\ED\sb i$ be the exceptional (reduced) divisor $\blu\inv (p\sb i)$.
We put
$$
L:=\blu\sp * \OOO\sb{\Pt} (5) \otimes \OOO\sb S (- 2\ED\sb 1-\cdots-2 \ED\sb 6).
$$
Let $\EED\sb i$ be the strict transform of $\EED\sprime\sb i$ by $\blu$.
We have $L^2=1$, $\EED\sb i\intsct L=0$ and $\EED\sb i^2=-1$.
\begin{proposition}
The complete linear system $|L|$ on $S$ has no base points,
and the morphism $\Phi\sb{|L|} : S\to\Pt$ defined by $|L|$ is the contraction of the curves $\EED\sb 1, \dots, \EED\sb 6$.
Let $p\sprime\sb i$ be the image of $\EED\sb i$ by $\Phi\sb{|L|}$.
Then $\Sigma\sprime=\{p\sprime\sb 1, \dots, p\sprime\sb 6\}$ satisfies the condition {\rm ($\Sigma$3)}.
\end{proposition}
\begin{proof}
By Proposition~\ref{prop:Lfcf}, the complete linear system $|L|$ on $S$ is  of dimension $2$ and has no fixed components.
Suppose that $|L|$ has a base point $p\in S$.
Let $\tilde{\blu}:\wtil{S}\to S$ be the blowing up of $S$ at $p$,
and let $\ED\sprime$ be the exceptional divisor $\tildeblu\inv (p)$ of $\tildeblu$.
Since $L^2=1$, 
the complete linear system $|\wtil{L}|$ of the line bundle
$\wtil{L}:=\tilde{\blu}\sp * L\otimes\OOO\sb{\wtil{S}} (- \ED\sprime)$ is of dimension $2$ and 
has no fixed components.
 We have  $K\sb{\wtil{S}}\intsct \wtil{L}=-2$, and hence 
$h\sp 2 (\wtil{S}, \wtil{L})= 
h\sp 0 (\wtil{S}, \OOO\sb{\wtil{S}} (K\sb{\wtil{S}}-\wtil{L}))=0$ 
follows. By Riemann-Roch theorem,
we have 
$h\sp 1 (\wtil{S}, \wtil{L})= h\sp 1 (\wtil{S}, K\sb{\wtil{S}}-\wtil{L})=1$.
Using the argument of Koszul complex as in the proof of Proposition~\ref{prop:fcf},
we see that $h\sp 0 (\wtil{S}, \OOO\sb{\wtil{S}} (K\sb{\wtil{S}}+\wtil{L}))>0$.
Hence there exists a conic curve in $\Pt$ that passes through the points of $\Sigma$,
which contradicts ($\Sigma$2).
Thus $|L|$ has no base points.
\par
Since $L^2=1$, the morphism $\Phi\sb{|L|}$ is of degree $1$.
Because $\EED\sb i\intsct L=0$,
the curves $\EED\sb i$ are contracted by $\Phi\sb{|L|}$.
Let $C$ be a reduced irreducible curve on $S$ that is contracted by $\Phi\sb{|L|}$.
Because $\ED\sb i\intsct L=2$,
we have $C\ne \ED\sb i$ and hence $C\sprime:=\blu (C)\st\Pt$
is a reduced  irreducible curve.
We will show that $C\sprime$ is equal to one of the conic curves $\EED\sprime\sb i$.
Let $d$ be the degree of $C\sprime$.
We have 
$$
\blu\sp * (C\sprime)=C+m\sb 1 \ED\sb 1 + \cdots + m\sb 6 \ED\sb 6,
$$
where $m\sb j$ is the multiplicity of $C\sprime$ at $p\sb j$.
The condition $C\intsct L=0$ implies 
$$
5d=2(m\sb 1+\cdots+m\sb 6).
$$
If $C\sprime$ is not equal to $\EED\sprime\sb i$ for any $i$,
then 
$$
C\sprime\intsct  \EED\sprime\sb i= 2d\ge (m\sb 1+\cdots+m\sb 6)-m\sb i =5d/2 -m\sb i
$$
holds for each $i$.
Hence $2m\sb i\ge d$ for $i=1, \dots, 6$.
Therefore $5d=2\sum m\sb j\ge 6 d$, which is absurd.
Thus we have proved that $\Phi\sb{|L|}$ is the contraction of the $(-1)$-curves $\EED\sb 1, \dots, \EED\sb 6$.
\par
Since $\ED\sb i\intsct L=2$,
the image of $\ED\sb i$ by $\Phi\sb{|L|}$ is a nonsingular conic curve.
Because $\ED\sb i\intsct \EED\sb j= 0$ if and only if $i=j$,
the conic curve $\Phi\sb{|L|}(\ED\sb i)$ satisfies
$$
\Phi\sb{|L|}(\ED\sb i)\cap \Sigma\sprime=\Sigma\sprime\sm \{ p\sprime\sb i\}.
$$
Hence $\Sigma\sprime$ satisfies ($\Sigma$3).
\end{proof}
\begin{corollary}
The rational map
$$
\cremona\sb{\Sigma} : \Pt \cdots\to \Pt
$$
defined by the  linear system $|\ideal \sb{\Sigma}\sp 2 (5)|$ is birational,
and the inverse map is given by $\cremona\sb{\Sigma\sprime}$.
\end{corollary}
We will write 
$$
\map{\bld}{S}{\Pt}
$$
instead of $\Phi\sb{|L|}$.
Let $\Line$ and $\Line\sprime$ be the pull-backs of a general line of $\Pt$ by $\blu$ and $\bld$,
respectively.
We put
$$
\ED\sb i\sprime:=\bld (\ED\sb i),
$$
which is a nonsingular conic curve containing $\Sigma\sprime\sm \{p\sb i\sprime\}$
and not passing through $p\sb i\sprime$.
We also put
$$
U:=S\sm \left(\bigcup\sb{i=1}\sp{6} \EED\sb i\cup \bigcup\sb{i=1}\sp{6} \ED\sb i\right).
$$
The morphisms $\blu$ and $\bld$ induce isomorphisms
\begin{equation}\label{eq:isomU}
\Pt\sm \cup \EED\sprime\sb i \;\;\cong\;\; U \;\;\cong\;\; \Pt\sm \cup \ED\sprime\sb i.
\end{equation}
The Picard group $\Pic S$  of $S$ is a free $\Z$-module of rank $7$,
and is generated by the linear equivalence classes 
$[\Line], [\ED\sb 1], \dots, [\ED\sb 6]$, or by the linear equivalence classes 
$[\Line\sprime], [\EED\sb 1], \dots, [\EED\sb 6]$.
They are related by
$$
[\Line\sprime]=5[\Line]-2\sum\sb{j=1}\sp {6} [\ED\sb j], 
\quad
[\EED\sb i]=2[\Line]-\sum\sb{j=1}\sp 6 [\ED\sb j]+ [\ED\sb i]\quad (i=1, \dots, 6).
$$
In particular, we have
\begin{equation}\label{eq:inPicS}
3[\Line]-\sum\sb{j=1}\sp{6} [\ED\sb j]= 3 [\Line\sprime] -\sum\sb{j=1}\sp 6 [\EED\sb j]
\end{equation}
in $\Pic S$.
\subsection{Cremona transformations of  supersingular $K3$ surfaces}
Let $G$ be a homogeneous polynomial in $\UUU$,
and let $\Sigma=\{ p\sb 1, \dots, p\sb 6\}$
be a subset of $\ZZ{dG}$ with $|\Sigma|=6$.
We assume that $\Sigma$ satisfies the condition ($\Sigma$1) and 
\begin{itemize}
\item[($\Sigma$2)${}\sprime$]
for each $p\sb i\in \Sigma$,
the nonsingular conic curve $\EED\sprime\sb i$ containing $\Sigma\sm \{ p\sb i\}$ 
satisfies $\EED\sprime\sb i \cap \ZZ{dG}=\Sigma\sm \{p\sb i\}$.
\end{itemize}
Then the subset
$$
Z\sprime := \cremona\sb{\Sigma} (\ZZ{dG}\sm \Sigma) \;\cup\; \Sigma\sprime
$$
of $\Pt$ is well-defined and consists of $21$ points.
\begin{proposition}\label{prop:Gsprime}
There exists $G\sprime\in \UUU$ such that $Z\sprime=\ZZ{dG\sprime}$.
\end{proposition}
For the proof of Proposition~\ref{prop:Gsprime}, we first show the following:
\begin{lemma}\label{lem:G1}
There exists $G\sb 1\in \UUU$ that satisfies  $\ZZ{dG}=\ZZ{dG\sb 1}$ and 
$G\sb 1 (p\sb i)=0$ for each $p\sb i\in \Sigma$.
\end{lemma}
\begin{proof}
By $(\Sigma1)$ and $(\Sigma2)\sprime$,
the points of $\Sigma$ impose independent linear conditions on 
the linear system $|\OOO\sb{\Pt} (3)|$.
(See~\cite[p.~715]{GH}.)
Hence
there exists a homogeneous cubic polynomial $H$ such that
$(G+H^2)(p\sb i)=0$ holds for each $p\sb i\in \Sigma$.
Then $G\sb 1:=G+H^2$ satisfies $\ZZ{dG}=\ZZ{dG\sb 1}$.
\end{proof}
We replace $G$ by $G\sb 1$ in Lemma~\ref{lem:G1}.
Then the sextic curve $D$ defined  by $G=0$ is reduced and has an ordinary node
at each point of $\Sigma$.
Hence 
$$
\wtil{D}:=\blu\sp * D -2\sum\sb{j=1}\sp 6  \ED\sb j
$$
is a reduced  effective divisor of $S$,
and it does not contain any of $\ED\sb j$.
\begin{proof}[Proof of Proposition~\ref{prop:Gsprime}]
Let $D\sprime$ be the image of $\wtil{D}$ by $\beta\sprime$.
Since $\wtil{D}\Line\sprime=6$,
$D\sprime$ is a reduced curve of degree $6$.
Let
$G\sprime=0$  be
the defining equation of $D\sprime$.
We will show that $Z\sprime=\ZZ{dG\sprime}$.
It is enough to show that 
$\ZZ{dG\sprime}$ is of dimension $0$ and that 
$Z\sprime\subseteq\ZZ{dG\sprime}$.
\par
Since $\wtil{D} \EED\sb j=2$ for each $\EED\sb j$,
we have
$$
\wtil{D}:=\bld\sp * D\sprime -2\sum\sb{j=1}\sp 6  \EED\sb j.
$$
Because $\wtil{D}$ is effective,
we have $\Sing D\sprime\supset \Sigma\sprime$.
Hence $\Sigma\sprime\subseteq\ZZ{dG\sprime}$.
We put 
$$
\HL:= \OOO\sb{S} (3\Line-\sum\sb{j=1}\sp{6} \ED\sb j)=\OOO\sb{S} ( 3 \Line\sprime -\sum\sb{j=1}\sp 6 \EED\sb j).
$$
Let $\wtil{G}$ be a global section of 
$$
\HL\sp{\otimes 2}=L(H)
=\OOO\sb{S} ( 6 \Line -2 \sum \ED\sb j)
=\OOO\sb{S} ( 6 \Line\sprime - 2\sum \EED\sb j)
$$
such that 
$\wtil{G}=0$ defines $\wtil{D}$.
Let $m$ and $n$ be global sections of $\OOO\sb{S} (\sum \ED\sb j)$ and $\OOO\sb{S} (\sum \EED\sb j)$
such that $\sum \ED\sb j=\{m=0\}$ and $\sum \EED\sb j =\{ n=0\}$.
We can choose them so that
\begin{equation}\label{eq:sq}
\blu\sp * G= \wtil{G}\cdot m^2
\quand
\bld\sp * G\sprime = \wtil{G}\cdot n^2
\end{equation}
hold.
We define isomorphisms
\begin{equation}\label{eq:isomHL}
\blu\sp * \OOO\sb{\Pt} (3) |U \;\;\cong\;\; \HL|U\;\;\cong\;\; \bld\sp * \OOO\sb{\Pt} (3) |U
\end{equation}
of line bundles on $U$
by multiplications by $m$ and $n$.
We can define $d\wtil{G}$, $d (\blu\sp * G)$ and $d (\bld\sp * G\sprime)$
as global sections of the vector bundles $\Omega\sb S \sp 1\otimes \HL\sp{\otimes 2}$,
$\Omega\sb S \sp 1\otimes \blu\sp *\OOO\sb{\Pt} (3) \sp{\otimes 2}$ and 
$\Omega\sb S \sp 1\otimes   \bld\sp *\OOO\sb{\Pt} (3) \sp{\otimes 2}$, respectively.
By~\eqref{eq:sq} and~\eqref{eq:isomHL},
we get
\begin{eqnarray*}
&&
\blu\inv(\ZZ{dG}\sm \Sigma)
=\blu\inv (\ZZ{dG} \cap(\Pt\sm \cup \EED\sprime\sb j) )
=\blu\inv (\ZZ{dG})\cap U \\
&=&\ZZ{d (\blu\sp * G|U ) } 
=\ZZ{d ( \wtil{G} |U ) } 
=\ZZ{d (\bld\sp * G\sprime |U ) }\\
&=&\bld\inv (\ZZ{dG\sprime})\cap U
=\bld\inv (\ZZ{dG\sprime} \cap(\Pt\sm \cup \ED\sprime\sb j) ).
\end{eqnarray*}
Hence we get $\ZZ{dG\sprime}\cap (\Pt\sm \cup \ED\sprime\sb j)=\cremona\sb{\Sigma} (\ZZ{dG}\sm \Sigma)$.
In particular, we have $\cremona\sb{\Sigma} (\ZZ{dG}\sm \Sigma) \st \ZZ{dG\sprime}$.
\par
If $\dim \ZZ{dG\sprime} >0$, then
one of the conic curves $\ED\sprime\sb j$ is contained in $\ZZ{dG\sprime}$.
Suppose that $\ED\sprime\sb k\st \ZZ{dG\sprime}$.
Then $\ED\sb k \st \ZZ{d\wtil{G}}$ holds.
We choose affine coordinates $(x, y)$ of $\Pt$ such that $p\sb k=\blu (\ED\sb k)$ is the origin.
Let
$$
g(x, y)=\sum\sb{i+j\le 6} a\sb{ij}x^i y^j
$$
be the inhomogeneous polynomial corresponding to $G$.
Since $p\sb k=(0,0)\in \Sigma$ is contained in $\Sing D$, we have
$$
a\sb{ij}=0\quad\textrm{for}\quad i+j\le 1.
$$
Let the blowing up $\blu$ be given by 
$$
(u, v) \mapsto (x, y)=(uv, v)
$$
around a point of $\ED\sb k$.
Then $\wtil{G}$ is written 
in terms of the coordinates $(u, v)$ as
$$
\tilde g (u, v) = \beta\sp * g/v^2=\sum\sb{2\le i+j\le 6} a\sb{ij} u^i v^{i+j-2}.
$$
Since $d\wtil{G}$ is zero along the curve $\ED\sb k=\{v=0\}$ by the assumption,
we have
$$
\frac{\partial \tilde g}{\partial u} (u, 0)=a\sb {11}=0.
$$
This contradicts the fact that $p\sb k$ is a reduced point of $\ZZ{dG}$.
\end{proof}
\begin{proposition}
Let $G$ and $G\sprime$ be as above.
Then the Cremona transformation $\cremona\sb{\Sigma}$ of $\Pt$
lifts to an isomorphism
$$
{\Cremona\sb{\Sigma}}\;\;:\;\;{X\sb G}\;\;\isomto\;\;{X\sb{G\sprime}}
$$
of supersingular $K3$ surfaces.
\end{proposition}
\begin{proof}
Let $\wtil{Y}$ be the subvariety of the total space of the line bundle $\HL$ defined by
$W^2=\wtil{G}$,
where $\wtil{G}$ is the global section of $\HL\sp{\otimes 2}=L(H)$ introduced in the proof of Proposition~\ref{prop:Gsprime},
and $W$ is the fiber coordinate of $\HL$.
From~\eqref{eq:sq} and~\eqref{eq:isomHL},
we obtain  isomorphisms
$$
X\sb{G}|(\Pt\sm \cup \EED\sprime\sb j )
\;\;\cong\;\; 
\wtil{Y}|U 
\;\;\cong\;\;
X\sb{G\sprime}|(\Pt\sm \cup \ED\sprime\sb j )
$$
that are compatible with the isomorphisms~\eqref{eq:isomU}.
Since  $K3$ surfaces are minimal,
the isomorphism between  Zariski open subsets of $X\sb{G}$ and $X\sb{G\sprime}$ extends to 
an isomorphism between $X\sb{G}$ and $X\sb{G\sprime}$.
\end{proof}
\begin{remark}\label{rem:CRonNS}
We describe the action of $\Cremona\sb{\Sigma}$ on 
the numerical N\'eron-Severi lattices of the supersingular $K3$ surfaces.
We number the points of $\ZZ{dG}$ and $\ZZ{dG\sprime}$ in such a way that
\begin{eqnarray*}
\Sigma =\{ p\sb 1, \dots, p\sb 6\}, && \ZZ{dG}=\Sigma \cup \{ p\sb 7, \dots, p\sb{21}\}, \\
\Sigma\sprime =\{ p\sprime\sb 1, \dots, p\sprime\sb 6\}, &&
\ZZ{dG\sprime}=\Sigma\sprime \cup \{ p\sprime\sb 7, \dots, p\sprime\sb{21}\},
\end{eqnarray*}
where $p\sb i\sprime =\cremona\sb{\Sigma} (p\sb i)$ for $i=7, \dots, 21$.
Let $E\sb i\st X\sb{G}$ be the $(-2)$-curve 
that is contracted to $p\sb i\in \Pt$,
and  $E\sprime\sb i\st X\sb{G\sprime}$  the $(-2)$-curve  
that is contracted to $p\sprime\sb i\in \Pt$.
Then $\NS (X\sb{G})\otimes\sb{\Z}\Q$ is generated by
$[E\sb 1], \dots, [E\sb{21}], [\pol\sb{G}]$,
and 
$\NS (X\sb{G\sprime})\otimes\sb{\Z}\Q$ is generated by
$[E\sprime\sb 1], \dots, [E\sprime\sb{21}], [\pol\sb{G\sprime}]$.
Since $\cremona\sb{\Sigma} (p\sb i)=p\sprime\sb i$ for $i>6$, 
we have
$$
\Cremona\sb{\Sigma}\sp * ([E\sprime\sb i])=[E\sb i] \quad\textrm{for} \quad i>6.
$$
The exceptional curve $\EED\sb i$ on $S$ contracted to $p\sprime\sb i$ by $\bld:S\to\Pt$
is mapped by $\blu: S\to \Pt$ to the nonsingular  conic curve $\EED\sprime\sb i$
such that  $\EED\sprime\sb i \cap\ZZ{dG}=\Sigma\sm\{p\sb i\}$.
Hence
$$
\Cremona\sb{\Sigma}\sp * ([E\sprime\sb i])= 2[\pol\sb G] -\sum\sb{j=1}\sp 6 [E\sb j] + [E\sb i] \quad\textrm{for}\quad 
i=1,\dots, 6.
$$
The pull-back of a general line of $\Pt$ by $\cremona\sb{\Sigma}: \Pt\cdots\to\Pt$
is a quintic curve $Q$ such that
$Q\cap \ZZ{dG}=\Sigma$ and that the multiplicity of $Q$ at each point of $\Sigma$ is $2$ (Remark~\ref{rem:multQ}).
Thus
$$
\Cremona\sb{\Sigma}\sp * ([\pol\sb{G\sprime}])= 5[\pol\sb G] -2\sum\sb{j=1}\sp 6 [E\sb j].
$$
These formula describe the homomorphism $\Cremona\sb{\Sigma}\sp *\otimes\sb{\Z}\Q$
from $\NS(X\sb{G\sprime})\otimes\sb{\Z}\Q$ to $\NS(X\sb G)\otimes\sb{\Z}\Q$.
\end{remark}
\begin{remark}\label{rem:CRasAut}
Suppose that the point $[G\sprime]\in \moduli$ corresponding to $G\sprime\in \UUU$ in Proposition~\ref{prop:Gsprime}
coincides with the point $[G]\in \moduli$.
Then the Cremona transformation $\Cremona\sb{\Sigma}$ defines 
a right coset in $ \Aut(X\sb G)$ 
with respect to the subgroup $\Aut (X\sb G, \pol\sb G)\st  \Aut(X\sb G)$.
Indeed, the assumption $[G]=[G\sprime]$
implies the existence of a linear isomorphism
$g:\Pt\isomto \Pt$ such that $g(\ZZ{dG\sprime})=\ZZ{dG}$.
Let $\hat g\in \GL (3, k)$ be a lift of $g\in \PGL(3, k)$.
Then there exists $c\in k\sptimes$ and $H\in H\sp 0 (\Pt, \OOO\sb{\Pt} (3))$ such that
$\hat g\sp * G=c\,G\sprime +H^2$.
Let $X\sb G$ and $X\sb{G\sprime}$ be defined by $W^2=G(X,Y,Z)$ and $W\sp{\prime 2}=G\sprime(X,Y,Z)$,
respectively.
We have a lift
$\tilde{g} :X\sb{G\sprime}\isomto X\sb{G}$
of $g$ given by
$$
\tilde g\sp * W=\sqrt{c}\, W\sprime +H.
$$
The composite $\tilde{g}\circ \Cremona\sb{\Sigma}$ is an automorphism of $X\sb G$.
Since the linear isomorphism $g$ is unique up to the group
$$
\set{h\in \PGL (3, k)}{ h(\ZZ{dG})=\ZZ{dG}}=\Aut(X\sb G, \pol\sb G),
$$
the automorphism $\tilde{g}\circ \Cremona\sb{\Sigma}\in \Aut(X\sb{G})$ is also unique up to  
$\Aut(X\sb G, \pol\sb G)$.
\end{remark}
\section{The isomorphism correspondences by Cremona transformations}\label{sec:isomcor}
\subsection{The action of Cremona transformations on the moduli space}
Let $\abscode$ be a code satisfying the conditions in (ii) of Theorem~\ref{thm:codes}.
For $\gamma\in\gs\sb{\abscode}$, we denote by $G\sb{\gamma}\in \UUU$ a  homogeneous polynomial 
such that $\gamma (\Ps)=\ZZ{dG\sb{\gamma}}$.
Let $c\in \Pow (\Ps)$ be a word of weight $6$.
Recall from Definition~\ref{def:center} that
$\gamma (c)$ is a center of Cremona transformation for $G\sb{\gamma}$
if no three points of $\gamma (c)$ are collinear and 
there are no nonsingular conic curves $C\st \Pt$ such that $|C\cap \gamma (c)|\ge 5$ and 
$|C\cap \gamma (\Ps)|\ge 6$.
By Propositions~\ref{prop:linewt},~\ref{prop:conicwt} and~\ref{prop:wordsbijection},
we see that the following conditions on a word $c\in \Pow(\Ps)$ of weight $6$ is equivalent to each other:
\begin{itemize}
\item the word $c$ satisfies the following:
\begin{itemize}
\item[(i)] $|c\cap l|\le 2$ for any linear word $l$ of $\abscode$, and
\item[(ii)] $|c\cap q|\le 4$ for any quadratic word $q$ of $ \abscode$,
\end{itemize}
\item there exists $\gamma\in \gs\sb{\abscode}$ such that
$\gamma (c)$ is a center of Cremona transformation for
$G\sb{\gamma}$, and 
\item for arbitrary $\gamma\in \gs\sb{\abscode}$,
$\gamma (c)$ is a center of Cremona transformation for
$G\sb{\gamma}$.
\end{itemize}
\begin{definition}
A word $c\in \Pow(\Ps)$ of weight $6$ is called a \emph{center of Cremona transformation with respect to $\abscode$}
if  the above conditions are satisfied.
\end{definition}
Let $c$ be a center of Cremona transformation with respect to $\abscode$.
For $\gamma\in \gs\sb{\abscode}$,
we put
$$
\Sigma:=\gamma (c),
$$
and consider the Cremona transformation $\cremona\sb{\Sigma}$.
We put 
$$
Z\sprime\sb{\gamma, c}:=
\set{\cremona\sb{\Sigma} (\gamma (P))}{P\in \Ps \sm c } 
\cup 
\{ p\sprime\sb 1, \dots, p\sprime\sb 6\},
$$
where $p\sprime\sb i$ is the image of the strict transform $\EED\sb i\st S$ of
the conic curve $\EED\sprime\sb i$ that contains $\gamma (c) \sm \{p\sb i\}$. 
By Proposition~\ref{prop:Gsprime},
there exists a polynomial $G\sprime\sb{\gamma, c}\in \UUU$ such that
$Z\sprime\sb{\gamma, c}=\ZZ{dG\sprime\sb{\gamma, c}}$.
Even though the polynomial $G\sprime\sb{\gamma, c}$ is not uniquely determined,
the corresponding point $[G\sprime\sb{\gamma, c}]\in \moduli$ is uniquely determined by $c$ and $\gamma$.
The map
$\gamma\mapsto [G\sprime\sb{\gamma, c}]$ gives a morphism from $\gs\sb{\abscode}$ to $\moduli$.
It is obvious that this morphism descends to the morphism 
$$
\map{\mcremona\sb{c}}{ \PGL(3, k)\backslash \gs\sb{\abscode}}{\moduli}.
$$
\subsection{The case where Artin invariant is $2$}
Let $T$ be  $A$, $B$ or $C$,
and let $c\in \Pow (\Ps)$ be a center 
of Cremona transformation with respect to $\CT$.
The image by $\mcremona\sb{c}$ of the connected component
$$
(\PGL(3,k)\backslash \gs\sb T)\sp+ =\set{[\gamma\sb{\lambda}]}{\lambda\in \A\sp 1 \sm \{0,1,\omega, \bar\omega\}}
$$
of $\PGL(3,k)\backslash \gs\sb T$ is a connected component of 
$\moduli\sb 2=\moduli\sb{A}\sqcup\moduli\sb{B}\sqcup \moduli\sb{C}$,
and hence there exists $T\sprime\in \{A, B, C\}$ such that
$\mcremona\sb{c}$ yields a morphism
$$
\map{\mcremona\sb{T, c}\sp{+}}{(\PGL(3,k)\backslash \gs\sb T)\sp+}{\moduli\sb{T\sprime}}.
$$
Using $\mcremona\sb{T, c}\sp{+}$ and the quotient morphism
$$
\map{p\sb{T}}{(\PGL(3,k)\backslash \gs\sb T)\sp+}{\moduli\sb{T}}
$$
by  $N\sb{T}=\Ker(\Aut(\CT)\to\wtil{T})$,
we obtain an irreducible isomorphism correspondence
$$
D\sb{T, T\sprime}[c]:=\set{(p\sb{T}([\gamma]),\mcremona\sb{T, c}\sp{+}([\gamma])) }{%
[\gamma]\in (\PGL(3,k)\backslash \gs\sb T)\sp+ } \;\;\st\;\;
\moduli\sb{T}\times
\moduli\sb{T\sprime}.
$$
For $\sigma\in N\sb T$, we have
$$
\mcremona\sb{c} ([\gamma\circ\sigma])=\mcremona\sb{\sigma(c)} ([\gamma]).
$$
Hence the type  $T\sprime$ and the correspondence $D\sb{T, T\sprime}[c]$
depends only on the orbit of $c$ under the action of $N\sb T$.
We present in Table~\ref{table:NTorbits} the decomposition of the set of centers of
Cremona transformation with respect to $\CT$ into the orbits under the action of $N\sb T$.
For each orbit,
the type $T\sprime$ and the defining equation of the isomorphism correspondence 
$D\sb{T, T\sprime}[c]$ are also given.
\begin{table}
{\small
$$
\renewcommand{\arraystretch}{1.3}
\begin{array}{|c|c|c|c|c|}
\hline
T & c & | N\sb T \cdot c | & T\sprime & \phantom{aaa}D\sb{T, T\sprime}[c]\phantom{aaa} \\
\hline
A& \{1, 2, 8, 10, 15, 16\}&  12&  A &  \trivial \\
A& \{1, 7, 8, 15, 18, 21\}&  144& A &		D1=0	 \\
A& \{2, 4, 8, 11, 14, 16\}&  576& B& D2=0 \\
A& \{1, 4, 6, 9, 12, 20\}&  72& A& D3=0\\
A& \{2, 8, 10, 12, 14, 21\}&  72& A& D3=0 \\
A& \{5, 8, 9, 10, 14, 16\}&  48& A&  \trivial \\
A& \{4, 9, 12, 16, 17, 18\}&  24& A&  \trivial \\
A& \{1, 2, 9, 10, 16, 19\}&  36& A&  \trivial\\
A& \{7, 12, 13, 14, 19, 20\}&  36& A&  \trivial\\
A& \{2, 6, 9, 10, 13, 21\}&  48& C& D4=0 \\
A& \{2, 5, 11, 13, 17, 21\}&  576& A& D1=0 \\
\hline
B& \{2, 7, 8, 9, 10, 17\}&  216& B&  \trivial \\
B& \{1, 2, 11, 12, 13, 18\}&  72& B& D5=0 \\
B& \{4, 5, 6, 10, 13, 19\}&  54& B&  \trivial\\
B& \{4, 7, 12, 15, 20, 21\}&  6& B&  \trivial\\
B& \{1, 2, 6, 10, 14, 16\}&  54& B&  \trivial\\
B& \{3, 5, 14, 16, 19, 20\}&  108& B&  \trivial\\
B& \{1, 3, 8, 12, 13, 17\}&  108& B&  \trivial\\
B& \{1, 5, 6, 16, 20, 21\}&  216& A& D6=0\\
B& \{2, 6, 9, 13, 16, 18\}&  36& B& \trivial\\
B& \{3, 7, 8, 10, 19, 21\}&  216& B& D5=0\\
B& \{2, 5, 9, 16, 18, 19\}&  108& B& \trivial\\
B& \{1, 3, 5, 15, 19, 21\}&  72& B& D5=0\\
B& \{2, 6, 7, 16, 20, 21\}&  108& B& \trivial\\
\hline
C& \{3, 5, 9, 13, 17, 21\}&  960& A& D7=0\\
C& \{3, 5, 10, 14, 17, 21\}&  64& C& D8=0\\
C& \{1, 5, 8, 10, 14, 18\}&  960& C& \trivial\\
C& \{1, 2, 5, 8, 18, 19\}&  240& C& \trivial \\
\hline
\end{array}
$$
\par
\medskip
\begin{eqnarray*}
\Delta &:=& \JS+\JT, \\
D1&:=&\JT^6 \JS ^2+\JT^4 \JS ^4+\JT^2  \JS ^6+\JT^4  \JS ^3+\JT^3  \JS ^4+\JT^4  \JS ^2+\JT^3  \JS ^3+ \JT^2  \JS ^4+\\
&&+\JT^4  \JS +\JT  \JS ^4+\JT^3  \JS +\JT^2 
\JS ^2+\JT  \JS ^3+\JT^3+\JT^2  \JS +\JT  \JS ^2+\JS ^3,\\ 
D2&:=&\JS ^4+\JT^2  \JS +\JT  \JS ^2+\JT  \JS +\JT,\\
D3&:=&\JT^3  \JS ^2+\JT^2  \JS ^3+\JT^2  \JS ^2+\JT  \JS +1,\\
D4&:=&\JT^2  \JS ^4+\JT  \JS ^2+\JT  \JS +\JT+\JS ,\\
D5&:=&\JT^4  \JS +\JT^3  \JS ^2+\JT^2  \JS ^3+\JT  \JS ^4+\JT^3  \JS +\JT^2  \JS ^2+\JT  \JS ^3+1,\\
D6&:=&\JT^4+\JT^2  \JS +\JT  \JS ^2+\JT  \JS +\JS ,\\
D7&:=&\JT^4  \JS ^2+\JT^2  \JS +\JT  \JS +\JT+\JS ,\\
D8&:=&\JT^4  \JS ^4+\JT^3  \JS +\JT^2  \JS ^2+\JT  \JS ^3+\JT^3+\JT^2  \JS +\JT  \JS ^2+\JS ^3.\\
\end{eqnarray*}
}
\caption{Isomorphism correspondences  by Cremona transformations}\label{table:NTorbits}
\end{table}
\par
\medskip
We will explain the algorithm for obtaining the defining equation of $D\sb{T, T\sprime}[c]$.
For example,
consider the case where $T=A$ and $c=\{P\sb{1}, P\sb{4}, P\sb{6}, P\sb{9}, P\sb{12}, P\sb{20}\}$.
The six points $\Sigma= \gamma\sb{\lambda} (c) =\{p\sb 1, \dots, p\sb 6\}$ are as follows:
$$
\begin{array}{lccl}
p\sb{1}:=\gamma\sb{\lambda} (P\sb 1) =  [1,  \omega,0], &\quad&\quad& p\sb{4}  :=\gamma\sb{\lambda}(P\sb{9}) = [1,1,1], \\
p\sb{2}:=\gamma\sb{\lambda} (P\sb 4) =  [1+\lambda,\lambda,1], &\quad&\quad&  p\sb{5}:=\gamma\sb{\lambda}(P\sb{12}) = [0,1,0], \\  
p\sb{3}:=\gamma\sb{\lambda}(P\sb{6}) = [\lambda,1+\lambda,1], &\quad&\quad& p\sb{6}:=\gamma\sb{\lambda}(P\sb{20}) = [0,0,1].
\end{array}
$$
Solving linear equations,
we see that the $3$-dimensional linear space
$H^0 (\Pt, \ideal \sb{\Sigma}^2 (5))$ is generated by the  homogeneous quintic polynomials in Table~\ref{table:F123}.
\begin{table}
{\small
\begin{eqnarray*}
F\sb 1&:=&\phantom{+}
{X}^{4}Z+
{X}^{3}{Z}^{2}{\lambda}^{2}{\omega}+
{X}^{3}{Z}^{2}\lambda\,{\omega}+
{X}^{2}{Y}^{2}Z{\omega}+
{X}^{2}Y{Z}^{2}{\lambda}^{2}{\omega}+
{X}^{2}Y{Z}^{2}\lambda\,{\omega}+\\&&
+{X}^{2}Y{Z}^{2}{\omega}+  
{X}^{2}{Z}^{3}{\lambda}^{4}{\omega}+
{X}^{2}{Z}^{3}{\lambda}^{4}+
{X}^{2}{Z}^{3}{\lambda}^{2}+
{X}^{2}{Z}^{3}\lambda\,{\omega}+
{X}^{2}{Z}^{3}+\\&&
+X{Y}^{2}{Z }^{2}{\lambda}^{2}{\omega}+
X{Y}^{2}{Z}^{2}\lambda\,{\omega}+
X{Y}^{2}{Z}^{2}{\omega}+XY{Z}^{3}{\omega}+
{Y}^{3}{Z}^{2}{\lambda}^{2}{\omega}+\\&&
+{Y}^{3}{Z}^{2}\lambda\,{\omega}+
{Y}^{2}{Z}^{3}{\lambda}^{4}{\omega}+
{Y}^{2}{Z}^{3}{\lambda}^{4}+
{Y}^{2}{Z}^{3}{\lambda}^{2}+
{Y}^{2}{Z}^{3}\lambda\,{\omega},\\
\\
F\sb 2&:=&\phantom{+}
X{Y}^{3}Z{\lambda}^{2}{\omega}+
{X}^{3}{Z}^{2}{\lambda}^{2}+
{X}^{2}Y{Z}^{2}\lambda+
{X}^{2}{Z}^{3}\lambda+
X{Y}^{2}{Z}^{2}{\lambda}^{2}+
XY{Z}^{3}+\\&&
+{X}^{2}{Y}^{2}Z+ 
{X}^{2}{Z}^{3}{\lambda}^{4}+
{Y}^{3}{Z}^{2}{\lambda}^{4}+
{Y}^{2}{Z}^{3}{\lambda}^{4}+
{Y}^{2}{Z}^{3}\lambda+
{X}^{2}Y{Z}^{2}{\lambda}^{4}+\\&&
+{X}^{3}{Z}^{2}\lambda+
X{Y}^{2}{Z}^{2}\lambda+
{X}^{4}Y+
X{Y}^{2}{Z}^{2}+
{Y}^{3}{Z}^{2}\lambda+\
{Y}^{2}{Z}^{3}{\lambda}^{3}{\omega}+\\&&
+{X}^{2}{Z}^{3}{\lambda}^{5}{\omega}+
{X}^{2}{Z}^{3}{\lambda}^{3}{\omega}+
{Y}^{2}{Z}^{3}{\lambda}^{5}{\omega}+
{Y}^{2}{Z}^{3}{\lambda}^{2}{\omega}+
{Y}^{2}{Z}^{3}\lambda\,{\omega}+
{Y}^{2}{Z}^{3}{\lambda}^{6}{\omega}+\\&&
+{X}^{2}{Y}^{3}{\omega}+
{X}^{3}YZ{\lambda}^{2}{\omega}+
{X}^{2}{Z}^{3}{\lambda}^{4}{\omega}+
{X}^{2}Y{Z}^{2}\lambda\,{\omega}+
{Y}^{3}{Z}^{2}\lambda\,{\omega}+
X{Y}^{2}{Z}^{2}{\omega}+\\&&
+{X}^{2}{Z}^{3}{\lambda} ^{6}{\omega}+
X{Y}^{3}Z{\omega}+
{Y}^{3}{Z}^{2}{\lambda}^{4}{\omega}+
X{Y }^{3}Z\lambda\,{\omega}+
{Y}^{2}{Z}^{3}{\lambda}^{4}{\omega}+
{X}^{3}YZ\lambda\,{\omega}+\\&&
+{X}^{2}{Y}^{2}Z\lambda\,{\omega}+
{X}^{2}{Y}^{2}Z{\omega}+
{X}^{2}{Y}^{2}Z{\lambda}^{2}{\omega}+
{X}^{2}Y{Z}^{2}{\lambda}^{4}{\omega}, \\
\\
F\sb 3&:=&\phantom{+}
X{Y}^{3}Z{\lambda}^{2}{\omega}+
{X}^{2}Y{Z}^{2}{\lambda}^{2}{\omega}+
{X}^{3}{Y}^{2}{\omega}+
{X}^{3}{Z}^{2}+
{X}^{2}Y{Z}^{2}\lambda+
{X}^{2}{Z}^{3}\lambda+\\&&
+XY{Z}^{3}+
{X}^{2}Y{Z}^{2}+
{X}^{2}{Y}^{2}Z+
{X}^{2}{Z}^{3}{\lambda}^{4}+
{Y}^{3}{Z}^{2}{\lambda}^{2}+
{Y}^{2}{Z}^{3}{\lambda}^{4}+ 
{Y}^{2}{Z}^{3}\lambda+\\&&
+X{Y}^{2}{Z}^{2}{\lambda}^{4}+
{X}^{3}{Z}^{2}{\lambda}^{4}+
{X}^{3}{Z}^{2}\lambda+
X{Y}^{2}{Z}^{2}\lambda+
{X}^{2}Y{Z}^{2}{\lambda}^{2}+
X{Y}^{2}{Z}^{2}+\\&&
+{Y}^{3}{Z}^{2}\lambda+ 
{Y}^{2}{Z}^{3}{\lambda}^{3}{\omega}+
{X}^{5}+
{X}^{2}{Z}^{3}{\lambda}^{5}{\omega}+
{X}^{2}{Z}^{3}{\lambda}^{3}{\omega}+
{Y}^{2}{Z}^{3}{\lambda}^{5}{\omega}+\\&&
+{X}^{2}{Z}^{3}{\lambda}^{2}{\omega}+ 
{Y}^{2}{Z}^{3}\lambda\,{\omega}+
{X}^{3}{Z}^{2}{\lambda}^{4}{\omega}+
{Y}^{2}{Z}^{3}{\lambda}^{6}{\omega}+
{X}^{3}YZ{\omega}+
{X}^{3}YZ{\lambda}^{2}{\omega}+ \\&&
+{X}^{2}Y{Z}^{2}\lambda\,{\omega}+ 
{Y}^{3}{Z}^{2}\lambda\,{\omega}+
{Y}^{3}{Z}^{2}{\lambda}^{2}{\omega}+
X{Y}^{2}{Z}^{2}{\omega}+XY{Z}^{3}{\omega}+
{X}^{2}{Z}^{3}{\lambda}^{6}{\omega}+ \\&&
+X{Y}^{3}Z\lambda\,{\omega}+ 
X{Y}^{2}{Z}^{2}{\lambda}^{4}{\omega}+
{X}^{3}YZ\lambda\,{\omega}+
{X}^{3}{Z}^{2}{\lambda}^{2}{\omega}+
{X}^{2}{Y}^{2}Z\lambda\,{\omega}+ \\&&
+X{Y}^{2}{Z}^{2}{\lambda}^{2}{\omega}+
{X}^{2}{Y}^{2}Z{\lambda}^{2}{\omega}.
\end{eqnarray*}
}
\caption{Generators of $H^0 (\Pt, \ideal \sb{\Sigma}^2 (5))$}\label{table:F123}
\end{table}
The Cremona transformation $\cremona\sb{\Sigma}: \Pt\cdots\to\Pt$ is given by
$$
[X,Y,Z] \mapsto [F\sb 1, F\sb 2, F\sb 3].
$$
The points $\gamma\sb{\lambda} (P\sb i)$ $(P\sb i \notin c)$ is mapped by $\cremona\sb{\Sigma}$ to the  points
in Table~\ref{table:qs1}.
\begin{table}
\newcommand{\qi}[3]{q\sb{#1}&:=&\cremona\sb{\Sigma}(\gamma\sb{\lambda}(P\sb{#2}))\;\;=\;\;#3}
\begin{eqnarray*}
\qi{1}{2}{[\;0,1,\omega\;],}\\
\qi{2}{3}{[\;\omega, \lambda ^2+\bar\omega \lambda + 1, \lambda (\lambda +\bar\omega)\;],}\\
\qi{3}{5}{[\;\omega, \lambda ^2+\bar\omega \lambda + \bar\omega, (\lambda +1)(\lambda +\omega)\;],}\\
\qi{4}{7}{[\;\omega, \lambda ^2+\omega \lambda + \omega, (\lambda +\bar\omega)^2\;],}\\
\qi{5}{8}{[\;\omega, (\lambda +1)^2, \lambda ^2+\omega\lambda  +1\;],}\\
\qi{6}{10}{[\;\omega, (\lambda +\omega)^2, (\lambda +1)(\lambda +\bar\omega)\;],}\\
\qi{7}{11}{[\;\omega, (\lambda +\omega)(\lambda +\bar\omega), \lambda ^2+\lambda +\omega\;],}\\
\qi{8}{13}{[\;0,1,1\;],}\\
\qi{9}{14}{[\;\omega, \lambda ^2+\omega \lambda +1, (\lambda +\omega)^2\;],}\\
\qi{10}{15}{[\;\omega, \lambda ^2, \lambda ^2+\omega\lambda +\omega\;],}\\
\qi{11}{16}{[\;\omega, (\lambda +\bar\omega)^2,\lambda (\lambda +\omega)\;],}\\
\qi{12}{17}{[\;\omega,\lambda  (\lambda +\omega), \lambda ^2\;],}\\
\qi{13}{18}{[\;0,0,1\;],}\\
\qi{14}{19}{[\;\omega, \lambda ^2+\lambda +\bar\omega,\lambda  (\lambda +1)\;],}\\
\qi{15}{21}{[\;\omega, (\lambda +1)(\lambda +\bar\omega), (\lambda +1)^2\;].}
\end{eqnarray*}
\newcommand{\qqi}[3]{q\sb{#1}&:=&\bld (\EED\sb{#2})\;\;=\;\;#3}
\begin{eqnarray*}
\qqi{16}{1}{[\;0,1,0\;],}\\
\qqi{17}{2}{[\;\omega, (\lambda+\bar\omega)\lambda, \lambda^2+\bar\omega\lambda+\bar\omega\;],}\\
\qqi{18}{3}{[\;\omega, (\lambda+1)(\lambda+\omega), \lambda^2+\bar\omega\lambda+1\;],}\\
\qqi{19}{4}{[\;\omega, \lambda^2+\lambda+\omega, \lambda^2+\lambda+\bar\omega\;],}\\
\qqi{20}{5}{[\;0,1,\bar\omega \;],}\\
\qqi{21}{6}{[\;\omega, \lambda(\lambda+1), (\lambda+\omega)(\lambda+\bar\omega)\;].}
\end{eqnarray*}
\caption{Points $q\sb i$}\label{table:qs1}
\end{table}
The conic curve $\EED\sb 1\sprime\st\Pt$ containing $\Sigma\sm \{p\sb 1\}$
is defined by
$$
E\sb 1:= X^2+(\lambda^2+\lambda) YZ +(\lambda^2+\lambda +1) ZX=0.
$$
Let $V\sb 1$ be the vector space of cubic homogeneous polynomials $C$ 
such that $E\sb 1 C$ is a member of $ H^0 (\Pt, \ideal \sb{\Sigma}^2 (5))$.
Then we have $\dim V\sb 1=2$, and the image of the linear map
$V\sb 1\to H^0 (\Pt, \ideal \sb{\Sigma}^2 (5))$
given by $C\mapsto E\sb 1 C$ is spanned by
$F\sb 1$ and $F\sb 3$.
Hence the image $\bld (\EED\sb 1)$ of the strict transform $\EED\sb 1\st S$ of $\EED\sb 1\sprime$ 
is $[0,1,0]$.
In the same way, we calculate  
$\bld (\EED\sb i)$   as in Table~\ref{table:qs1}.
\begin{table}
{\small
\begin{eqnarray*}
&LW=&
\{ \{1, 3, 5, 11, 17\}, 
\{4, 5, 6, 8, 12\}, 
\{1, 2, 6, 10, 18\}, 
\{2, 3, 8, 19, 21\},\\
&& \{1, 4, 9, 14, 21\}, 
\{8, 9, 10, 11, 15\}, 
\{6, 7, 11, 20, 21\}, 
\{1, 7, 12, 15, 19\}, \\
&&\{3, 9, 12, 18, 20\}, 
\{1, 8, 13, 16, 20\}, 
\{7, 8, 14, 17, 18\},
\{5, 10, 14, 19, 20\},\\ 
&&\{2, 4, 15, 17, 20\}\}.
\end{eqnarray*}
\begin{eqnarray*}
&QW=&
\{\{2, 3, 4, 5, 9, 10, 13, 16\}, 
\{2, 3, 4, 7, 10, 11, 12, 14\}, 
\{2, 3, 5, 6, 7, 9, 14, 15\}, \\&&
\{2, 4, 5, 7, 9, 11, 18, 19\}, 
\{3, 4, 5, 7, 10, 15, 18, 21\},
\{3, 4, 6, 7, 9, 10, 17, 19\}, \\&&
\{2, 5, 6, 7, 13, 16, 17, 19\}, 
\{2, 3, 6, 11, 12, 13, 15, 16\}, 
\{3, 4, 6, 11, 14, 15, 18, 19\},\\&& 
\{3, 5, 6, 13, 14, 16, 18, 21\}, 
\{2, 5, 7, 9, 10, 12, 17, 21\}, 
\{4, 6, 7, 10, 13, 14, 15, 16\},\\&& 
\{3, 4, 7, 12, 13, 16, 17, 21\}, 
\{5, 7, 9, 11, 12, 13, 14, 16\}, 
\{2, 6, 9, 11, 12, 14, 17, 19\}, \\&&
\{4, 6, 9, 11, 13, 16, 17, 18\}, 
\{2, 7, 9, 13, 15, 16, 18, 21\}, 
\{5, 6, 9, 15, 17, 18, 19, 21\}, \\&&
\{3, 7, 10, 11, 13, 16, 18, 19\}, 
\{6, 9, 10, 12, 13, 16, 19, 21\}, 
\{3, 6, 10, 12, 14, 15, 17, 21\}, \\&&
\{2, 5, 11, 12, 14, 15, 18, 21\}, 
\{4, 10, 11, 12, 17, 18, 19, 21\}, 
\{2, 10, 11, 13, 14, 16, 17, 21\},\\&& 
\{4, 5, 11, 13, 15, 16, 19, 21\}, 
\{2, 4, 12, 13, 14, 16, 18, 19\}, 
\{5, 10, 12, 13, 15, 16, 17, 18\}, \\&& 
\{3, 9, 13, 14, 15, 16, 17, 19\}\}.
\end{eqnarray*}
}
\caption{Sets $LW$ and $QW$}\label{table:LWQW}
\end{table}
The set $LW$ of collinear $5$-tuples of the points in 
$Z\sprime=\{q\sb 1, \dots, q\sb{21}\}$
and the set $QW$ of  $8$-tuples of the points in 
$Z\sprime$ that are on a nonsingular conic curve are  given in Table~\ref{table:LWQW},
where $\{1, 3, 5, 11, 17\}$ means $\{q\sb{1}, q\sb{3},q\sb{5},q\sb{11},q\sb{17}\}$, for example.
Since $|LW|=13$ and $|QW|=28$,
we see that the type $T\sprime$ of the target moduli curve  is $A$.
Let  $\sigma$ be the following permutation:
\begin{multline*}
\left(
\begin{array}{ccccccccccc}
1&2&3&4&5&6&7&8&9&10&11\\
13& 16& 2& 6& 21& 4& 18& 17& 3& 7& 12\\
\end{array}
\right. \\
\left.
\begin{array}{cccccccccccccccccccc}
&&&&&&&&&&12&13&14&15&16&17&18&19&20&21\\
&&&&&&&&&&1& 8& 9&  11& 15& 14& 20& 5& 19& 10
\end{array}
\right).
\end{multline*}
Then the map
$$
\map{\gamma\sprime}{\Ps}{\Pt}
$$
defined by $\gamma\sprime(P\sb i)=q\sb{\sigma (i)}$
yields 
bijections from the set of linear words in $\CA$ (see Table~\ref{table:linesA}) to $LW$ and 
from the set of quadratic words in $\CA$ (see Tables~\ref{table:qsprimeA} and~\ref{table:qabcA}) to $QW$.
Hence the map $\gamma\sprime$ is an element of $\gs\sb{A}$.
We make the linear change of homogeneous coordinates of $\Pt$ so that
\begin{eqnarray*}
&& \gamma\sprime(P\sb{18})=q\sb{20}=[1,0,0], \\
&&\gamma\sprime (P\sb{12})=q\sb{1}=[0,1,0],\;\;\;\;
\gamma\sprime(P\sb{13})=q\sb{8}=[1,1,0],\;\\
&&\gamma\sprime(P\sb{20})=q\sb{19}=[0,0,1],\;\;\;
\gamma\sprime(P\sb{19})=q\sb{5}=[1,0,1],
\end{eqnarray*}
hold (see~\eqref{eq:gammaPFIX}); that is, we multiply the matrix
$$
\left[
\begin{matrix}
\lambda^2+\lambda & 1 &\bar\omega \\
\lambda^2+\lambda+1 &\bar\omega & 1\\
\omega\lambda+1 & 0 & 0
\end{matrix}
\right]
$$
from the left to the vectors $\gamma\sprime (P\sb i)=q\sb{\sigma(i)}$.
Then we have
$$
\gamma\sprime (P\sb 1)=q\sb{13}=[1, \omega, 0].
$$
Therefore the projective equivalence class $[\gamma\sprime]\in \PGL(3, k)\backslash \gs\sb{A}$
of $\gamma\sprime$ 
is contained in the connected component $(\PGL(3, k)\backslash \gs\sb{A})\sp +$,
because otherwise we would have $\gamma\sprime (P\sb 1)=[1, \bar\omega, 0]$.
Since
$$
\gamma\sprime (P\sb {10})=q\sb{7}=[0, 1, \lambda+\bar\omega],
$$
the point $[\gamma\sprime]$ corresponds to
$1/(\lambda+\bar\omega)$ under the isomorphism 
$(\PGL(3, k)\backslash \gs\sb{A})\sp +\cong\A\sp 1\sm\{0,1,\omega,\bar\omega\}$.
Substituting $1/(\lambda+\bar\omega)$ for $\lambda$ in
$$
J\sb A=\frac{(\lambda^2+\lambda+1)^3}{\lambda^2 (\lambda+1)^2},
$$
we see that 
the $J\sb A$-invariant of  $[\gamma\sprime]$ is  equal to
$$
J\sprime\sb A=\frac{\lambda^3 (\lambda+1)^3}{(\lambda^2+\lambda+1)^2}.
$$
Eliminating $\lambda$ from $J\sb A$ and $J\sb A\sprime$,
we obtain the defining equation
$$
1+J\sb A J\sb A\sprime + J\sb A ^2 J\sb A\sp{\prime 2} + J\sb A ^3 J\sb A\sp{\prime 2} + J\sb A ^2 J\sb A\sp{\prime 3}=0
$$
of the isomorphism correspondence given by the Cremona transformation with the center 
$c=\{P\sb{1}, P\sb{4}, P\sb{6}, P\sb{9},
P\sb{12}, P\sb{20}\}$.
\begin{table}
\begin{center}
\hbox{
\hskip -3.5cm
\vtop{
{\tiny
\begin{eqnarray*}
D_{A,A,1} * D_{A,A,1} &=&\Delta\sb{A}+ D_{A,A,2}, \\ 
D_{A,A,1} * D_{A,A,2} &=&D_{A,A,1}+ D_{A,A,2}, \\ 
D_{A,A,2} * D_{A,A,1} &=&D_{A,A,1}+ D_{A,A,2}, \\ 
D_{A,A,2} * D_{A,A,2} &=&\Delta\sb{A}+ D_{A,A,1}+ D_{A,A,2}, \\ 
D_{A,A,1} * D_{A,B,1} &=&D_{A,B,2}, \\ 
D_{A,A,1} * D_{A,B,2} &=&D_{A,B,1}+ D_{A,B,2}, \\ 
D_{A,A,2} * D_{A,B,1} &=&D_{A,B,1}+ D_{A,B,2}, \\ 
D_{A,A,2} * D_{A,B,2} &=&D_{A,B,1}+ D_{A,B,2}, \\ 
D_{A,A,1} * D_{A,C,1} &=&D_{A,C,2}, \\ 
D_{A,A,1} * D_{A,C,2} &=&D_{A,C,1}+ D_{A,C,2}, \\ 
D_{A,A,2} * D_{A,C,1} &=&D_{A,C,1}+ D_{A,C,2}, \\ 
D_{A,A,2} * D_{A,C,2} &=&D_{A,C,1}+ D_{A,C,2}, \\ 
D_{B,B,1} * D_{B,B,1} &=&\Delta\sb{B}+ D_{B,B,1}, \\ 
D_{B,B,1} * D_{B,A,1} &=&D_{B,A,1}+ D_{B,A,2}, \\ 
D_{B,B,1} * D_{B,A,2} &=&D_{B,A,1}+ D_{B,A,2}, \\ 
D_{B,B,1} * D_{B,C,1} &=&D_{B,C,2}, \\ 
D_{B,B,1} * D_{B,C,2} &=&D_{B,C,1}+ D_{B,C,2}, \\ 
D_{C,C,1} * D_{C,C,1} &=&\Delta\sb{C}+ D_{C,C,1}, \\ 
D_{C,C,1} * D_{C,B,1} &=&D_{C,B,2}, \\ 
D_{C,C,1} * D_{C,B,2} &=&D_{C,B,1}+ D_{C,B,2}, \\ 
D_{C,C,1} * D_{C,A,1} &=&D_{C,A,1}+ D_{C,A,2}, \\ 
D_{C,C,1} * D_{C,A,2} &=&D_{C,A,1}+ D_{C,A,2}, \\ 
D_{A,B,1} * D_{B,B,1} &=&D_{A,B,1}+ D_{A,B,2}, \\ 
D_{A,B,2} * D_{B,B,1} &=&D_{A,B,1}+ D_{A,B,2}, \\ 
D_{A,B,1} * D_{B,A,1} &=&\Delta\sb{A}+ D_{A,A,2}, \\ 
D_{A,B,1} * D_{B,A,2} &=&D_{A,A,1}+ D_{A,A,2}, \\ 
D_{A,B,2} * D_{B,A,1} &=&D_{A,A,1}+ D_{A,A,2}, \\ 
D_{A,B,2} * D_{B,A,2} &=&\Delta\sb{A}+ D_{A,A,1}+ D_{A,A,2}, \\ 
D_{A,B,1} * D_{B,C,1} &=&D_{A,C,1}, \\ 
D_{A,B,1} * D_{B,C,2} &=&D_{A,C,1}+ D_{A,C,2}, \\ 
D_{A,B,2} * D_{B,C,1} &=&D_{A,C,2}, \\ 
D_{A,B,2} * D_{B,C,2} &=&D_{A,C,1}+ D_{A,C,2}, \\ 
D_{B,A,1} * D_{A,A,1} &=&D_{B,A,2}, \\ 
D_{B,A,1} * D_{A,A,2} &=&D_{B,A,1}+ D_{B,A,2}, \\ 
D_{B,A,2} * D_{A,A,1} &=&D_{B,A,1}+ D_{B,A,2}, \\ 
D_{B,A,2} * D_{A,A,2} &=&D_{B,A,1}+ D_{B,A,2}, \\ 
D_{B,A,1} * D_{A,B,1} &=&\Delta\sb{B}+ D_{B,B,1}, \\ 
D_{B,A,1} * D_{A,B,2} &=&D_{B,B,1}, \\ 
D_{B,A,2} * D_{A,B,1} &=&D_{B,B,1}, \\ 
D_{B,A,2} * D_{A,B,2} &=&\Delta\sb{B}+ D_{B,B,1}, \\ 
D_{B,A,1} * D_{A,C,1} &=&D_{B,C,1}+ D_{B,C,2}, \\ 
D_{B,A,1} * D_{A,C,2} &=&D_{B,C,2}, \\ 
D_{B,A,2} * D_{A,C,1} &=&D_{B,C,2}, 
\end{eqnarray*}
}
}
\hskip -6cm
\vtop{
{\tiny
\begin{eqnarray*}
D_{B,A,2} * D_{A,C,2} &=&D_{B,C,1}+ D_{B,C,2}, \\ 
D_{B,C,1} * D_{C,C,1} &=&D_{B,C,2}, \\ 
D_{B,C,2} * D_{C,C,1} &=&D_{B,C,1}+ D_{B,C,2}, \\ 
D_{B,C,1} * D_{C,B,1} &=&\Delta\sb{B}, \\ 
D_{B,C,1} * D_{C,B,2} &=&D_{B,B,1}, \\ 
D_{B,C,2} * D_{C,B,1} &=&D_{B,B,1}, \\ 
D_{B,C,2} * D_{C,B,2} &=&\Delta\sb{B}+ D_{B,B,1}, \\ 
D_{B,C,1} * D_{C,A,1} &=&D_{B,A,1}, \\ 
D_{B,C,1} * D_{C,A,2} &=&D_{B,A,2}, \\ 
D_{B,C,2} * D_{C,A,1} &=&D_{B,A,1}+ D_{B,A,2}, \\ 
D_{B,C,2} * D_{C,A,2} &=&D_{B,A,1}+ D_{B,A,2}, \\ 
D_{C,B,1} * D_{B,B,1} &=&D_{C,B,2}, \\ 
D_{C,B,2} * D_{B,B,1} &=&D_{C,B,1}+ D_{C,B,2}, \\ 
D_{C,B,1} * D_{B,A,1} &=&D_{C,A,1}, \\ 
D_{C,B,1} * D_{B,A,2} &=&D_{C,A,2}, \\ 
D_{C,B,2} * D_{B,A,1} &=&D_{C,A,1}+ D_{C,A,2}, \\ 
D_{C,B,2} * D_{B,A,2} &=&D_{C,A,1}+ D_{C,A,2}, \\ 
D_{C,B,1} * D_{B,C,1} &=&\Delta\sb{C}, \\ 
D_{C,B,1} * D_{B,C,2} &=&D_{C,C,1}, \\ 
D_{C,B,2} * D_{B,C,1} &=&D_{C,C,1}, \\ 
D_{C,B,2} * D_{B,C,2} &=&\Delta\sb{C}+ D_{C,C,1}, \\ 
D_{C,A,1} * D_{A,A,1} &=&D_{C,A,2}, \\ 
D_{C,A,1} * D_{A,A,2} &=&D_{C,A,1}+ D_{C,A,2}, \\ 
D_{C,A,2} * D_{A,A,1} &=&D_{C,A,1}+ D_{C,A,2}, \\ 
D_{C,A,2} * D_{A,A,2} &=&D_{C,A,1}+ D_{C,A,2}, \\ 
D_{C,A,1} * D_{A,B,1} &=&D_{C,B,1}+ D_{C,B,2}, \\ 
D_{C,A,1} * D_{A,B,2} &=&D_{C,B,2}, \\ 
D_{C,A,2} * D_{A,B,1} &=&D_{C,B,2}, \\ 
D_{C,A,2} * D_{A,B,2} &=&D_{C,B,1}+ D_{C,B,2}, \\ 
D_{C,A,1} * D_{A,C,1} &=&\Delta\sb{C}+ D_{C,C,1}, \\ 
D_{C,A,1} * D_{A,C,2} &=&D_{C,C,1}, \\ 
D_{C,A,2} * D_{A,C,1} &=&D_{C,C,1}, \\ 
D_{C,A,2} * D_{A,C,2} &=&\Delta\sb{C}+ D_{C,C,1}, \\ 
D_{A,C,1} * D_{C,C,1} &=&D_{A,C,1}+ D_{A,C,2}, \\ 
D_{A,C,2} * D_{C,C,1} &=&D_{A,C,1}+ D_{A,C,2}, \\ 
D_{A,C,1} * D_{C,B,1} &=&D_{A,B,1}, \\ 
D_{A,C,1} * D_{C,B,2} &=&D_{A,B,1}+ D_{A,B,2}, \\ 
D_{A,C,2} * D_{C,B,1} &=&D_{A,B,2}, \\ 
D_{A,C,2} * D_{C,B,2} &=&D_{A,B,1}+ D_{A,B,2}, \\ 
D_{A,C,1} * D_{C,A,1} &=&\Delta\sb{A}+ D_{A,A,2}, \\ 
D_{A,C,1} * D_{C,A,2} &=&D_{A,A,1}+ D_{A,A,2}, \\ 
D_{A,C,2} * D_{C,A,1} &=&D_{A,A,1}+ D_{A,A,2}, \\ 
D_{A,C,2} * D_{C,A,2} &=&\Delta\sb{A}+ D_{A,A,1}+ D_{A,A,2}.
\end{eqnarray*}
}
}
}
\end{center}
\caption{Relations between non-trivial isomorphism correspondences}\label{table:correls}
\end{table}
\par
\medskip
Putting 
\begin{eqnarray*}
&&
D\sb{A,A, 1}:=\{ D3=0\},
\quad
D\sb{A,A, 2}:=\{ D1=0\},
\quad
\\
&&
D\sb{B,B,1}:=\{ D5=0\},
\quad
D\sb{C,C,1}:=\{ D8=0\},
\quad
\\
&&
D\sb{A, B, 1}:=\{ D2=0\}=\spT D\sb{B,A, 1}=\spT\{D6=0\}, \\
&&
D\sb{A,C,1}:=\{ D4=0\}=\spT D\sb{C,A, 1}=\spT\{D7=0\},
\end{eqnarray*}
we obtain Theorem~\ref{thm:correspondence}.
The composite $D\sb 1 * D\sb 2$ of correspondences 
\begin{eqnarray*}
D\sb 1 &=&\{f\sb 1(J\sb{T}, J\sb{T\sprime})=0\}\;\;\st\;\;\moduli\sb{T}\times \moduli\sb{T\sprime}\quand \\
D\sb 2 &=&\{f\sb 2(J\sb{T\sprime}, J\sb{T\spprime})=0\}\;\;\st\;\;\moduli\sb{T\sprime}\times \moduli\sb{T\spprime}
\end{eqnarray*}
is obtained by eliminating the variable $J\sb{T\sprime}$ from
$f\sb 1(J\sb{T}, J\sb{T\sprime})=f\sb 2(J\sb{T\sprime}, J\sb{T\spprime})=0$.
Starting from the eight isomorphism correspondences above and making composites,
we obtain irreducible isomorphism correspondences listed in Table~\ref{table:coreqs},
which have the relations in Table~\ref{table:correls}. 
This table  also shows that the isomorphism correspondences
$\Delta\sb A$, $\Delta\sb B$, $\Delta\sb C$ and the ones in Table~\ref{table:coreqs}
are closed under compositions of correspondences.

\bibliographystyle{amsplain}

\providecommand{\bysame}{\leavevmode\hbox to3em{\hrulefill}\thinspace}

\end{document}